\documentclass[12pt]{article}

\usepackage{amssymb}
\usepackage[leqno]{amsmath}
\usepackage{amsthm}

\usepackage{graphicx,color}

\usepackage{subfigure}

\newcommand\W{{\rm W}}
\newcommand\MW{{\rm {MW}}}
\newcommand\Def{{\overset {\rm {def}}{\ =\ }}}

\newcommand\out{{\rm {out}}}
\newcommand\inn{{\rm {inn}}}
\newcommand\per{{\rm {per}}}

\newcommand\const{{\rm {const}}}

\newcommand\eff{{\rm {eff}}}
\newcommand\for{{\rm {for}\;}}
\newcommand\corr{{\rm {corr}}}
\newcommand\w{{\rm w}}

\newcommand\bx{{\bar {\mathstrut x}}}

\newcommand\bxi{\bar {\mathstrut \xi}}

\newcommand\blangle{{\boldsymbol \langle }}
\newcommand\brangle{{\boldsymbol \rangle }}
\newcommand\bv{{\boldsymbol | }}

\newcommand\sign{\operatorname{sign}}
\newcommand\Spec{\operatorname{Spec}}
\newcommand\supp{\operatorname{supp}}

\newcommand\Tr{\operatorname{Tr}}

\newcommand\bR{{\mathbb R}}

\newcommand\bH{{\mathbb H}}
\newcommand\bK{{\mathbb K}}

\newcommand\bZ{{\mathbb Z}}
\newcommand\cA{{\mathcal A}}
\newcommand\cB{{\mathcal B}}

\newcommand\cE{{\mathcal E}}

\newcommand\cL{{\mathcal L}}

\newcommand\cR{{\mathcal R}}

\newcommand\cT{{\mathcal T}}

\newcommand\cZ{{\mathcal Z}}

\newtheorem{theorem}{Theorem}[section]

\newtheorem{corollary}[theorem]{Corollary}
\newtheorem{proposition}[theorem]{Proposition}

\theoremstyle{definition}

\newtheorem{remark}[theorem]{Remark}

\newenvironment{claim}[1][(\theequation)]{\refstepcounter{equation}\smallskip
\begin{trivlist}
\item[{\hskip\labelsep#1}]}{\smallskip\end{trivlist}}

\usepackage[bookmarks,pdfnewwindow]{hyperref}

\newcommand{\sect}[1]{\setcounter{equation}{0}\section{#1}}

\renewcommand{\theequation}{\thesection.\arabic{equation}}

\begin{document}

\title{%
Sharp Spectral Asymptotics for two-dimensional Schr\"odinger operator with a strong degenerating magnetic field.}

\author{%
Victor Ivrii
\footnote{Work was partially supported by NSERC grant OGP0138277.}
}
\maketitle

{\abstract
This paper is a continuation of \cite{IRO1} and \cite{Ivr1, IRO2,IRO3, IRO4, IRO5}. I consider two-dimensional Schr\"odinger operator with  degenerating magnetic field and in the generic situation I derive spectral asymptotics 
as $h\to +0$ and $\mu\to +\infty$ where $h$ and $\mu$ are Planck and coupling parameters respectively. The remainder estimate is $O(\mu^{-{\frac 1 2}}h^{-1})$
which is between $O(\mu^{-1}h^{-1})$ valid as magnetic field non-degenerates and $O(h^{-1})$ valid as magnetic field is identically 0. As $\mu$ is close to its maximal reasonable value $O(h^{-2})$ the principal part contains correction terms associated with short periodic trajectories of the corresponding classical dynamics.
\endabstract}

%

\sect{Introduction}

\setcounter{subsection}{-1}

\subsection{Preface}
We consider spectral asymptotics of
\begin{equation}
A= {\frac 1 2}\Bigl(\sum_{j,k}P_jg^{jk}(x)P_k -V\Bigr),\qquad P_j=D_j-\mu V_j
\label{0-1}
\end{equation}
where $g^{jk}$, $V_j$, $V$ are smooth real-valued functions of $x\in \bR^2$ and
$(g^{jk})$ is positive-definite matrix, $0<h\ll 1$ is a Planck parameter and 
$\mu \gg1$ is a coupling parameter. We assume that $A$ is a self-adjoint operator.

In contrast to my recent papers \cite{IRO3, IRO4, IRO5} I assume that all the coefficients are very smooth and in contrast to \cite{IRO4} I consider only two-dimensional case here. However degeneration of the magnetic field makes result much more interesting and difficult than I expected when I started this work.

I hope to investigate similar even-dimensional case in the future; odd-dimensional results are much easier and of no special interest because only in the even-dimensional full-rank case magnetic field can improve the remainder estimate.

\subsection{Assumptions and notations}
Let $g=\det (g^{jk})^{-1}$, $F_{12}=\partial_{x_1}V_2-\partial_{x_2}V_1$ and
$f=F_{12}g^{-{\frac 1 2}}V^{-1}$. Note that both $F_{12}g^{-{\frac 1 2}}$ and $f$ are coordinate independent. In \cite{IRO3} and in earlier papers I assumed that
$F$ disjoint from 0 but I will not assume this anymore; however I assume that
degenerations are generic, namely
\begin{align}
&V\ge \epsilon_0,\label{0-2}\\
&|F_{12}|+|\nabla F_{12}|\ge \epsilon_0,\label{0-3}
\end{align}
(condition (\ref{0-2}) will be dropped in the very end by rescaling arguments).
Then $\Sigma =\{x: F_{12}(x)=0\}$ is a smooth curve. One can introduce a coordinate $x_1$ as a Riemannian distance from $x$ to $\Sigma$; then $g^{11}=1$. One can always change $x_2$ to keep $g^{12}=0$ and then
\begin{equation}
A= {\frac 1 2}\Bigl(P_1^2 +P_2g^{-1}P_2-V\Bigr),\qquad P_j=D_j-\mu V_j.
\label{0-4}
\end{equation}
Furthermore, without any loss of the generality one can assume locally that \begin{equation}
V_1=0, \qquad V_2\asymp |x_1|^\nu
\label{0-5}
\end{equation}
with $\nu =2$; one can achieve this by the gauge transformation. Actually instead I consider operator (\ref{0-4}) with positive integer $\nu\ge 2$ and in the classical dynamics I even take any real $\nu\ge2$.

What I am interested is an asymptotics as $h\to +0$, $\mu\to +\infty$ of
\begin{equation}
\int e(x,x,0)\psi (x)\,dx
\label{0-6}
\end{equation}
where $e(x,y,\tau)$ is the Schwartz kernel of the spectral projector $E(\tau)$ of operator $A$. I will prove later that (\ref{0-6}) will be $O(\mu^{-s})$ as
$\mu \ge C_0h^{-\nu}$ and therefore I will assume that
\begin{equation}
\mu \le C_0h^{-\nu}.
\label{0-7}
\end{equation}
The natural answer coming from the non-degenerate case is
$\int \cE^\MW(x,0)\psi(x)\, dx$ where
\begin{equation}
 \cE^\MW={\frac 1 {2\pi}} \sum_{n\ge 0}
\theta \bigl(2\tau +V-(2n+1)\mu h F\bigr) \mu h ^{-1}F\sqrt{g}
\label{0-8}
\end{equation}
(with $\theta(\tau)=0,1$ as $\tau\le 0$, $\tau>0$ respectively) which implies that there are actually two different cases 
\begin{align}
&1\le \mu \le h^{-1}\label{0-9},\\
\intertext{and}
&h^{-1}\le \mu \le C_0h^{-\nu}\label{0-10}
\end{align}
when the answer will be of magnitude $h^{-2}$ and 
$(\mu h)^{-{1/(\nu-1)}}h^{-2}$ respectively; in the latter case the contribution to the answer will be given by the strip
\begin{equation}
Z=\{x: |x_1|\le {\bar\gamma}_1= C_0(\mu h)^{-{\frac 1 {\nu-1}}}\}.
\label{0-11}
\end{equation}

The standard rescaling procedure applied to results of \cite{IRO3} implies
\begin{equation}
\cR \Def |\int \Bigl(e(x,x,0)-\cE^\MW (x,0)\Bigr)\psi(x)\,dx|\le
Ch^{-1}
\label{0-12}
\end{equation}
but an aim of this paper is much better estimate, up to 
$O(\mu^{-{\frac 1 \nu}}h^{-1})$.

I am going to improve this remainder estimate.

\subsection{Results}

\begin{theorem}\label{thm-0-1} Let conditions $(\ref{0-2})$, $(\ref{0-3})$ and $(\ref{0-5})$ be fulfilled and let $\psi$ be supported in the small enough neighborhood of $\{x_1=0\}$. Then 

\smallskip
\noindent
(i) As $1\le \mu \le \mu^*_1= C h^{-\nu/3}$ estimate $\cR\le C\mu^{-1/\nu}h^{-1}$ holds; under some nondegeneracy assumptions one can push $\mu^*_1$ up;

\smallskip
\noindent
(ii)  As $ h^{-\nu/3}\mu \le \mu ^*_2= C h^{-\nu}$  estimate 
\begin{multline}
\cR^* \Def |\int \Bigl(e(x,x,0)-\cE^\MW (x,0)\Bigr)\psi(x)\,dx-\int \cE^\MW _\corr (x_2,0)\psi (0,x_2,0)\,dx_2 |\le\\
C\mu^{-1/\nu}h^{-1} +C h^{-\delta}
\label{0-13}
\end{multline}
holds  with arbitrarily small exponent $\delta>0$.
Here and below correction term $\cE^\MW_\corr$ is defined by $(\ref{3-42})$ in the terms of the eigenvalue counting function ${\bf n}_0$ of the axillary 1-dimensional Schr\"odinger operator with semiclassical parameter $\hbar= \mu^{1/\nu}h$.

\smallskip
\noindent
(iii) Under non-degeneracy condition \ref{2-105} for $W(x_2)=V(0,x_2)$ (see below) as 
$\mu \le \epsilon h^{-\nu}$  with small enough constant $\epsilon=\epsilon(\epsilon_0,m)$ (where $\epsilon_0$ is a constant in \ref{2-105}) and under non-degeneracy condition $(\ref{4-15})$ as $\epsilon h^{-\nu}\le \mu \le Ch^{-\nu}$ one can skip the last term in the right-hand expression of $(\ref{0-13})$;

\smallskip
\noindent
(iv) As $\mu \ge Ch^{-\nu}$\;   $|e(x,x,0)|\le C'\mu^{-s}$ with arbitrarily large exponent $s$.
\end{theorem}

\begin{remark}\label{rem-0-2} With some error making estimate sometimes not that sharp one can replace ${\bf n}_0$  by Bohr-Sommerfeld approximation for it, providing formulae (\ref{3-52}), \ref{3-52-*}, \ref{3-52-**} for $\cE^\MW_\corr$.
\end{remark}

\subsection{Plan}
I start from Classical Dynamics which is useful both for understanding and proofs of our main results. In Section 2 I consider corresponding quantum (semiclassical) dynamics and the remainder estimates. It appears that there will be \emph{forbidden zone\/} $\{|x_1|\ge {\bar \gamma}_1=C(\mu h)^{1/(\nu -1)}\}$ 
\footnote{\label{foot-1} Only as $Ch^{-1}\le \mu \le Ch^{-\nu}$; otherwise ${\bar\gamma}_1$ is artificially set to $\epsilon$.} \footnote{\label{foot-2} Here and below coordinate $x_1$ is defined so that $\Sigma=\{x_1=0\}$.}, \emph{outer zone\/} 
$\cZ_\out=\{{\bar \gamma}_0=C\mu^{-1/\nu}\le |x_1|\le {\bar \gamma}_1\}$ and inner zone $\cZ_\inn=\{|x_1|\le {\bar\gamma}_0\}$ divided in \emph{non-periodic zone\/} and \emph{periodic zone\/} with the analysis in the latter one the most difficult and interesting. The width (in the phase space) of the periodic zone is defined by the uncertainty principle and it increases as $\mu$ grows. 

In Section 3 I instead of implicit main part of asymptotics derived in Section 2 provide much more explicit answer $\cE^\MW+\cE^\MW_\corr$. 

All the analysis in sections 2,3 is under assumption 
$\mu\le h^{-\nu+\delta}$ while  section 4 is devoted to the similar analysis in the case of \emph{superstrong magnetic field\/} 
$h^{-\nu+\delta}\le \mu \le \epsilon h^{-\nu}$ which is in some sense easier and in some sense more difficult than the previous one because $\hbar$ is not very small anymore. Also in Section 4 I consider \emph{ultra-strong magnetic field\/} $\epsilon h^{-\nu}\le \mu\le Ch^{-\nu}$ when $\hbar$ is not even small anymore; as I mentioned case $\mu\ge Ch^{-\nu}$ is trivial.

Finally, Section 5 is devoted to generalization for vanishing $V$ and discussion of other generalizations.

\sect{Classical dynamics}
In this section I consider classical dynamics on energy level 0 described by Hamiltonian
\begin{equation}
a= {\frac 1 2}\Bigl(\sum_{j,k}g^{jk}(x)p_jp_k -V\Bigr),\qquad p_j=\xi_j-\mu V_j
\label{1-1}
\end{equation}
where all conditions of subsection 0.1 are assumed to be fulfilled.

\subsection{Pilot-model}

Let us consider first a model
\begin{equation}
a={\frac 1 2}\Bigl(\xi_1^2 + \bigl(\xi_2- {\frac 1 \nu}|x_1|^\nu\bigr)^2-1\Bigr)
\label{1-2}
\end{equation}
with $\nu\ge 2$, corresponding to $\mu=1$, $V=1$.

\begin{proposition}\label{prop-1-1} For symbol $(\ref{1-2})$ 

\smallskip
\noindent
(i) Along trajectories  $\xi_2=\const$, and 
\begin{description}

\item{(a)} if $\xi_2 >1$   on energy level $0$ these trajectories  oscillate \underline{either} between $x_1=b_1$ and  $x_1=b_2$ \underline{or}
between $x_1=-b_1$ and $x_1=-b_2$ with $b_1=\bigl((\xi_2-1)\nu\bigr)^{1/\nu}$, $b_2=\bigl((\xi_2+1)\nu\bigr)^{1/\nu}$; 

\item{(b)} if $\xi_2=1$ these trajectories vary between $b_1=0$ and $b_2=1$ 
or between $-1$ and $0$; there is also exceptional trajectory $x_1=\xi_1=0$, 
$x_2=t+\const$;

\item{(c)} if $-1<\xi_2<1$ these trajectories oscillate between $x_1=b_1=-b_2$ and $x_1=b_2$;  

\item{(d)} if $\xi_2=-1$ energy level 0 degenerates into $\{x_1=\xi_1=0\}$ and this trajectory is $x_1=\xi_1=0$, $x_2=\xi_2t+\const$; 

\item{(e)} energy level 0 is empty as $\xi_2<-1$;

\end{description}

\smallskip
\noindent
(ii) Along each of these trajectories  
$(x_1,\xi_1)$ are periodic with period  
\begin{equation}
T(\xi_2)=2\int_{b_1}^{b_2} {\frac 1 {\sqrt{1-(k-|y|^\nu/\nu)^2}}}\,dy
\label{1-3}
\end{equation}
and 
$x_2(t)=v(\xi_2)t+ {\tilde x}_2(t)$ where $v=T^{-1}I$,
\begin{equation}
I(\xi_2)=\int_{b_1}^{b_2} {\frac {k- |y|^\nu/\nu} {\sqrt{1-(k-|y|^\nu/\nu)^2}}}\,dy,
\label{1-4}
\end{equation}
and  ${\tilde x}_2(t)$ is $T(\xi_2)$-periodic function.
\end{proposition}

\begin{proof} Proof is obvious consequence of the fact that evolution in $(x_1,\xi_1)$ is described by 1-dimensional Hamiltonian (\ref{1-2}) while $\xi_2=k=\const$. Potential 
$W=W(k)=(k- {\frac 1 \nu}|x|_1^\nu)^2-1$ has two wells as $k>0$ and one well as 
$k\le 0$; as $k=0$ it has a flat bottom. Note that the bottoms of the potential are in $\pm b(k)$, $b(k)=(k\nu)^{1/\nu}$. 

\begin{figure}[h]
\centering
\subfigure[$k=1.5$, two wells]{
\label{fig-1a}
Ê\includegraphics[width=.3\textwidth]{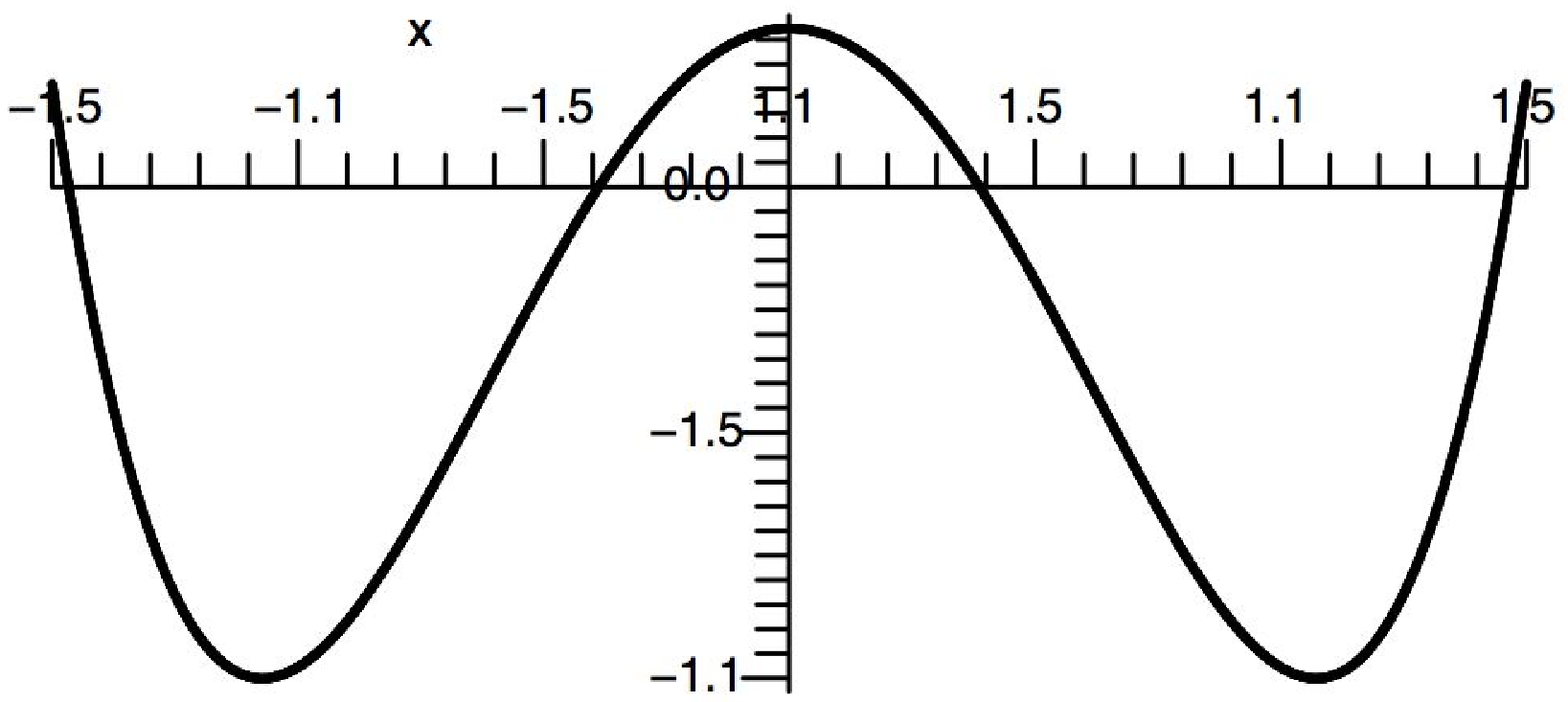}}
\subfigure[$k=1$, two touching wells]{
ÊÊÊÊÊ\label{fig-1b}
ÊÊÊÊÊÊ\includegraphics[width=.3\textwidth]{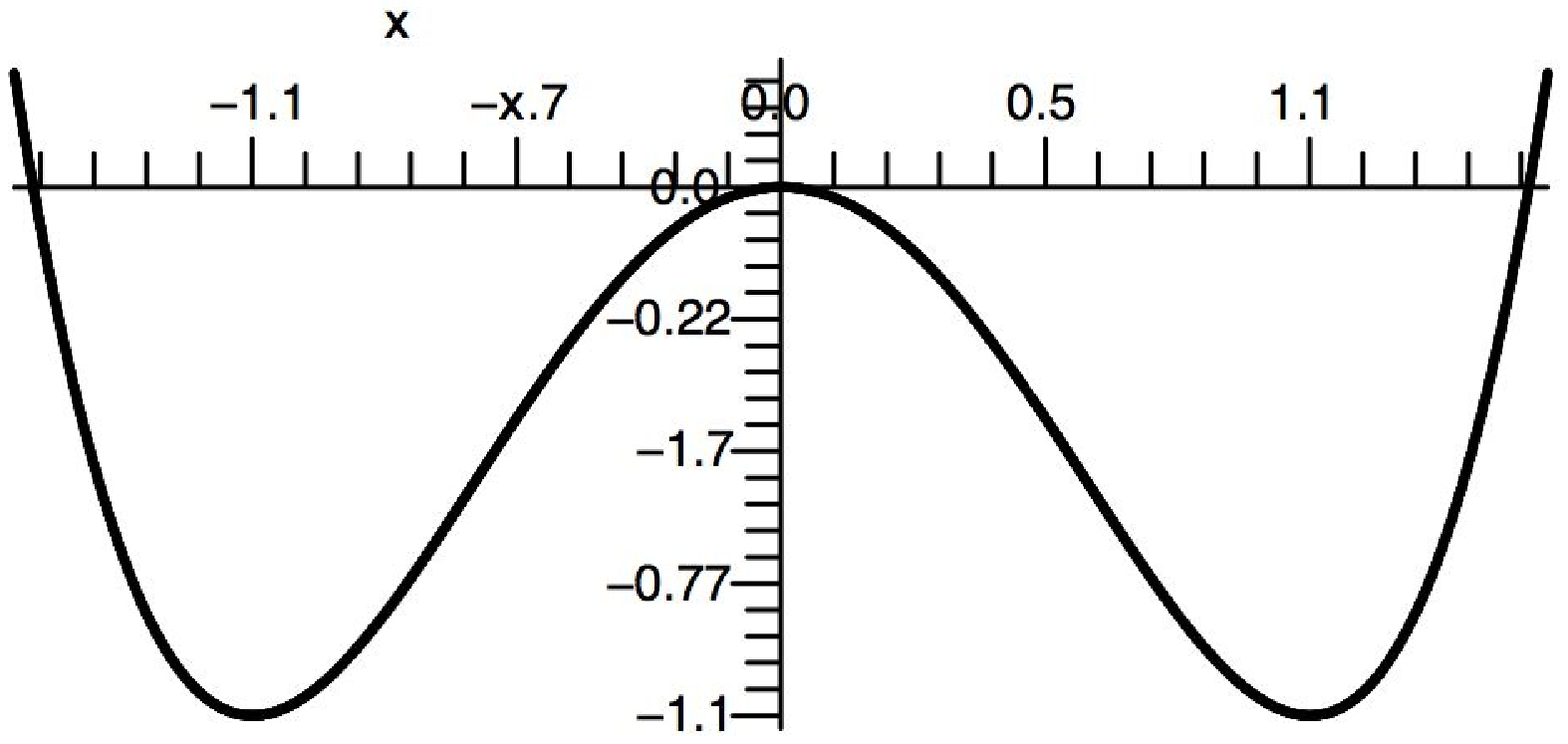}}
\subfigure[$k=0.9$, one well]{
ÊÊÊÊÊ\label{fig-1c}
ÊÊÊÊÊÊ\includegraphics[width=.3\textwidth]{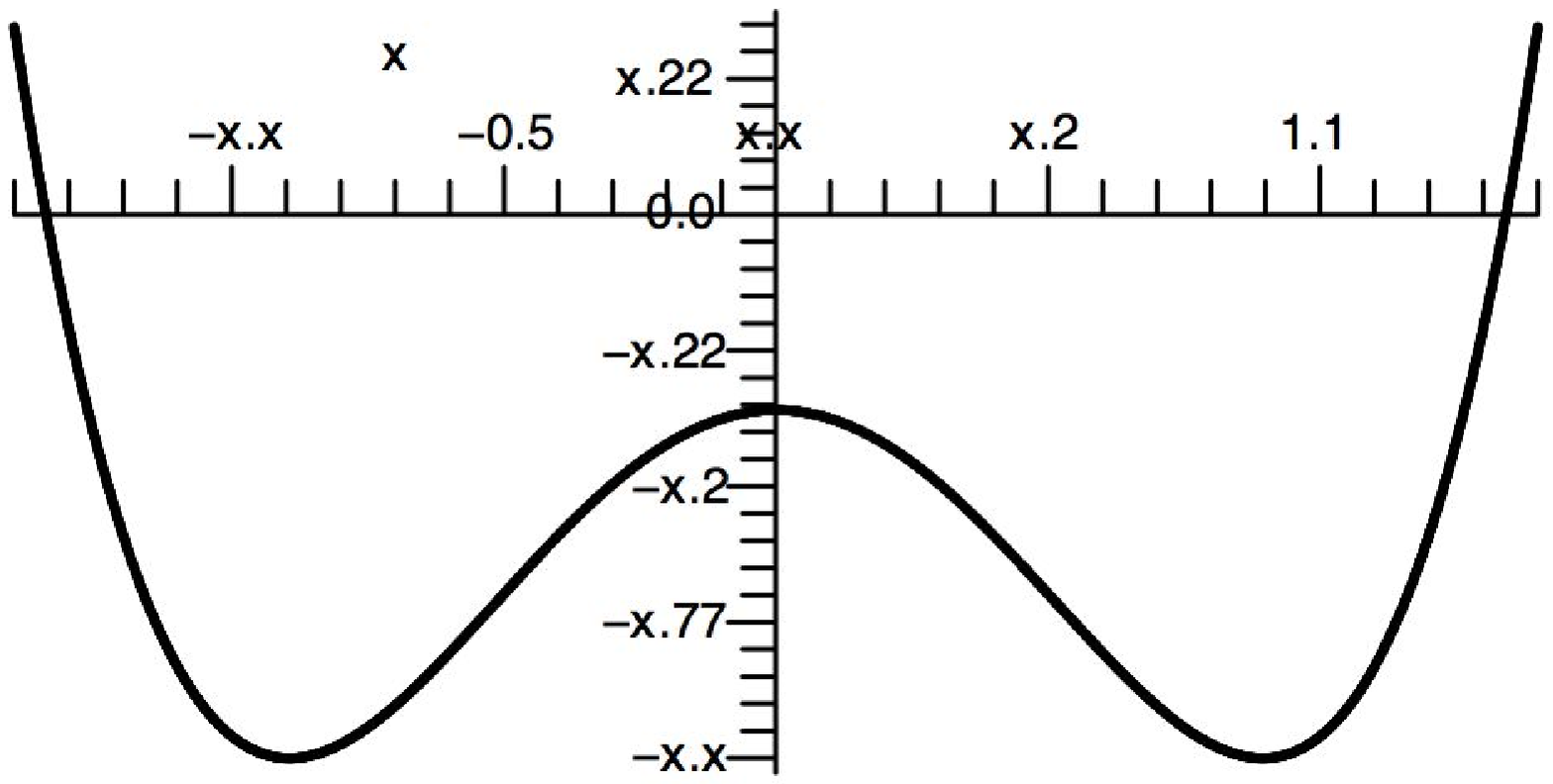}}
\caption{Graphs of $W(k)$, 1-st pilot-model}
\end{figure}

Further, as $k>1$ we have $W(0)>0$ and thus the particle lives in one well
oscillating either between $b_1(k)=( (k-1)\nu)^{1/\nu}$ and 
$b_2(k)=( (k+1)\nu)^{1/\nu}$ or between $-b_2(k)$ and $-b_1(k)$. As $k<1$ we have either one-well potential or $W(0)<0$ so particle oscillates between $b_1(k)=-b_2(k)$ and $b_2(k)$. As $k=1$ particle moves between $b_1(k)=0$ (as 
$t\to \pm\infty$) and $\pm b_2(k)$ (as long as $\nu \ge 2$).
\end{proof}

\begin{proposition}\label{prop-1-2} (i)  $v(k)>0$ and
both $T(k)$ and $v(k)$ are decreasing functions of $k>1$. 

\smallskip
\noindent
(ii)   The width of the wells is $\sim 2b^{(1-\nu)/\nu}$ as $k\to +\infty$;
 $T(k)\sim 2\pi b^{1-\nu}$, $v(k)\sim {\frac 1 2}(\nu-1)b^{-\nu}$ as 
$k\to +\infty$, $b(k)= (k\nu )^{1/\nu}$;

\smallskip
\noindent
(iii)  $T\sim 2\bigl|\log |k-1|\bigr|$ (for $\nu=2$),
$T\sim \const |k-1|^{{\frac 1 \nu} -{\frac 1 2}}$ (for $\nu>2$) and 
$v\sim 1$ as $k\to 1+0$.
\end{proposition}

\begin{proof} Changing $z=y^\nu/\nu-k$ we arrive to
\begin{align}
T=&2\int_{-1}^1 ((k+z)\nu)^{{\frac 1 \nu}-1}(1-z^2)^{-{\frac 1 2}}\,dz,
\label{1-5}\\
I=&-2\int_{-1}^1 z((k+z)\nu)^{{\frac 1 \nu}-1}(1-z^2)^{-{\frac 1 2}}\,dz=
\label{1-6}\\
&\int_{-1}^1 
z\Bigl( ((k-z)\nu)^{{\frac 1 \nu}-1} - ((k+z)\nu)^{{\frac 1 \nu}-1}\Bigr) 
(1-z^2)^{-{\frac 1 2}}\,dz,\notag
\end{align}
which implies that $I>0$ as $k>1$.  Asymptotics 
$T\sim 2\pi (k\nu)^{{\frac 1 \nu}-1}$,
$I\sim \pi (\nu -1)(k\nu)^{{\frac 1 \nu}-2}$ as $k\to +\infty$ follow from
(\ref{1-5}),(\ref{1-6}) as well.

Also (\ref{1-5}),(\ref{1-6}) imply that
\begin{align}
{\frac {\partial T}{\partial k}}=&
-2(\nu-1)\int_{-1}^1 \bigl((k+z)\nu\bigr)^{{\frac 1 \nu}-2}
(1-z^2)^{-{\frac 1 2}}\,dz, \label{1-7}\\
{\frac {\partial I}{\partial k}}=&
-2(\nu-1)\int_{-1}^1 z\bigl((k+z)\nu\bigr)^{{\frac 1 \nu}-2}
(1-z^2)^{-{\frac 1 2}}\,dz,
\label{1-8}\\
{\frac {\partial v}{\partial k}}=& 
-\nu(\nu-1)T^{-2}\int_{-1}^1\int_{-1}^1(z-z_1)^2 
\bigl((k+z)\nu\bigr)^{{\frac 1 \nu}-2}
\bigl((k+z_1)\nu\bigr)^{{\frac 1 \nu}-2}\times\label{1-9}\\
&\hskip220pt(1-z^2)^{-{\frac 1 2}}\,dz,\notag
\end{align}
implying monotonicity properties.

Asymptotics as $k\to 1+0$ follow directly from (\ref{1-3}),(\ref{1-4}).
\end{proof}

\begin{remark}\label{rem-1-3} In the strong magnetic field the drift is described by ${\frac 1 2}\cL f^{-1}=(\partial_{x_2}f,-\partial_{x_1}f)=(0,(\nu-1)x_1^{-\nu})$
(where ${\frac 1 2}$ comes from this factor in the front of Hamiltonian)
and this matches to asymptotics of $v$ as $\xi_2\to +\infty$.
\end{remark}

As $-1<k<1$ one can rewrite
\begin{align}
&T=4\int_0^{b_2}\Bigl(1-(k-y^\nu/\nu)^2\Bigr)^{-{\frac 1 2}}\, dy=\label{1-10}\\
&\qquad\qquad\qquad
4(k+1)^{{\frac  1 \nu}-{\frac 1 2}} \nu^{\frac 1 \nu}
\int_0^1(1-y^\nu)^{-{\frac 1 2}}\Bigl(2-(k+1)(1-y^\nu)\Bigr)^{-{\frac 1 2}}\,dy=
\notag\\
&\qquad\qquad\qquad\qquad\qquad\qquad
4(k+1)^{{\frac  1 \nu}-1} \nu^{\frac 1 \nu}
\int_0^1(1-y^\nu)^{-{\frac 1 2}}\varphi (u,y^\nu)\,dy
\notag\\
\intertext{and}
&I=4\int_0^{b_2}(k-y^\nu/\nu)
\Bigl(1-(k-y^\nu/\nu)^2\Bigr)^{-{\frac 1 2}}\, dy=\label{1-11}\\
&\qquad\qquad\qquad\qquad\qquad\qquad
4((k+1)\nu)^{\frac 1 \nu}
\int_0^1  (1-y^\nu)^{-{\frac 1 2}}\varphi_1 (u,y^\nu)\,dy
\notag
\end{align}
with $u= k /(k+1) \in(-\infty,{\frac 1 2})$, 
$\varphi (u,z)=(1-2u+z)^{-{\frac 1 2}}$, 
$\varphi_1 (u,z)=(u-z )(1-2u+z)^{-{\frac 1 2}}$.

Let $T_1(u)$, $I_1(u)$ denote integrals in the right-hand expressions
of (\ref{1-10}), (\ref{1-11}). Then 
${\frac{\partial T_1}{\partial u}}>0$, 
${\frac{\partial I_1}{\partial u}}=  (1-u) {\frac{\partial T_1}{\partial u}}>0$
because ${\frac{\partial \varphi}{\partial u}}=(1-2u+z)^{-{\frac 3 2}}$,
${\frac{\partial \varphi_1}{\partial u}}=(1-u)(1-2u+z)^{-{\frac 3 2}}$.

Also
\begin{equation*}
{\frac{\partial\ }{\partial u}}{\frac {I_1}{T_1}}=
4T_1^{-2}{\frac{\partial T_1}{\partial u}}\times 
\int_0^1 (1-y^\nu)^{-{\frac 1 2}} \bigl(1-2u+y^\nu\bigr)^{\frac 1 2}\,dy
> 1
\end{equation*}
due to Cauchy inequality and therefore
\begin{equation*}
{\frac{\partial v}{\partial k}}=
(1-u){\frac{\partial\ }{\partial u}}{\frac {I_1}{T_1}} + 
{\frac {I_1}{T_1}} > {\frac 1{T_1}}\bigl(I_1+(1-u)T_1\bigr)=
{\frac 4 {T_1}}
\int_0^1 (1-y^\nu)^{\frac 1 2} \bigl(1-2u+y^\nu\bigr)^{-{\frac 1 2}}\,dy >0.
\end{equation*}

\begin{proposition}\label{prop-1-4} (i) There exists $k^*=k^*_\nu \in (0,1)$
such that $v(\xi_2)(\xi_2-k^*)^{-1}>0$ and it is both bounded and disjoint from $0$ as $k\in (-1+\epsilon, 1)$;

\smallskip
\noindent
(ii) $v(k)$ is monotone increasing at $(-1,1)$;

\smallskip
\noindent
(iii)  $T(k)\sim  \const (k+1)^{\frac 1 \nu} $  and $v(k)\sim -1$ as 
$k\to -1+0$; 

\smallskip
\noindent
(iv)  $T\sim 4\bigl|\log |k-1|\bigr|$ (for $\nu=2$),
$T\sim 2\const |k-1|^{{\frac 1 \nu} -{\frac 1 2}}$ (for $\nu>2$) and 
$v\sim 1$ as $k\to 1+0$.
\end{proposition}

\begin{proof} (ii) is already proven; (iii),(iv) are easy and (i) follows from
(ii)-(iii) and from $v(0)<0$.
\end{proof}

\begin{remark}\label{rem-1-5}
(i) $T(k)$ is monotone decreasing function of $k$ for $\nu=2$; not sure if this is  the case for larger $\nu$;

\smallskip
\noindent
(ii) ${\frac{dx_2}{dt}}>0$ iff  $|x_1|< (k\nu)^{\frac 1 \nu}$;
in particular  ${\frac{dx_2}{dt}}<0$ on the whole trajectory iff
$-1<k<0$;

\smallskip
\noindent
(iii) One can prove easily that 
${\frac{\partial v}{\partial k}}\sim 
\mp \const |k-1|^{-{\frac 1 \nu}-{\frac 1 2}}|\log |k-1||^{-2\delta_{\nu 2}}$ as $k\to 1\pm 0$;

\smallskip
\noindent
(iv) In particular, $x_2(t)$ is periodic iff $\xi_2=k^*$;

\smallskip
\noindent
(v) Since 
\begin{equation*}
{\frac{\partial I_1}{\partial \nu}}=
\int_0^1 (u-1)^2(1-y^\nu)^{-{\frac 3 2}}(1-2u+y^\nu)^{-{\frac 3 2}}
y^\nu (-\log y)\,dy >0
\end{equation*}
we conclude that $k^*_\nu$ monotonically decreases; one can see easily that 
$k^*_\nu\to +0$ as $\nu \to +\infty$;

\smallskip
\noindent
(vi) Maple experiments show that $k^*_2\approx 0.65$ (but not $2/3$).
\end{remark}

\begin{figure}
\centering
\subfigure[$k=10$,  $w\approx 0.34$]{
ÊÊÊÊÊ\label{fig-2a}
ÊÊÊÊÊ\includegraphics[width=.48\textwidth]{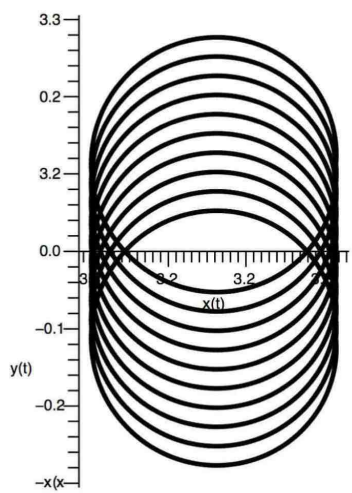}}
\subfigure[$k=2$,  $w\approx 0.73$ ]{
ÊÊÊÊÊ\label{fig-2b}
ÊÊÊÊÊÊ\includegraphics[width=.49\textwidth]{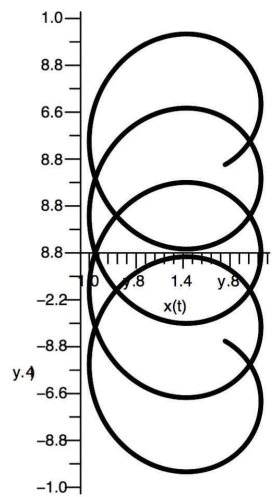}}\\

\vspace{.3in}
\subfigure[$k=1.1$,  $w\approx 1.13$]{
ÊÊÊÊÊÊ\label{fig-2c}
ÊÊÊÊÊÊ\includegraphics[width=.48\textwidth]{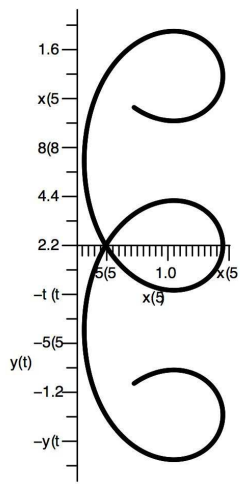}}
\subfigure[$k=1$,  $w\approx 1.13$, first critical case]{
\label{fig-2d}
\includegraphics[width=.49\textwidth]{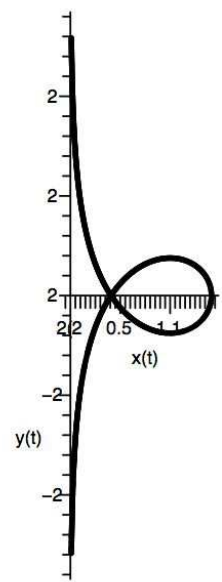}}

\caption{\label{fig-2} Drift up; and rotation clockwise,
each figure has its mirror with respect to $x_2$-axis with drift up and rotation counter-clockwise; $w=\sqrt{k+1}-\sqrt{k-1}$ is the width;}
\end{figure}

The following series of figures \ref{fig-2} show $(x_1,x_2)$-trajectories at level $0$ for different values of $\xi_2=k$ as $\nu=2$. In particular figure \ref{fig-2g} shows periodic curve.

\begin{proposition}\label{prop-1-6}
There exist  functions $Z=Z(x_1,\xi_1)$, $\alpha=\alpha(x_1,\xi)$
and $\beta=\beta (x_1,\xi_2)$, such that
\begin{equation}
\{a, x_2-Z\}=\xi_2-\{a,Z\}=\alpha (\xi_2-k^*)+\beta a(x_1,\xi).
\label{1-12}
\end{equation}
Further, $Z$ is odd and $\alpha$, $\beta$ are even with respect  to each of $x_1$, $\xi_1$.
\end{proposition}

\begin{proof}
Consider now periodic curve as $\xi_2=k^*$ with the ``center'' at 0. It is described by equations in $(x,\xi)$-phase space
\begin{equation}
\xi_2=k^*,\quad a(x_1,\xi)=0,\quad x_2=Z(x_1,\xi_1)
\label{1-13}
\end{equation}
with smooth odd with respect to each $x_1$, $\xi_1$ function $Z$. We cannot express it via $x_1$ or $\xi_1$ alone. Then
\begin{equation*}
\{a, x_2-Z\}=\xi_2-\{a,Z\}=0\qquad {\text as\ } \xi_2=k^*,\ a(x_1,\xi)=0.
\end{equation*}
Note that $d a$ and $d\xi_2$ are linearly independent unless $\xi_1=0$ and
$(\xi_2-|x_1|^\nu/\nu)x_1=0$ which for $a=0$ means that $\xi_2=\pm 1$. Therefore 
we arrive to (\ref{1-12}).
\end{proof}

\begin{remark}\label{rem-1-7}
For even integer $\nu$ all these functions are analytic.
\end{remark}

\setcounter{figure}{1}
\setcounter{subfigure}{4}

\begin{figure}
\centering
\subfigure[$k=0.9$,  $w\approx  4.32$]{
ÊÊÊÊÊ\label{fig-2e}
ÊÊÊÊÊ\includegraphics[width=.48\textwidth]{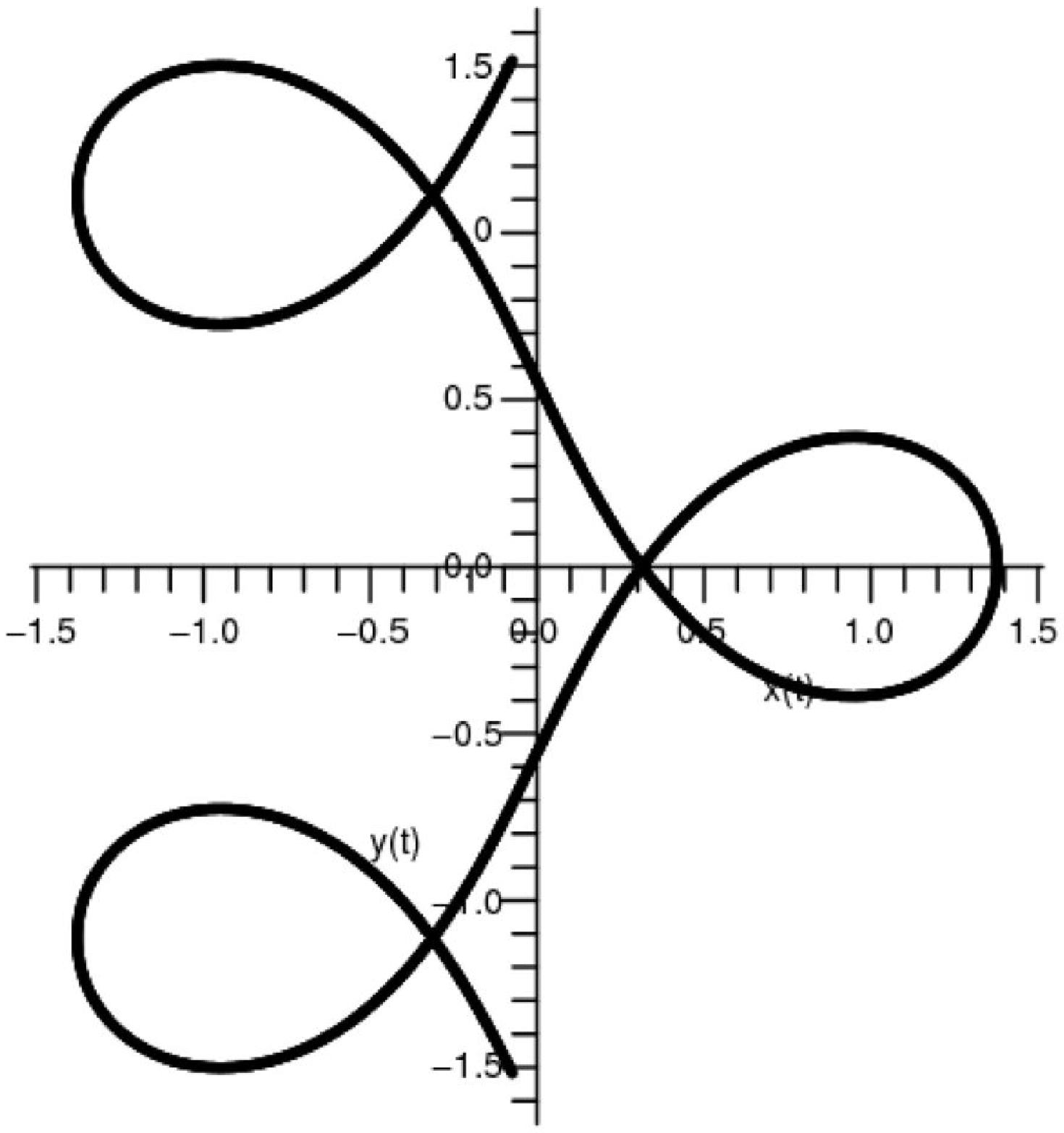}}
\subfigure[$k=0.7$,  $w\approx 3.78$]{
ÊÊÊÊÊ\label{fig-2f}
ÊÊÊÊÊÊ\includegraphics[width=.49\textwidth]{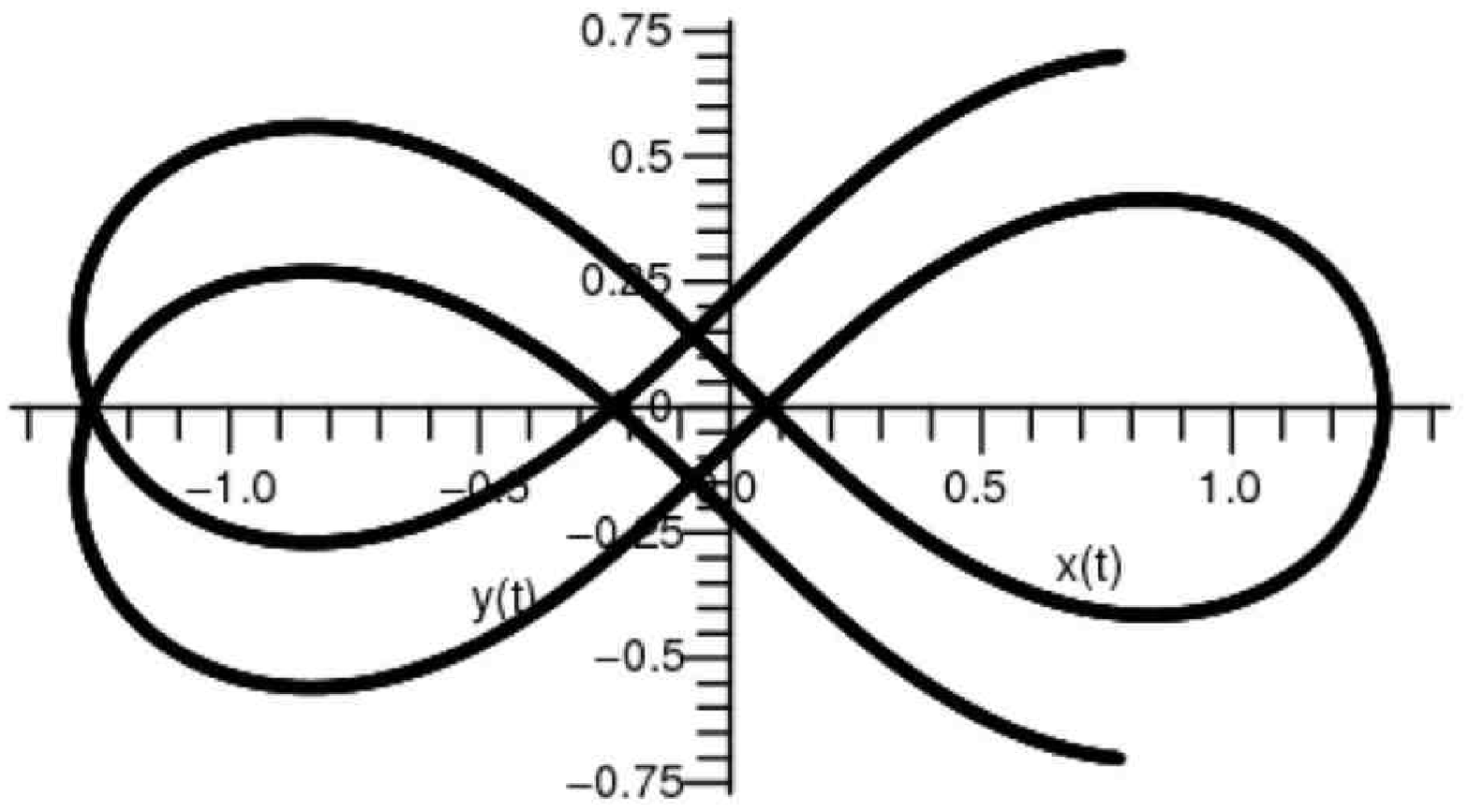}}\\

\vspace{.3in}
\subfigure[$k \approx 0.65$,  $w\approx  3.65$]{
ÊÊÊÊÊÊ\label{fig-2g}
ÊÊÊÊÊÊ\includegraphics[width=.48\textwidth]{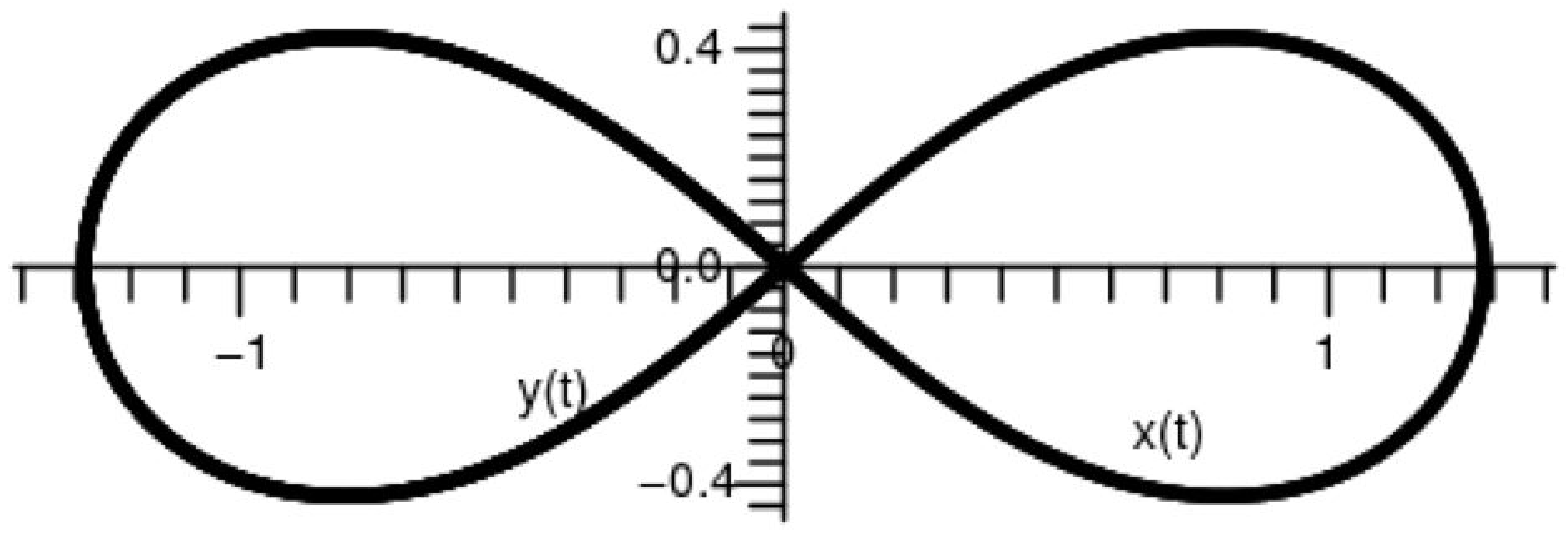}}
\subfigure[$k=0.6$,  $w\approx 3.52$, second critical case]{
\label{fig-2h}
\includegraphics[width=.48\textwidth]{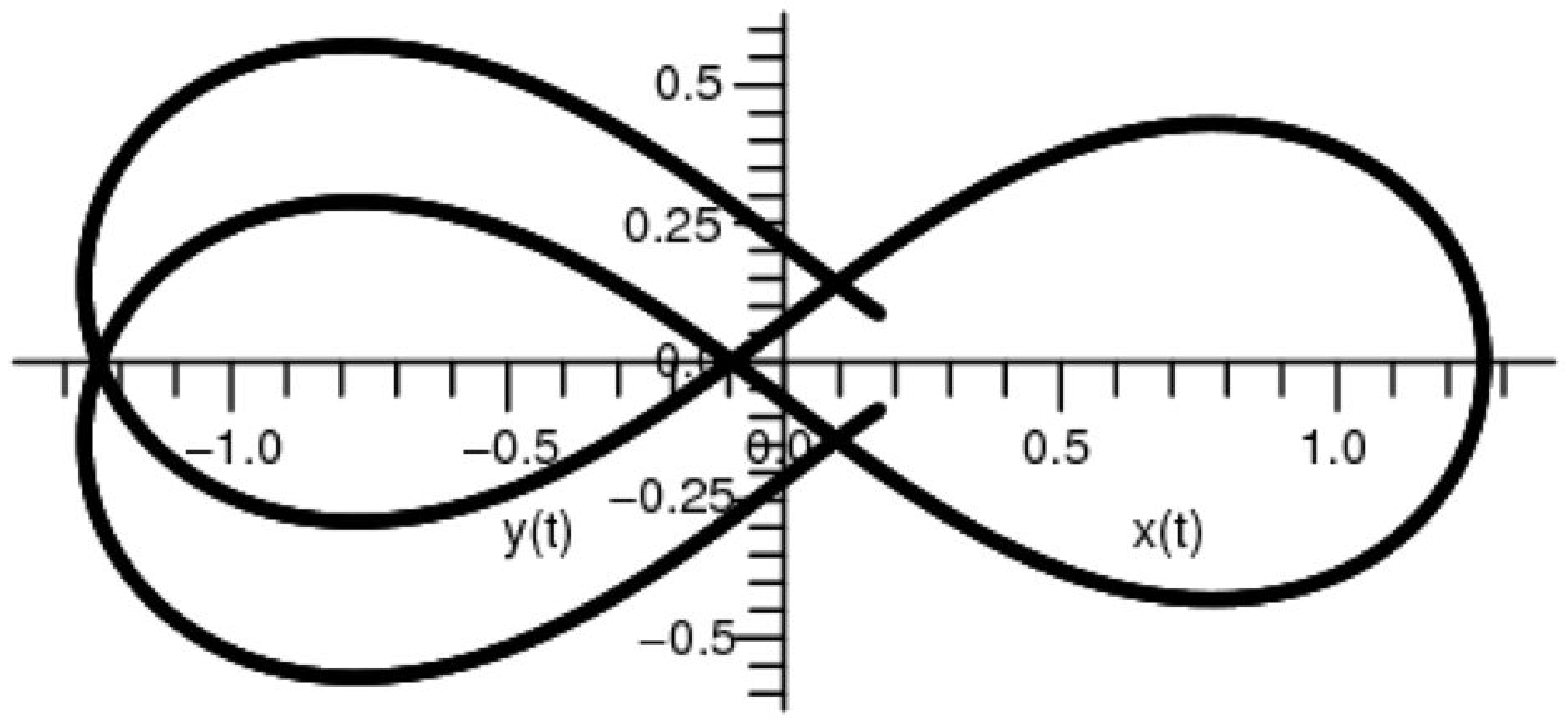}}

\caption{\label{fig-2+} (continued) As $k$ decreases drift up changes to drift down; on
\ref{fig-2g} movement is periodic; $w=2\sqrt{k+1}$ is the width;}
\end{figure}

\setcounter{figure}{1}
\setcounter{subfigure}{8}

\begin{figure}
\centering
\subfigure[$k=0.4$,  $w\approx 2.28$]{
ÊÊÊÊÊ\label{fig-2i}
ÊÊÊÊÊ\includegraphics[width=.48\textwidth]{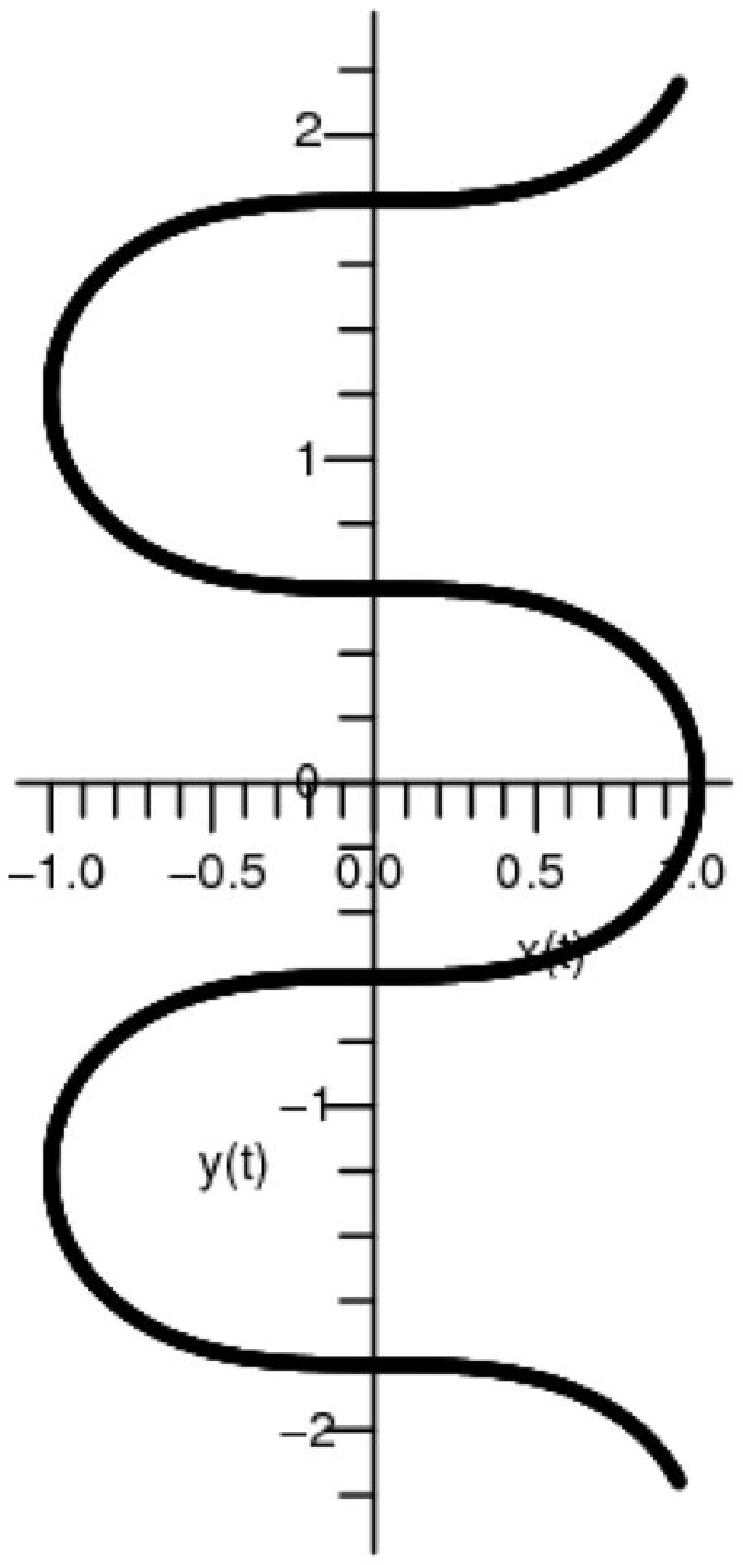}}
\subfigure[$k=0$,  $w=2$]{
ÊÊÊÊÊ\label{fig-2j}
ÊÊÊÊÊÊ\includegraphics[,width=.49\textwidth]{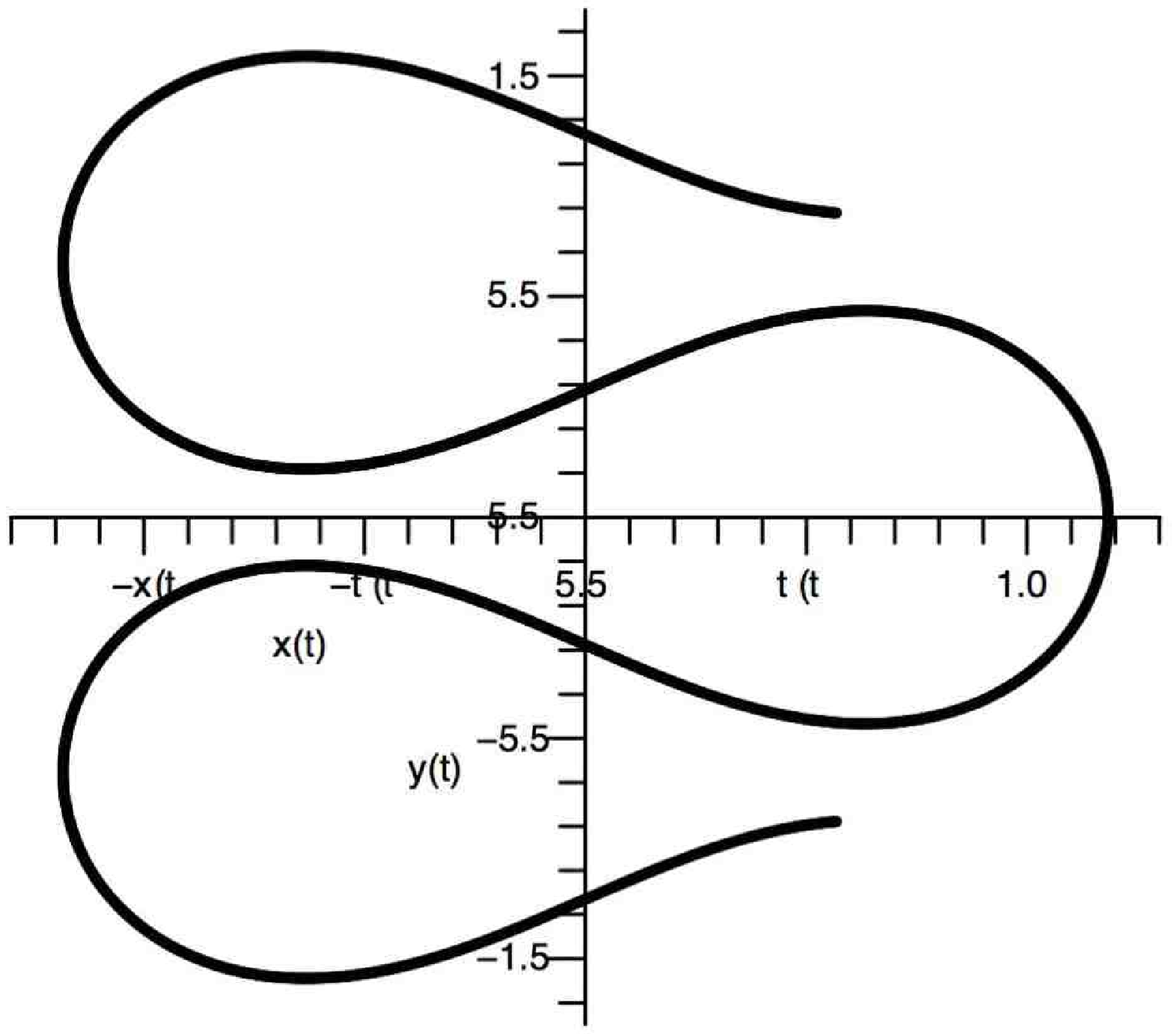}}\\

\vspace{.3in}
\subfigure[$k=-0.5$,  $w\approx 1.41$]{
ÊÊÊÊÊÊ\label{fig-2k}
ÊÊÊÊÊ\includegraphics[width=.48\textwidth]{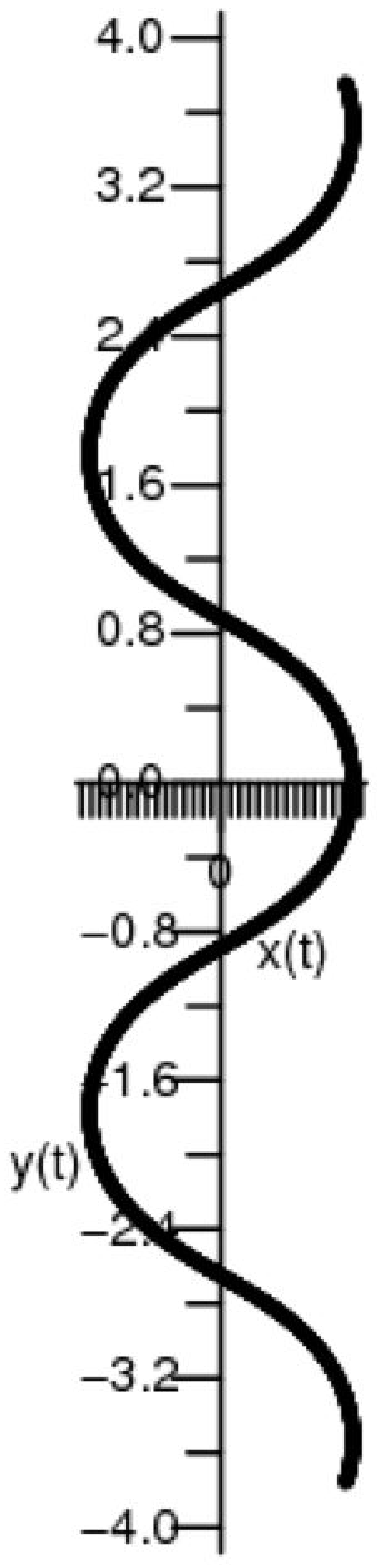}}
\subfigure[$k=-0.9$,  $w\approx 0.63$]{
\label{fig-2l}
\includegraphics[width=.48\textwidth]{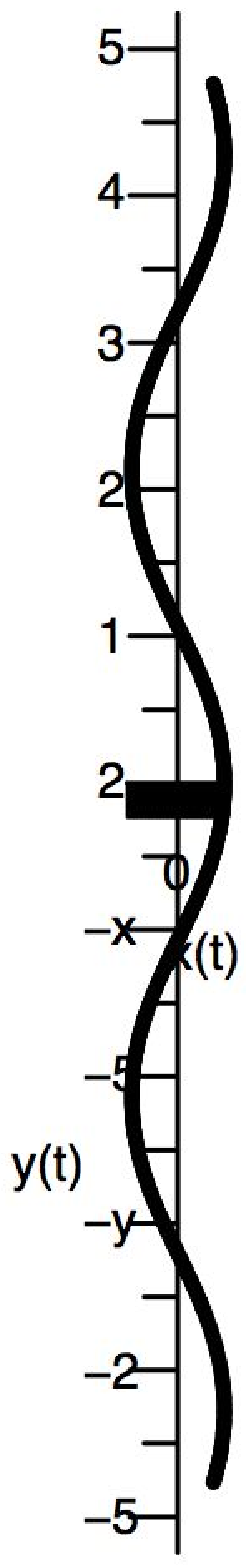}}

\caption{\label{fig-2++}(continued) As $k$ decays movement is straighten; drift is down.}
\end{figure}

\subsection{Another pilot-model}

Consider now another model
\begin{equation}
a={\frac 1 2}
\Bigl(\xi_1^2 + \bigl(\xi_2- {\frac 1 \nu}|x_1|^\nu\sign (x)\bigr)^2-1\Bigr)
\label{1-14}
\end{equation}
with $\nu\ge 2$, corresponding to $\mu=1$, $V=1$. Again we are looking at energy level 0.  Consider again the graph of potential.
\begin{figure}[h]
\centering
\subfigure[$k=1.5$,  well in $x_1>0$]{
\label{fig-3a}
Ê\includegraphics[width=.3\textwidth]{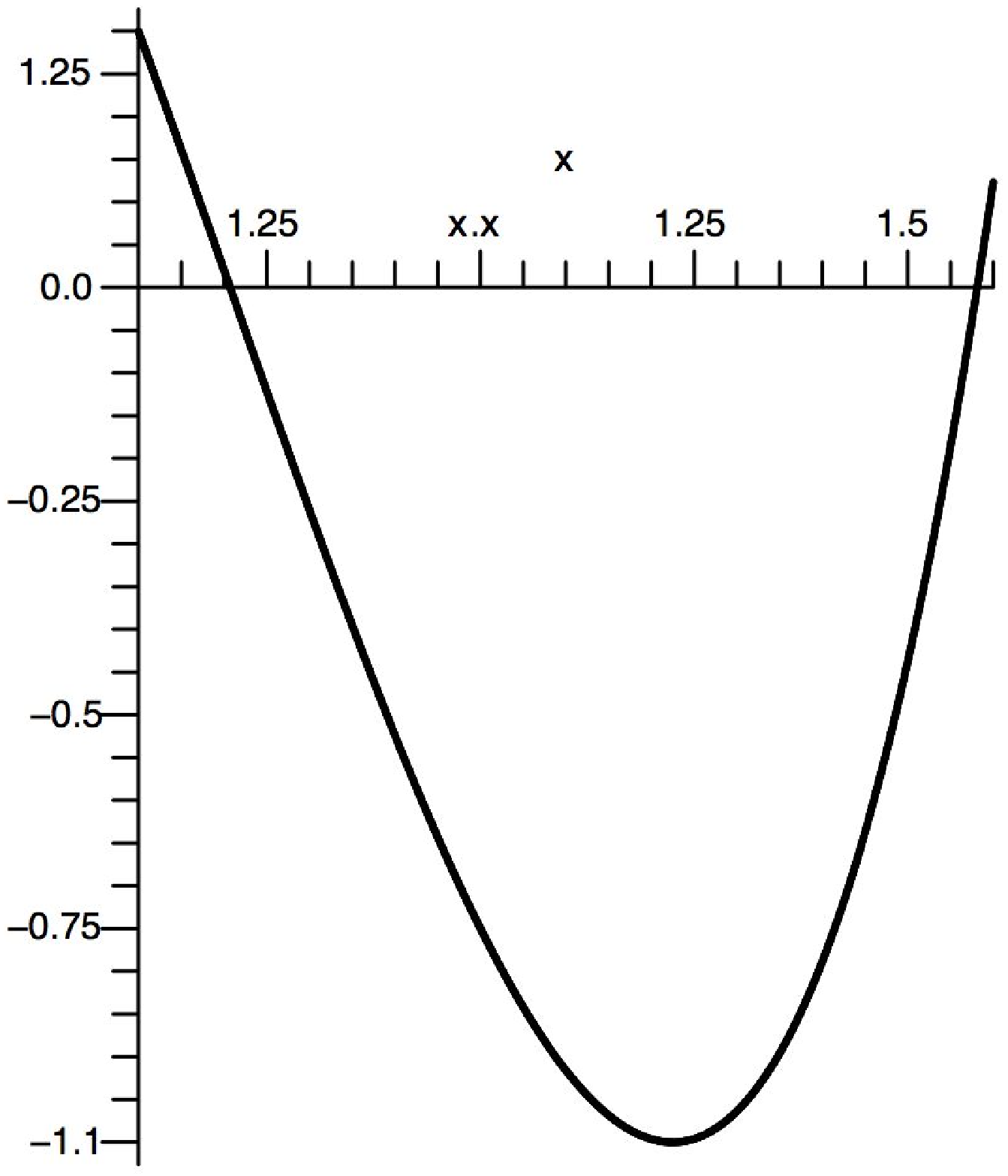}}
\subfigure[$k=1$, well touches $x_1=0$]{
ÊÊÊÊÊ\label{fig-3b}
ÊÊÊÊÊÊ\includegraphics[width=.3\textwidth]{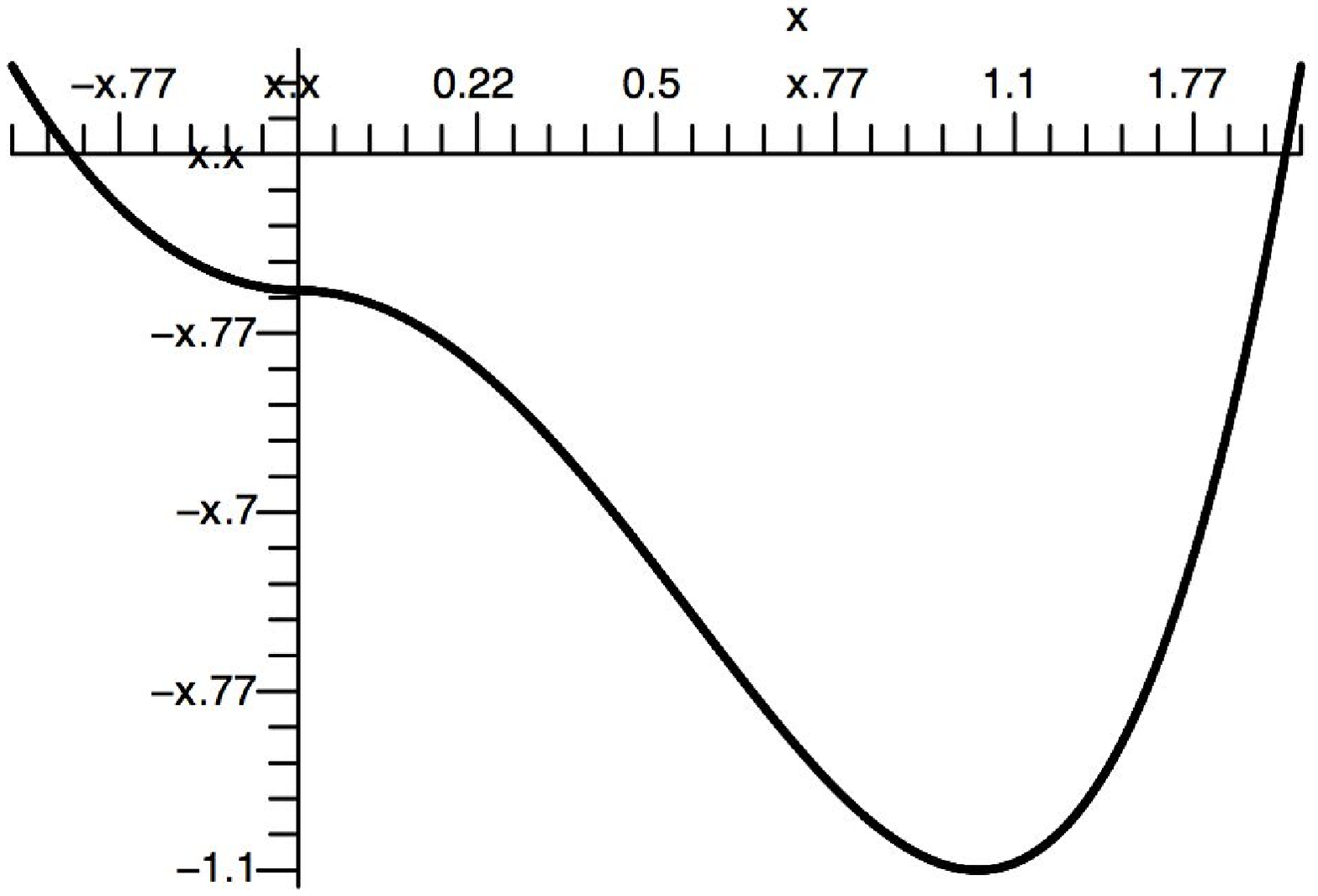}}
\subfigure[$k=0.9$, well contains $x_1=0$]{
ÊÊÊÊÊ\label{fig-3c}
ÊÊÊÊÊÊ\includegraphics[width=.3\textwidth]{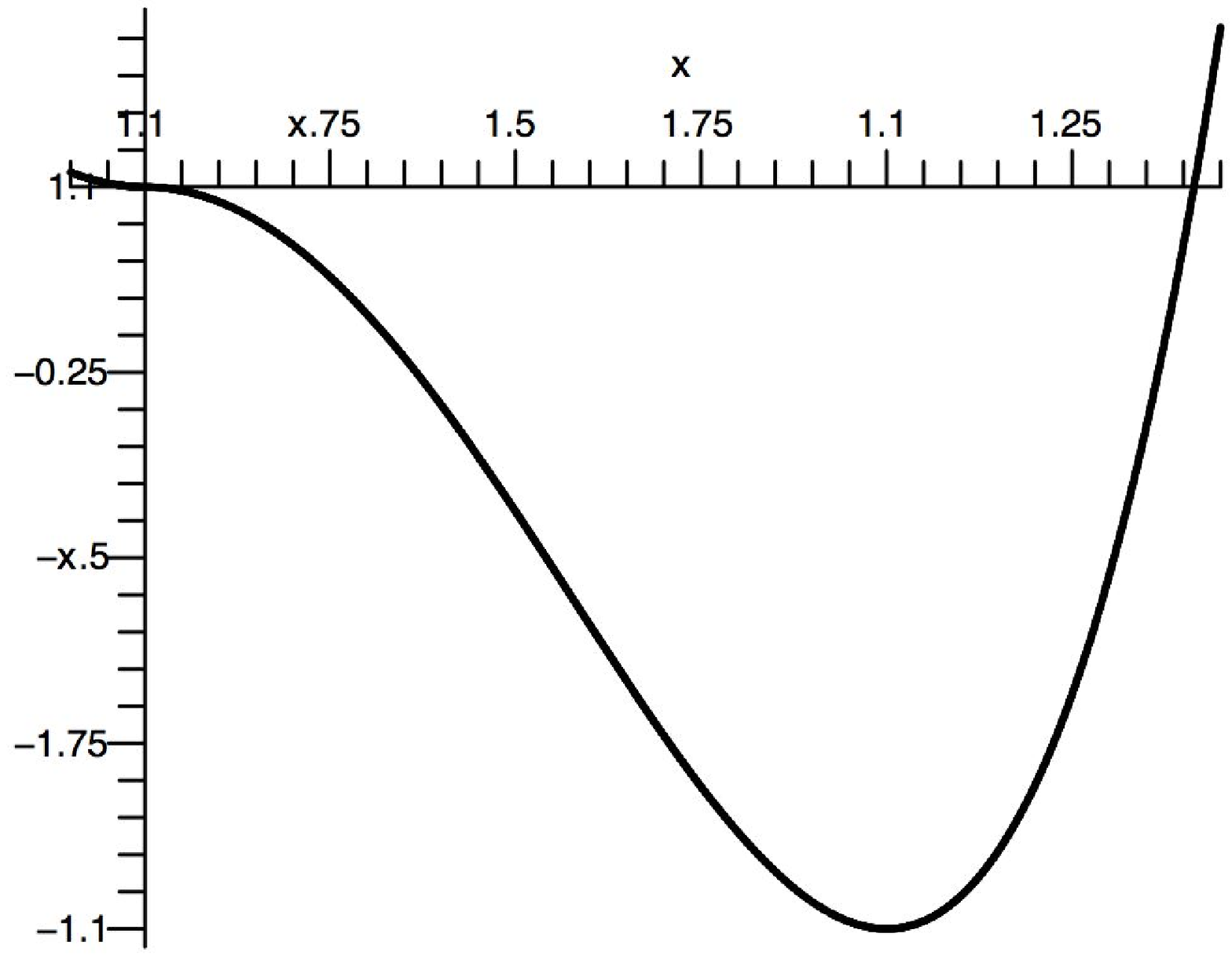}}
\caption{Graphs of $W(k)$, 2-nd pilot-model}
\end{figure}

\begin{figure}
\centering
\subfigure[$k=0.9$]{
ÊÊÊÊÊ\label{fig-4e}
ÊÊÊÊÊ\includegraphics[width=.48\textwidth]{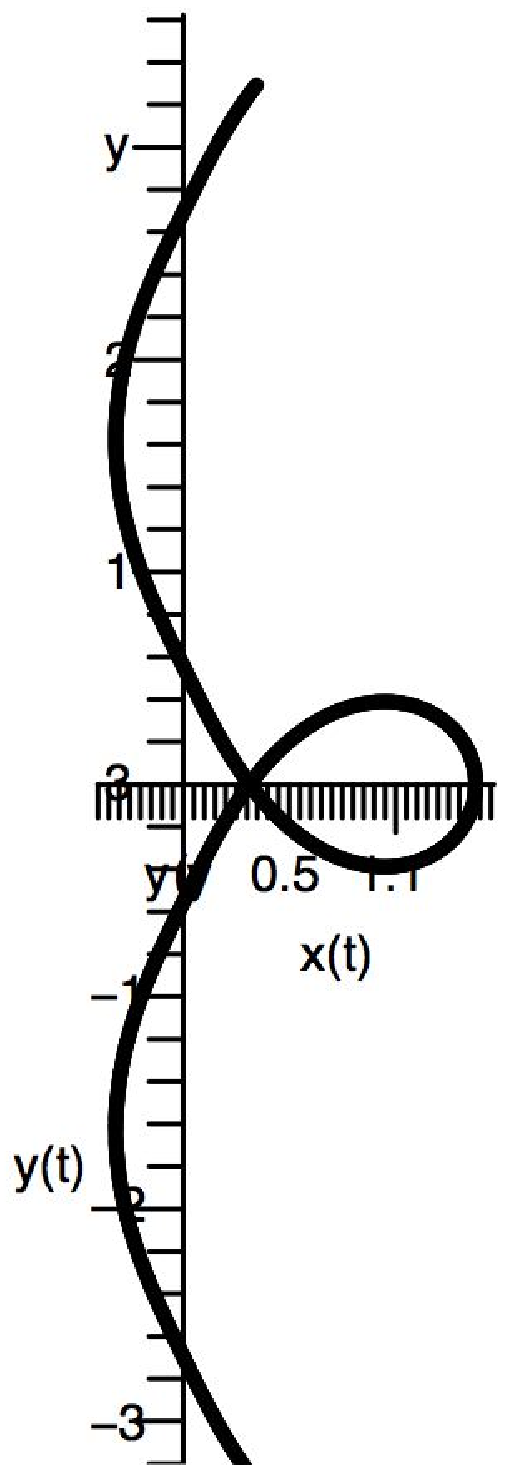}}
\subfigure[$k=0.5$]{
ÊÊÊÊÊ\label{fig-4f}
ÊÊÊÊÊÊ\includegraphics[width=.49\textwidth]{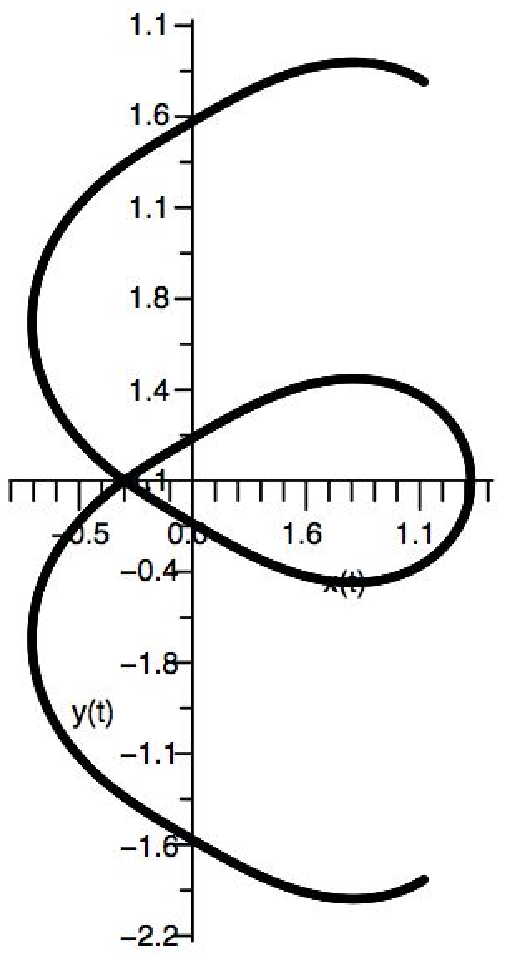}}\\

\vspace{.3in}
\subfigure[$k=0.1$]{
ÊÊÊÊÊÊ\label{fig-4g}
ÊÊÊÊÊÊ\includegraphics[width=.48\textwidth]{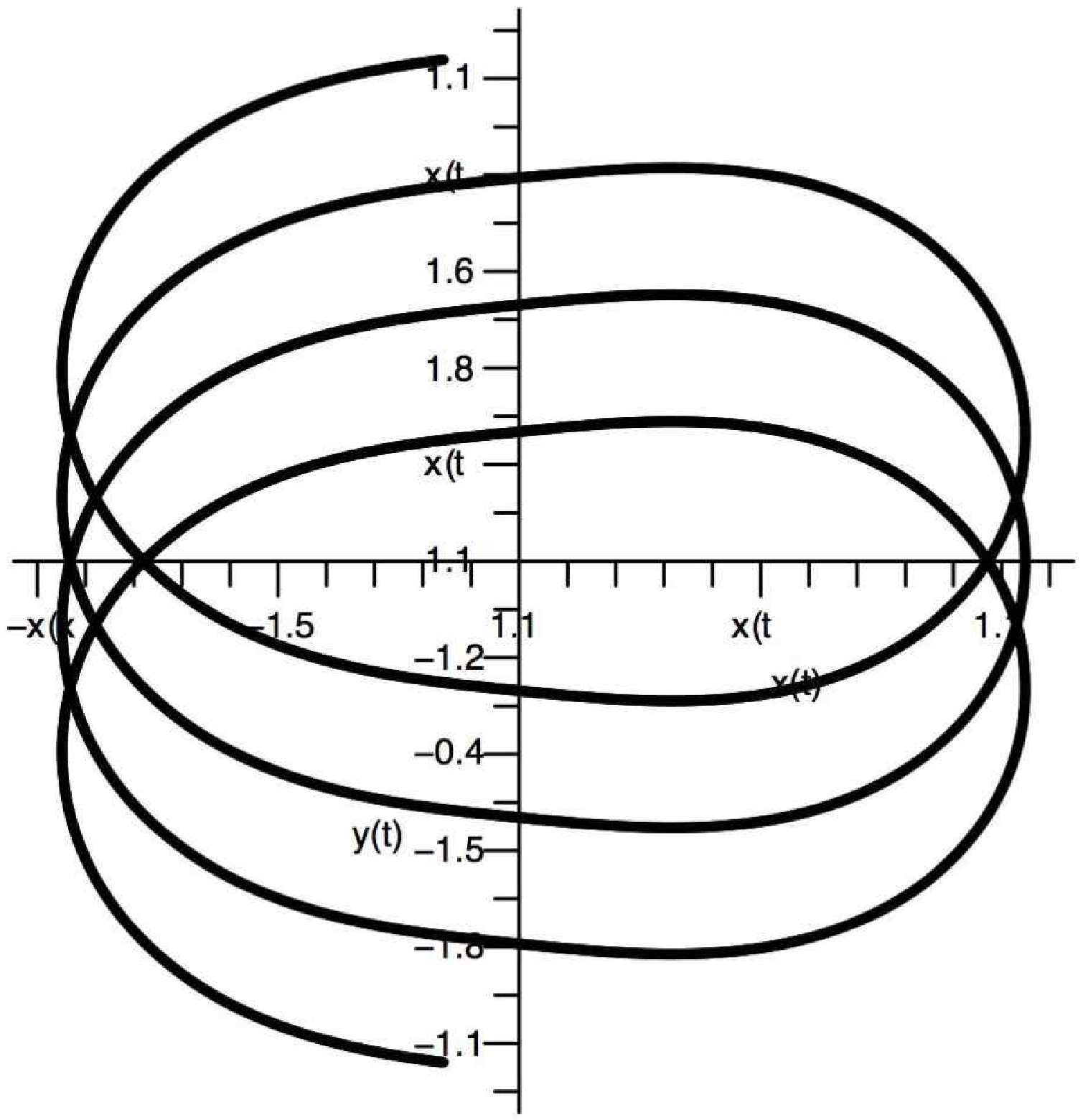}}
\subfigure[$k=0$]{
\label{fig-4h}
\includegraphics[width=.48\textwidth]{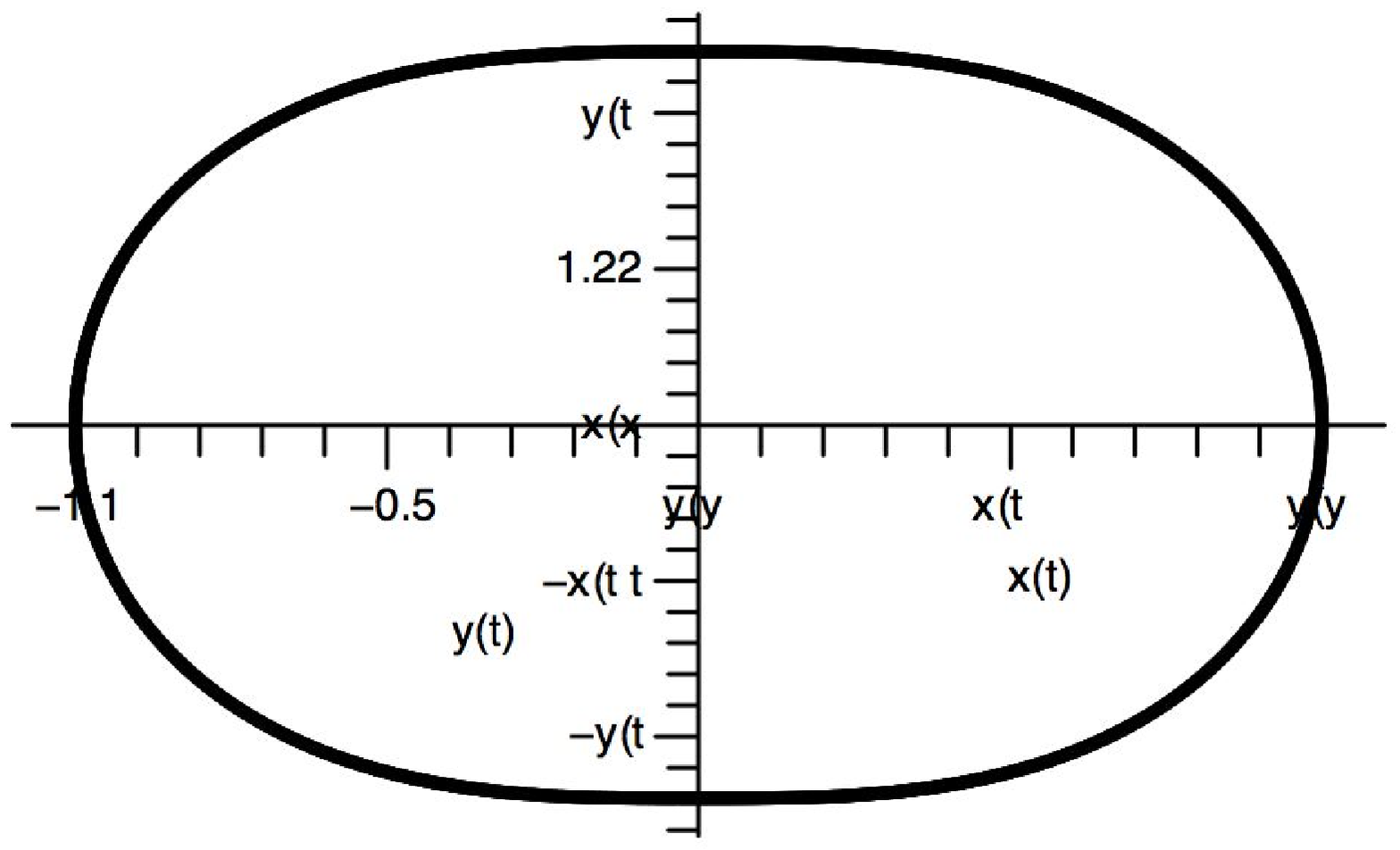}}

\caption{Drift is up and movement clockwise. As $k$ decays from $1$ to $0$  drift up slows down; \ref{fig-4h} corresponds to periodic case; we skip $k\ge 1$ which repeat those at figure \ref{fig-2}a-d, $k\le -1$ which will be mirror-symmetric to them with movement down and also clockwise and $-1<k<0$ as symmetric to given there.}
\label{fig-4}
\end{figure}

In the similar way movement in $(x_1,\xi_1)$ is described by Hamiltonian with one-well potential and if $k\ge 1$ or $k\le -1$ this particle moves only in $\{x_1\ge 0\}$ or $\{x_1\le 0\}$ respectively and in the former case the evolution coincides with given in the previous subsection while in the latter it one needs to replace $k$ by $-k$, take a mirror-symmetric picture but now the directions ``up'' and ``down'' are switched while direction of rotation does not change. So outer zone is described exactly as we did (see figures \ref{fig-2}a-d) but in the inner zone for $-1<k<1$ we get a kind of combined picture (see figure \ref{fig-4}).

Then
\begin{equation*}
T(k)={\frac 1 2}\bigl(T_+(k)+T_+(-k)\bigr),\qquad
I(k)={\frac 1 2}\bigl(I_+(k)-I_+ (-k)\bigr)
\end{equation*}
where $T_+(k)$, $I_+(k)$ are exactly what we introduced for Hamiltonian (\ref{1-2}). However due to monotonicity of $|y|^\nu\sign(y)$ we can apply idea of proposition \ref{prop-1-2}: introducing $z=|y|^\nu\sign(y)-k$ we arrive to
\begin{align}
T(k)=\hphantom{-} 2&\int_{-1}^1 |(k+z)\nu|^{{\frac 1 \nu}-1} (1-z^2)^{-{\frac 1 2}}\,dz, \label{1-15}\\
I(k)=-2&\int_{-1}^1 z|(k+z)\nu|^{{\frac 1 \nu}-1} (1-z^2)^{-{\frac 1 2}}\,dz=
\label{1-16}\\
&2\int_0^1 z\Bigl( |(k-z)\nu|^{{\frac 1 \nu}-1}-|(k+z)\nu|^{{\frac 1 \nu}-1}\Bigr) (1-z^2)^{-{\frac 1 2}}\,dz\notag
\end{align}
and obviously $\epsilon_0  I(k)/k\le c_0$ as $k\ne 0$. 

Thus we arrive to 

\begin{proposition}\label{prop-1-8} Consider on level $0$ trajectories of Hamiltonian $(\ref{1-14})$ with $\xi_2=k\in (-1,1)$. Along them

\smallskip
\noindent
(i) Variables $(x_1,\xi_1)$ are $T(k)$-periodic and $x_1$ oscillates between $b_1=-((1-k)\nu)^{1/\nu}$ and $b_2=((1+k)\nu)^{1/\nu}$;

\smallskip
\noindent
(ii) $x_2(t)={\tilde x}_2(t)+v(k)t$ with $T$-periodic ${\tilde x}_2(t)$ and $v(k)=I(k)/T(k)$;

\smallskip
\noindent
(iii) $\epsilon_0  I(k)/k\le c_0$ as $k\ne 0$;

\smallskip
\noindent
(iv) Therefore the only periodic trajectory is with $k=0$;

\smallskip
\noindent
(v) Proposition \ref{prop-1-6} holds with $k^*=0$ but $Z$ is even and $\alpha,\beta$ are odd with respect to $\xi_1$;

\smallskip
\noindent
(vi) For odd integer $\nu\ge 3$, $Z$, $\alpha$, $\beta$ are analytic.
\end{proposition}

\subsection{Pilot-models perturbed}
The periodic trajectories of our pilot-models are very fragile and one can destroy them easily. Consider \emph{heuristically\/} the same Hamiltonians with $V=\alpha x_1$ with $\alpha >0$ instead of $V=1$. Then 
${\frac {d\xi_2}{dt}}={\frac 1 2}\partial_{x_2}V= {\frac 1 2}\alpha$ and thus
$\xi_2 =c_0 + {\frac 1 2}\alpha t$. Then the averaged movement along $y=x_2$ is described by 
\begin{equation*}
{\frac {d\ }{dt}}y = V^{\frac 1 2} v\bigl(\xi_2(t)V^{-{\frac 1 2}}\bigr)\sim
(\alpha y)^{\frac 1 2}
v\bigl( {\frac 1 2} \alpha t (\alpha y)^{-{\frac 1 2}}\bigr) \iff 
{\frac {d\ }{dt}}z\sim v(tz^{-1})
\end{equation*}
with $z=2\sqrt {y/\alpha t}$. This ``equation'' has a solution $z=\beta_\pm^{-1}  t$ as $\pm t>0$ as
$ \beta_\pm v(\beta_\pm )=1$ and $\pm \beta_\pm >0$. 

\begin{figure}[t]
\centering
\subfigure[First pilot-model perturbed]{
ÊÊÊÊÊ\label{fig-5a}
ÊÊÊÊÊ\includegraphics[width=.48\textwidth]{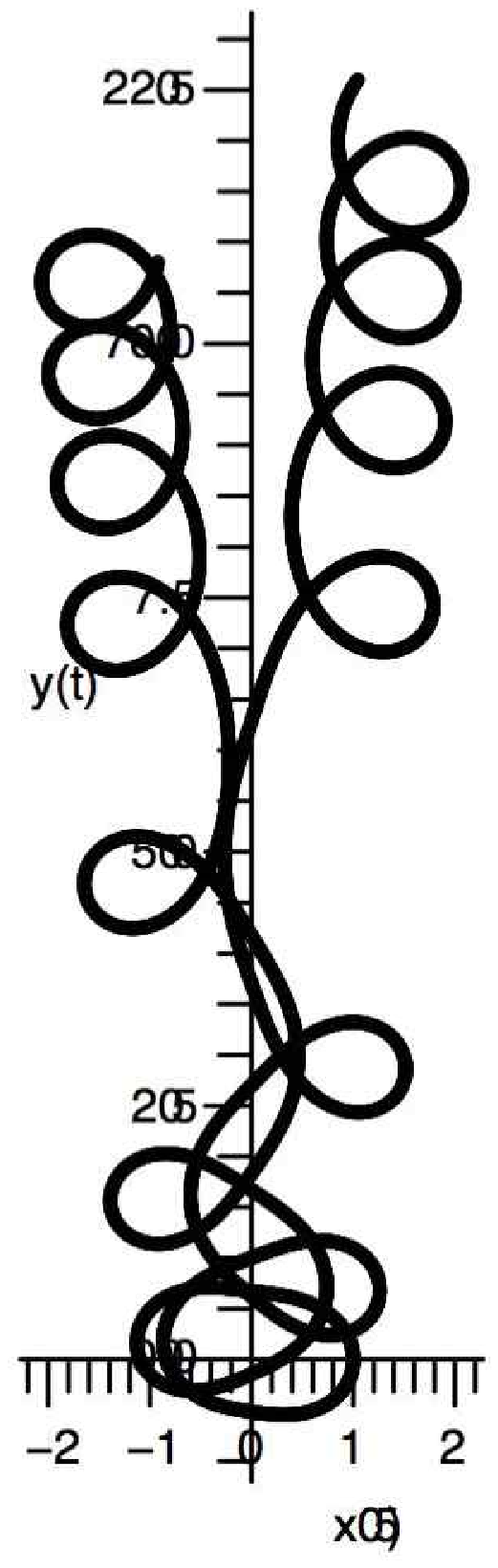}}
\subfigure[Second pilot-model perturbed]{
ÊÊÊÊÊ\label{fig-5b}
ÊÊÊÊÊÊ\includegraphics[width=.49\textwidth]{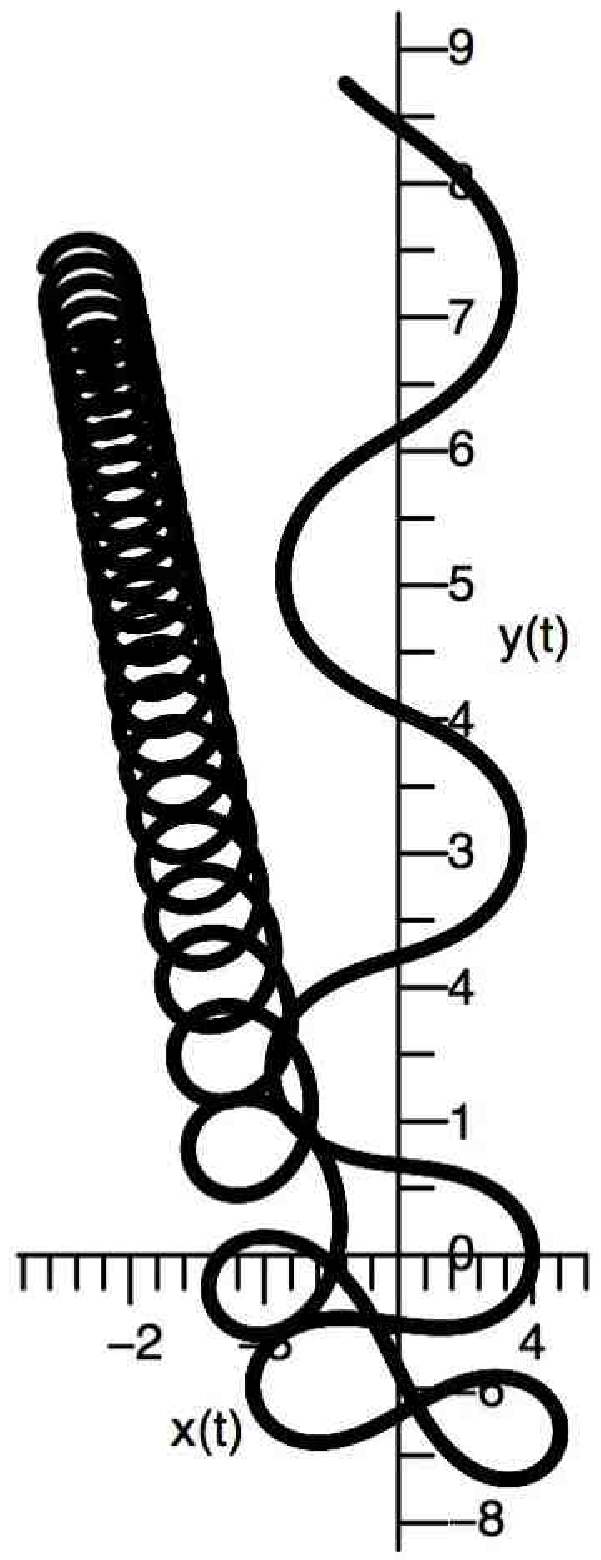}}\\

\caption{Both pilot-models are perturbed by a weak electric field  directed up
($V=\alpha x_2$)}
\label{fig-5}
\end{figure}

In the second pilot-model this equation  has solutions with $-\beta_-=\beta_+>1$ and therefore  on them ``in average'' $x_1 \sim (\nu \xi_2V^{-{\frac 1 2}})^{\frac 1 \nu}V^{\frac 1 {2\nu}}\sim \rho y^{\frac 1 \nu}$. In the first pilot-model $\beta_-=-1$. These observations 
explain the evolutions in vertical direction on figures \ref{fig-5}a-b.

\subsection{General case. Remarks}

We are going to consider
\begin{equation}
a={\frac 1 2}\Bigl(\xi_1^2 + \sigma(x)^2 \bigl(\xi_2-\mu \phi (x) \varrho_\nu(x_1)/\nu \bigr)^2 -V(x)\Bigr)
\label{1-17}
\end{equation}
with coefficients
\begin{equation}
\sigma \ge \epsilon_0,\quad\phi\ge \epsilon_0,\quad V\ge \epsilon_0
\label{1-18}
\end{equation}
and all of them are smooth and with $\varrho_\nu=|x_1|^\nu$ or $\varrho_\nu=|x_1|^\nu\sign x_1$, $\nu\ge 2$.

Let
\begin{equation}
p_1=\xi_1,\quad p_2=\xi_2-\mu \phi (x) \varrho_\nu(x_1)/\nu;
\label{1-19}
\end{equation}
then
\begin{equation}
\{p_1,p_2\}= -\mu \partial_{x_1}\bigl(\phi (x) \varrho_\nu(x_1)/\nu\bigr)=
-\mu \phi(x)\varrho_{\nu -2}(x_1)x_1 \bigl(1+O(x_1)\bigr).
\label{1-20}
\end{equation}

We want to consider trajectories on level $\{a=0\}$ (and thus on any level $\{a=\tau\}$) and thus we can multiply $a$ by any smooth factor disjoint from 0; 
then we can assume that 
\begin{equation}
\sigma |_\Sigma =\phi|_\Sigma =1.
\label{1-21}
\end{equation}
Really, starting from (\ref{1-17}) we can multiply $a$ by smooth function $\omega^2$ and redefine $x_1$ according to the new metrics; then we redefine $x_2$ to get rid of the mixed term. Then (\ref{1-17}) will still hold but with
$V\mapsto \omega^2V$, $\sigma |_\Sigma\mapsto (\sigma \omega^2)|_\Sigma$, 
$\phi|_\Sigma \mapsto (\phi \omega ^\nu)|_\Sigma$ and picking up 
$\omega = (\sigma \phi)^{-1/(\nu +1)}$ we arrive to the case
$(\sigma \phi)|_\Sigma =1$ with
\begin{equation}
V\mapsto V^*= (\sigma \phi)^{-2/(\nu +1)}V.
\label{1-22}
\end{equation}
Finally, replacing $x_2$ by $\rho(x_2)$ we can achieve $\sigma |_\Sigma=1$ and thus (\ref{1-21}).

\subsection{General case. Outer zone}

\begin{proposition}\label{prop-1-9} For general Hamiltonian $(\ref{1-17})$  in $B(0,1)$ satisfying $(\ref{1-18})$

\smallskip
\noindent
(i) There exists $b=b(x ,p_1,p_2)$, such that as $|a|\le \epsilon_0$ and 
$|f|\ge {\bar\gamma}=C\mu^{-1/\nu}$
\begin{align}
&C_0^{-1} \le {\frac b f}\le C_0,\label{1-23}\\
&|\{a,b\}|\le C_0\mu^{-1}|f|^{-1} \label{1-24};
\end{align}

\smallskip
\noindent
(ii) Trajectories starting at $t=0$ from 
$(\bx,\bxi) \in B(0,{\frac 1 2})\cap  \{x, \pm f(x)=\gamma \ge {\bar\gamma}\}$ for time 
$|t|\le T_1 = \epsilon _0\mu f({\bar x})^2$ remain in $B(0,1)$ and along them
$\epsilon \le \pm f(x(t))/\gamma \le C$;

\smallskip
\noindent
(iii) Trajectories starting at $t=0$ from 
$(\bx,\bxi) \in B(0,{\frac 1 2})\cap  \{x, |f(x)|\le \gamma = {\bar\gamma}\}$ for time $|t|\le T_1 = \epsilon $ remain in $B(0,1)$ and along them
$|f(x(t))| \le C_1 {\bar\gamma}$.
\end{proposition}

\begin{proof}
Let us consider trajectory starting at 
$(\bx,\bxi)$ with ${\bar x}\in B(0,{\frac 1 2})$, $a\le c$, 
$\gamma = |{\bar x}_1|\ge {\bar\gamma}=C_0\mu^{-1/\nu}$.  Consider this trajectory as long as it is contained in $B(0,1)$ but not further than 
$T_1=\epsilon \mu \gamma^\nu$\,\footnote{\label{foot-3}Then 
$T_1\asymp 1$ as $|{\bar x}_1|\asymp {\bar\gamma}$ and $T_1\asymp \mu$ as $|{\bar x}_1|\asymp 1$.}

Along this trajectory
\begin{equation}
{\frac {d\xi_2} {dt}}=
\{a,\xi_2\}=\sigma^2(x)p_2 \phi_{x_2} (x)\times \mu \varrho_\nu(x_1)/\nu + O(1).
\label{1-25}
\end{equation}
Then for $X=\xi_2-\beta (x) p_1$
\begin{equation}
{\frac {dX} {dt}}=\{a,X\}=\sigma^2 p_2 \mu \varrho_\nu (x_1) 
\bigl({\frac 1 \nu}\phi_{x_2}(x)  -\phi(x) \beta x_1^{-1} \bigr) + O(1)
\label{1-26}
\end{equation}
(as $a=0$)  and therefore ${\frac {d\ }{dt}}X=O(1)$ 
for \begin{equation}
X=\xi_2-\beta (x) p_1,\qquad \beta =  {\frac 1 \nu} x_1\phi_{x_2}\phi ^{-1}=
{\frac 1 \nu} x_1\partial_{x_2}\log |\phi(x)|.
\label{1-27}
\end{equation}
Therefore $|X(x)-X({\bar x})|\le CT\le C\epsilon \mu\gamma ^\nu$,  which together with inequality $|p_2|\le c$ imply that 
$|\xi_2-{\bar\xi}_2 |\le C\epsilon \mu \gamma^\nu $;
 for   small enough constant $\epsilon>0$ we also get
$\epsilon_1\mu \gamma^\nu \le |\xi_2|\le  c_1\mu \gamma^\nu $ 
(I remind that $\mu \gamma^\nu \ge C$) and thus
$\epsilon_1^{-1}\gamma \le |x_1|\le c_1\gamma$. 

So, trajectory remains  in the strip described above and since the drift speed according to \cite{IRO3} does not exceed 
$C_0(\mu \gamma^\nu )^{-1}$ there we are insured that trajectory remains in $B(0,1)$ for time $T$. This proves (ii).

On the other hand, if trajectory starts in $(\bx,\bxi)$ with 
$|{\bar x}_1|\le {\bar\gamma}$, it remains in $B(0,1)$ for time $T=\epsilon$ since the speed does not exceed $C_0$ and due to above arguments it cannot get into zone $\{x_1\ge C{\bar\gamma}\}$. This proves (iii).

To prove (i) consider $b=(\mu^{-1}X)^{1/\nu}$; then as
$|x_1|\ge {\bar\gamma}$ obviously $b\asymp |x_1|$ and furthermore
$|\{a,b\}|\le C(\mu^{-1}|\xi_2|)^{1/ \nu} |\xi_2|^{-1}\le
c\mu^{-1}|x_1|^{1-\nu}$ which is exactly (\ref{1-24}).
\end{proof}

Obviously ${\frac {d\ }{dt}}x_2=\sigma^2 p_2$ and for 
$Y=x_2-\mu^{-1}\varrho_{-\nu}(x_1)x_1\beta'p_1$ we have
\begin{equation*}
{\frac {dY}{dt}}=\sigma^2\bigl(1-\beta'\phi \bigr)p_2 +
 (\nu -1)\beta'\mu^{-1}\varrho_{-\nu}(x_1)p_1^2 +
O(\mu^{-1}|x_1|^{1-\nu}).
\end{equation*}
So, let $\beta'=\phi^{-1}$; then 
\begin{equation*}
{\frac {dY}{dt}}=
{\frac 1 2}(\nu -1)\beta' \mu^{-1}\varrho_{-\nu}(x_1) (a -V) + 
{\frac 1 2}(\nu -1)\beta' \mu^{-1}\varrho_{-\nu}(x_1)(p_1^2-\sigma^2p_2^2)+
O(\mu^{-1}|x_1|^{1-\nu})
\end{equation*}
and redefining 
\begin{equation}
Y= x_2-\mu^{-1}\varrho_{-\nu}(x_1) x_1\beta'p_1 - \mu^{-2}\varrho_{-2\nu}(x_1) x_1p_1p_2
\label{1-28}
\end{equation}
we arrive to 
\begin{equation}
{\frac {dY}{dt}}=-{\frac 1 2}(\nu -1)\beta' \mu^{-1}\varrho_{-\nu} (x_1)+ O\bigl(\mu^{-2}|x_1|^{-2\nu}+\mu^{-1}|x_1|^{1-\nu}\bigr),
\label{1-29}
\end{equation}
and we arrive to
\begin{proposition}\label{prop-1-10}
Along trajectories in outer zone $\{\epsilon \ge |x_1|\ge {\bar\gamma}\}$ 
\begin{equation}
\epsilon_0 (\mu \gamma ^\nu)^{-1}|t|\le \pm \varrho_0({\bar x}) \bigl(Y(t)-Y(0)\bigr)\le 
C (\mu \gamma ^\nu)^{-1}|t|
\label{1-30}
\end{equation}
as $0\le \pm t \le \epsilon \mu\gamma^\nu$.
\end{proposition}

\subsection{General case. Inner zone}

Now we need to consider trajectories lying in inner zone 
$\cZ_\inn=\{x:|x_1|\le C{\bar\gamma}\}$. Here function $X$ is still defined but  $Y$ is not and according to subsection 1.1 the important role is played by $X(x,p)V(0,x_2)^{-{\frac 1 2}}$.

\begin{proposition}\label{prop-1-11}
Let conditions  $(\ref{1-18})$, $(\ref{1-21})$,  $|\xi_2(0)|\le C$
be fulfilled and let condition
\begin{equation}
\varsigma\bigl( \xi_2-k^*|V|^{\frac 1 2}\bigr)\ge \epsilon _1, \qquad \varsigma=\pm 1
\label{1-31}
\end{equation}
be satisfied as $t=0$. Then 

\smallskip
\noindent
(i) For $|t|\le T=\epsilon$ this condition  is satisfied (with  $\epsilon_1$ replaced by $\epsilon_1/2$);

\smallskip
\noindent
(ii) Moreover
\begin{equation}
\epsilon_2 |t|\le \pm \varsigma \bigl(x_2(t)-x_2(0)\bigr)\le c_2 |t|
\qquad\qquad {\text as}\quad
c{\bar\gamma }\le \pm t\le T
\label{1-32}
\end{equation}
and 
\begin{equation}
\epsilon_2|t|\le |x(t)-x(0)|+|p(t)-p(0)|\le c_2|t|\qquad{\text as}\quad |t|\le c{\bar\gamma}.
\label{1-33}
\end{equation}
\end{proposition}

\begin{proof} Proof of (i) is obvious; (\ref{1-32}) and (\ref{1-33}) then follow from comparison of the general and pilot-model system. 
\end{proof}

So there are no $T$-periodic trajectories with $0<T\le \epsilon$ unless
\begin{equation}
\rho=| \xi_2-k^* V ^{\frac 1 2}| \le \epsilon _1.
\label{1-34}
\end{equation}
Note that now
\begin{equation}
{\frac {d\xi_2}{dt}}= {\frac 1 2}V_{x_2}+O({\bar\gamma})
\label{1-35}
\end{equation}
which immediately yields

\begin{proposition}\label{prop-1-12} Let conditions $(\ref{1-18})$, $(\ref{1-21})$, $(\ref{1-34})$ be fulfilled and
\begin{equation}
\varsigma \partial_{x_2}V \ge \epsilon_0\qquad \text{with\quad} \varsigma =\pm 1
\label{1-36}
\end{equation}
Then trajectories on level $0$ satisfy
\begin{equation}
\pm \varsigma\bigl(\xi_2(t)-\xi_2(0)\bigr)\ge \epsilon |t|\qquad  {\text as}\quad 0\le \pm  T\le \epsilon.
\label{1-37}
\end{equation}
\end{proposition}

In what follows we will consider more precisely trajectories satisfying (\ref{1-34}) and also $|\partial_{x_2}V |\le \epsilon_0$.

We still want to consider such trajectories as long as they are in a certain vicinity of the original point $x(0)$. For the pilot-models this vicinity is $B(0,1)$
and $T=\epsilon \rho^{-1}$ as $|\xi_2 - k^*|=\rho $ but for other Hamiltonians
both vicinity and time are smaller.

\begin{proposition}\label{prop-1-13} For any of two pilot-models with coupling constant $\mu \ge 1$ and potential $V=1\pm \zeta t$, $0\le\zeta\le \epsilon$ 

\smallskip
\noindent
(i) Trajectory starting at $(\bx,\bxi)$ satisfying $(\ref{1-34})$  remains in $B(0,1)$ and satisfies  $(\ref{1-34})$ with $2\epsilon_1$ as long
$|t|\le T= \epsilon \min\bigl(\zeta ^{-{\frac 1 2}}, \rho^{-1}\bigr)$;

\smallskip
\noindent 
(ii) Along it 
\begin{equation}
|x_2(t)-x_2(0)|\le C(\zeta |t|^2 + \rho |t|)+C{\bar \gamma},\qquad 
\xi_2(t) -\xi_2(0)= \pm {\frac 1 2}\zeta t
\label{1-38}
\end{equation}
and in at least one time direction 
\begin{equation}
|x_2(t)-x_2(0)|\ge \epsilon (\zeta |t|^2 + \rho |t|)-C{\bar \gamma};
\label{1-39}
\end{equation}
\end{proposition}

\begin{proposition}\label{prop-1-14} For general Hamiltonian $(\ref{1-17})$, satisfying $(\ref{1-18})$, $(\ref{1-21})$ and
\begin{equation}
\bigl|\partial_{x_2}V|_\Sigma \bigr|\le \zeta
\label{1-40}
\end{equation}

\smallskip
\noindent  
(i) Statement (i) of proposition \ref{prop-1-13} holds as 
\begin{equation}
\rho \ge {\bar\gamma}, \qquad \zeta\ge {\bar\gamma} 
\label{1-41}
\end{equation}
and $|\xi_2- k^*V^{\frac 1 2}|\le \rho$;

\smallskip
\noindent 
(ii) Assuming $(\ref{1-41})$, if
\begin{equation}
\bigl|\partial_{x_2}V|_\Sigma \bigr|\ge \epsilon_1\zeta
\label{1-42}
\end{equation}
then 
\begin{equation}
|\xi_2 (t)-\xi_2(0)|\ge \epsilon \zeta |t|;
\label{1-43}
\end{equation}

\smallskip
\noindent 
(iii) Assuming $(\ref{1-41})$, if
\begin{equation}
\bigl| \xi_2 (0)-k^* V(x(0))^{\frac 1 2} \bigr|\ge \epsilon_1\rho
\label{1-44}
\end{equation}
then in at least one time direction
\begin{equation}
|x_2 (t)-x_2(0)|\ge \epsilon \rho|t|.
\label{1-45}
\end{equation}
\end{proposition}

\begin{proof} Proofs are obvious.
\end{proof}

\subsection{Hamiltonian maps}

Now we need to analyze more precisely what happens in the periodic zone
$\cZ_\per=\{x, |x_1|\le c{\bar\gamma}, |\xi_2 - k^*V^{\frac 1 2}|\le \rho\}$
where $\rho$ will be specified in the next section.
Let us consider pilot-models first.

If on each energy level $\tau$, $|\tau|\le \epsilon$ all the trajectory were periodic with period $T=T(\tau)$ then replacing Hamiltonian $a$ by $g(a)$ with
$g(\tau)= \int ^\tau T(\tau)\,d\tau$ we would get 1-periodic Hamiltonian flow. 

However it is not the case: on each energy level $\tau$ the only periodic trajectory is one with $\xi_2= \xi_2(\tau)=k^* (1+2a\tau)^{\frac 1 2}$ and period is $T(\tau)=T^*(1+2\tau)^{(\nu -1)/2\nu}$. Still defining $g(\tau)$ accordingly let us consider the Hamiltonian $g(a)$.

\begin{proposition}\label{prop-1-15} For a pilot-model symbol $a$ given by $(\ref{1-2})$ or $(\ref{1-14})$
\begin{align}
&e^{T^*H_{g(a)}}=e^{H_b}, \label{1-46}\\
&g(\tau)={\frac 1 {2\kappa}}(1+2\tau )^\kappa - {\frac 1 {2\kappa}},\qquad \kappa={\frac {\nu +1}{2\nu}},
\label{1-47}\\
&b={\bar b}(\xi_2,a)=(1+2a)^\kappa \Bigl(\frac {\xi_2} {(1+2a)^{\frac 1 2}}-k^*\Bigr)^2 
\omega \Bigl(\frac {\xi_2} {(1+2a)^{\frac 1 2}}\Bigr)
\label{1-48}
\end{align}
where $T^*$ is an elementary period and  
\begin{equation}
\omega(k^*)= {\frac 1 2} \partial_k I(k)\bigr|_{k=k^*}.
\label{1-49}
\end{equation}
\end{proposition}

\begin{proof} As $a$ is a pilot-model Hamiltonian then taking $\xi_2=k^*$, $a=0$ we would conclude that the Hamiltonian trajectories  are periodic with period $T^*$. As $a=\tau$ disjoint from $-{\frac 1 2}$ we have both period $T(\tau)$ and $K^*(\tau)$ defined by
\begin{equation}
T (\tau)= T^* \cdot (1+2\tau)^{(\nu -1)/(2\nu)},\qquad 
K^*(\tau)= k^* \cdot (1+2\tau)^{\frac 1 2}.
\label{1-50}
\end{equation}
To get period which is independent on energy level one must replace $a$ by $g(a)$  with $g(\tau)=\int ^\tau T (\tau)\,d\tau$: then $e^{T^*H_{g(a)}}=I$ as $\xi_2= K^*(a)$;
period is $T^*$ rather than 1  since we define $g(a)$ without factor $T^*$.  

As $\xi_2\ne K^*(a)$ periodicity is broken. Since $H_a$ and $H_{\xi_2}$ near $\{a=0\}$ are linearly independent we conclude that
$e^{T^*H_{g(a)}}=I+2(\xi_2-K^*(a)) \omega _1(x_1,\xi_1,\xi_2)$. However since it is a symplectic map it must be of the form $e^{H_b}$ with 
\begin{equation*}
b= (\xi_2-K^*(a))^2 \omega _2(x_1,\xi_1,\xi_2)+ \lambda (x_1,\xi_1).
\end{equation*}
We know that as $\xi_2=K^*(a)$ we have $e^{H_b}=e^{H_\lambda}=I$ and since $\lambda$ is ``small'' we conclude that $\lambda=\const$ (and thus we can take it 0). 

Since $e^{H_b}$ commutes with $e^{H_a}$ we conclude that $\omega=\omega(\xi_2,a)$ and due to homogeneity properties 
$\omega_2 = \omega \bigl(\xi_2(1+2a)^{-{\frac 1 2}}\bigr)(1+2a)^{\kappa-1}$.
Then $x_2$-shift as $a=0$ and $\xi_2\approx k^*$ is 
$2\omega (k^*)(\xi_2-k^*)+O\bigl( (\xi_2-k^*)^2\bigr)$ which implies (\ref{1-50}).
\end{proof}

If we consider pilot-model but with $\mu\ne 1$ one needs to replace (\ref{1-46})
by
\begin{equation}
e^{{\bar\gamma}T^*g(a)}=e^{{\bar\gamma}H_b}.
\label{1-51}
\end{equation}
Consider now trajectories  residing in $B(0,C_0{\bar\gamma})$ and corresponding to these Hamiltonians but with potential $V=1+\zeta x_2$. One can see easily that (\ref{1-51}) will remain true with $b={\bar b}- {\frac 1 2}T^*W + b'$, $W=V(0,x_2)$  and $b'$satisfying (\ref{1-53}) below. Then for the general Hamiltonians one can prove

\begin{proposition}\label{prop-1-16} For general symbol $(\ref{1-17})$ with $\sigma, \phi$ satisfying $(\ref{1-21})$ and with $V$, $V(0)=1$, satisfying in $C_0{\bar\gamma}$-vicinity of \ $0$
\begin{equation}
|\partial _{x_2} V|\le \zeta 
\label{1-52}
\end{equation}
with $\zeta \ge C{\bar\gamma}$, equality $(\ref{1-51})$ holds with symbol 
$b={\bar b}- {\frac 1 2}T^*W(x_2) + b'$
where ${\bar b}={\bar b}(\xi_2,a)$ is defined by $(\ref{1-48})$,  and 
\begin{equation}
|\partial_x ^\alpha  b' |\le 
C \bigl(\zeta (\zeta+\rho)+{\bar\gamma}\bigr){\bar\gamma}^{1-|\alpha|}, \qquad |\alpha|\le 1.
\label{1-53}
\end{equation}
\end{proposition}

\sect{Quantum dynamics}

\subsection{Preliminary notes. Forbidden zone}

Starting from this section $\nu\ge2$ is an integer and magnetic potential is $(0, \phi (x) \varrho_\nu(x_1) /\nu)$ where  $\varrho_\nu(x_1) = x_1^\nu$.
Thus we have different pilot-models for even and odd $\nu$. 

In this subsection we consider $\mu \ge h^{-1}$ and \emph{forbidden\/} zone is defined by (\ref{2-1}). Also we consider $\xi_2$-localization which is useful 
when $\mu$ is close to $h^{-\nu}$.

\begin{proposition}\label{prop-2-1} Let $\mu \ge h^{-1}$ and 
\begin{equation}
\mu h f -V \ge \epsilon
\label{2-1}
\end{equation}
in $B({\bar x}, \gamma)$ with
\begin{equation}
\gamma \ge \max \bigl(h ,{\bar\gamma}\bigr),\qquad 
{\bar\gamma}=C\mu^{-{\frac 1 \nu}}.
\label{2-2}
\end{equation}
Then 
\begin{equation}
|e(x,x,0)|\le C\gamma^{-2}\bigr(\mu ^{-1}h\gamma^{-\nu -1}\bigr)^s\qquad 
{\text in}\quad B({\bar x}, \gamma).
\label{2-3}
\end{equation}
\end{proposition}

\emph{Proof\/} is just by rescaling. In particular, as 
\begin{equation}
h^{-1}\le \mu \le C_0h^{-\nu}
\label{2-4}
\end{equation}
the total contribution  to 
$\int e(x,x,0)\psi (x)\,dx$ of the \emph{forbidden zone\/} defined by (\ref{2-1}) does not exceed  $C{\bar\gamma}_1^{-1}(\mu h^\nu)^s$. Thus we need to consider only \emph{allowed zone\/} $\mu h f-V\le \epsilon$ and there
\begin{equation}
\gamma \le {\bar\gamma}_1=C_0(\mu h)^{-{\frac 1 {\nu -1}}}.
\label{2-5}
\end{equation}

This result is not very useful as 
$\mu\ge h^{-\nu+\delta}$. However using logarithmic uncertainty principles and related arguments one can prove easily

\begin{proposition}\label{prop-2-2}
Let conditions $(\ref{2-1})$, $(\ref{2-4})$ be fulfilled and
\begin{equation}
\gamma \ge \max \bigl(Ch|\log h| ,{\bar\gamma}\bigr).
\label{2-6}
\end{equation}
with $C=C_s$. Then 
\begin{equation}
|F_{t\to h^{-1}\tau} {\bar\chi}_T(t)\psi u|\le CT h^s
\qquad \forall \tau\le \epsilon\quad \forall T\ge Ch|\log h|
\label{2-7}
\end{equation}
and therefore
$e(x,x,0)\le Ch^s$ in $B({\bar x},\gamma)$.
\end{proposition}

In particular, (\ref{2-6}) is fulfilled automatically as 
$|x_1|\ge {\bar\gamma}_1$ provided 
\begin{equation}
\mu \le \epsilon (h|\log h|)^{-\nu}
\label{2-8}
\end{equation}
while the \emph{super-strong magnetic field\/} case when (\ref{2-8}) is violated requires some special consideration based on the analysis operators with operator-valued symbols. We need axillary

\begin{proposition}\label{prop-2-3} Let $\lambda=\lambda_n(z)$ be eigenvalues of one-dimensional Schr\"odinger operator 
\begin{equation}
{\mathbf a}(z)^0=D_1^2 + (z-\mu \varrho_\nu(x_1)\nu /\nu)^2
\label{2-9}
\end{equation}
in $L^2(\bR)$; here we assume only that $2\le \nu \in \bR$ and either
$\varrho_\nu(x_1)= |x_1|^\nu$ or $\varrho_\nu(x_1)= |x_1|^\nu\cdot\sign(x_1)$. 
Then

\smallskip
\noindent
(i) $\lambda_n(z)> 0$;

\smallskip
\noindent
(ii) Assume that either $z\ge 0$ or $\varrho_\nu = |x_1|^\nu\sign (x_1)$. Then for $n\le \epsilon |z|^{\frac {\nu +1} \nu}$
\begin{align}
&\lambda_n \asymp |z|^{\frac {nu -1} \nu } n\label{2-10},\\
&{\frac {\partial\ }{\partial z}}\lambda_n \cdot \sign (z) \asymp 
|z|^{-{\frac 1 \nu }} n;\label{2-11}
\end{align}

\smallskip
\noindent
(iii) Assume that $z<0$ and $\varrho_\nu = |x_1|^\nu$. Then 
\begin{equation}
\lambda_n \asymp z^2+n^{\frac {2\nu}{\nu +1}}.
\label{2-12}
\end{equation}
\end{proposition}

\begin{proof} Proof of (i) is obvious; (ii),(iii) easily follow from the semiclassical character of the spectrum as $|z|+n\gg 1$, $n\gg 1$ respectively. Easy details are left to the reader.
\end{proof}

Furthermore, using operators with operator-valued symbols one can prove easily

\begin{proposition}\label{prop-2-4}  Let $\mu \ge Ch^{-\nu}$. Then
\begin{equation}
 |e(x,x,\tau)| \le C\mu ^{-s}\qquad \for 
 \tau \le \epsilon \bigl(\mu h^\nu\bigr)^{\frac 2 {\nu+1}}.
 \label{2-13}
\end{equation}
\end{proposition}

So, \emph{in what follows  can assume that $\mu \le Ch^{-\nu}$}. In what follows
$u=u(x,y,t)$
is the \emph{propagator\/}, i.e. the Schwartz kernel of operator $e^{ih^{-1}tA}$.

\begin{proposition}\label{prop-2-5} Let $\mu \le Ch^{-\nu}$ and 
let $\varphi \in C^\infty(\bR)$, $\varphi (t)=1$ as $t\ge 2$, $\varphi(t)=0$ as $t\le 1$. Then 

\smallskip
\noindent
(i) The following inequality holds
\begin{multline}
|F_{t\to h^{-1}\tau} {\bar\chi}_T(t) 
\varphi  \bigl(\pm\epsilon (\mu h^\nu)^{\frac 1 {\nu -1}}hD_{x_2} \bigr)\psi u| \le C(T+1)h^s\qquad \forall \tau\le \epsilon\quad \forall T\ge Ch|\log h|
 \label{2-14}
\end{multline}
and therefore
\begin{equation}
|\varphi  \bigl(\pm\epsilon (\mu h^\nu)^{\frac 1 {\nu -1}}hD_2 \bigr)e(x,y,0)| \le Ch^s;
 \label{2-15}
 \end{equation}
 here and below $x,y\in B(0,1-\epsilon)$;
 
\smallskip
\noindent
(iii) Furthermore, for even $\nu$ 
 \begin{equation}
|F_{t\to h^{-1}\tau} {\bar\chi}_T(t) 
\varphi  \bigl(-\epsilon  hD_2 \bigr)\psi u| \le C(T+1)h^s\qquad
\forall \tau\le \epsilon\quad \forall T\ge Ch|\log h|
 \label{2-16}
 \end{equation}
and therefore 
 \begin{equation}
|\varphi   \bigl(-\epsilon hD_2 \bigr)e(x,y,0)| \le Ch^s.
 \label{2-17}
 \end{equation}
\end{proposition}

Easy proofs are left to the reader. Inequalities (\ref{2-14}) and (\ref{2-16}) are proved by the standard elliptic methods. Then inequalities (\ref{2-15}) and (\ref{2-17}) are proved by the standard Tauberian methods.

Thus in what follows we localized $u$ with respect to  $hD_{x_2}$ (and thus with respect to $hD_{y_2}$ as well due to symmetry. This localization is more precise than $x_1,y_1$ localization as $\epsilon (h|\log h|)^{-\nu}\le \mu \le Ch^{-\nu}$; as $\mu\le \epsilon (h|\log h|)^{-\nu}$ both localizations are equivalent. 

\subsection{Outer zone}

In this subsection we consider \emph{outer zone\/} 
$\cZ_\out=\{C_0{\bar\gamma}\le|x_1|\le {\bar\gamma}_1\}$ with 
\begin{equation*}
{\bar\gamma}_1=\min \bigl(\epsilon, C(\mu h)^{-1/(\nu -1)}\bigr).
\end{equation*}
This definition works fine under (\ref{2-8})
but should be modified as 
\begin{equation}
\epsilon (h|\log h|)^{-\nu}\le \mu \le C_0h^{-\nu}
\label{2-18}
\end{equation}
to 
\begin{equation*}
\cZ_\out= \{ C_0 \le |\xi_2|\le C(\mu h^\nu)^{-\frac 1{\nu -1}}\}
\end{equation*}
(and  in addition $\{\xi_1 \ge C_0\}$ as $\nu$ is even)\footnote
{\label{foot-4}Restriction to $|x_1|$ will be $|x_1|\le Ch|\log h|$.}.

Our first statement is that in the quantum dynamic in this zone the magnitude of $x_1$ (and $\xi_2$) persists.

\begin{proposition}\label{prop-2-6} (i)  Let 
${\bar x}\in B(0,{\frac 1 2})\cap  \{x, \pm x_1\ge C{\bar\gamma}\}$,
$\gamma =\epsilon |{\bar x}_1|$ satisfy $(\ref{2-6})$ and 
$\psi\in C_0^\infty (B(\bx,\gamma))$ be a rescaling of the standard function. Let $\psi_1$ be $\gamma$-admissible and supported in 
$\{1/(2C_0)\le x_1/{\bar x}_1\le 2C_0\}$ and equal $1$ in 
$\{1/C_0\le x_1/{\bar x}_1 \le 2C_0\gamma\}$. Let $T_1=\epsilon \mu\gamma^\nu$. Then for $T=T_1$
\begin{equation}
|F_{t\to h^{-1}\tau} {\bar\chi}_T(t)(1-\psi_{1\,x})u\psi_y |\le Ch^s.
\label{2-19}
\end{equation}

\smallskip
\noindent
(ii) Let $ {\bar\xi}_2=\pm \rho$, $\rho \ge C$ and 
$\varphi\in C_0^\infty ([{\frac 3 4},{\frac 3 4}])$ and 
$\varphi_1\in C_0^\infty$, $\varphi_1=1$ on $[{\frac 1 2}, {\frac 3 2}]$ be  the standard functions. Let $T_1= \epsilon \rho$. Then for $T=T_1$
\begin{equation}
|F_{t\to h^{-1}\tau} {\bar\chi}_T(t)\Bigl
(1-\varphi _1\bigl({\frac 1 {{\bar \xi}_2}} hD_x\bigr)\Bigr)u
\varphi \bigl( {\frac 1 {{\bar \xi}_2}} hD_x\bigr) |\le Ch^s.
\label{2-20}
\end{equation}
\end{proposition}

\begin{proof} Note that both statements are equivalent as 
$\rho \asymp \mu \gamma^\nu$ so we will prove (ii). As $|x_1|\ge {\bar\gamma}$ we constructed by $(\ref{1-27})$ symbol which after rescaling becomes $X=\xi_2-\mu^{-1}x_1^{1-\nu}\beta (x) p_1$ such that 
\begin{equation}
|\{a,X\}-\alpha a |\le C_0
\label{2-21}
\end{equation}
Then using arguments of the proof of theorem 3.1 \cite{IRO1} with symbol
\begin{equation}
\varpi (\varsigma Ct \pm X )
\label{2-22}
\end{equation} 
where $\varpi$ is the same function $\chi$ as in this proof and 
$\varsigma=\pm 1$, one can easily prove (ii) as long as symbol (\ref{2-22}) is quantizable; this condition is equivalent to (\ref{2-6}).

Surely one needs to check that in time $T_1$ we stay in $B(0,1)$ but it will be done in the next proposition.

If (\ref{2-6}) is violated (and thus (\ref{2-8}) is violated as well) we can quantize with respect to $(x_2,\xi_2)$ but not with respect to $(x_1,\xi_1)$ unless $x_1\ge Ch|\log h|$ and instead we note that 
\begin{equation*}
\{a, \xi_2\}=\alpha a+\sigma ^2 \bigl(\xi_2 - \phi \varrho_\nu/\nu\bigr) +O(1)
\end{equation*}
and then $X= \omega \xi_2$ with $\omega =\phi (0,x_2)^{-1}$ satisfies (\ref{2-21}) and therefore the same arguments of the proof of theorem 3.1 \cite{IRO1} remain valid but we consider $h$-pdo with operator-valued symbols (in the axillary space $\bK=L^2(\bR^1)$).
\end{proof}

The following  proposition estimates by $C_0\mu^{-1}\gamma^{-\nu}$ and
$C_0\rho^{-1}$ $x_2$-speed of the propagation from above; the same result will hold for inner zone as well.

\begin{proposition}\label{prop-2-7} (i) In frames of proposition \ref{prop-2-6}(i) let $\psi $ be supported in $\ell$-vicinity of ${\bar y}$ and $\psi_2=\psi_2(x_2)$ satisfy
\begin{equation}
\psi_2=1\qquad {\text in}\quad
\bigl\{ | y_2-{\bar y}_2|\le C_0 (\mu \gamma^\nu)^{-1}T+C_0\ell\bigr\}
\label{2-23}
\end{equation}
with $\ell=C_0\mu^{-1}\gamma ^{1-\nu}$.

Then  for ${\bar T}_0= C_0 h|\log h| \le T\le T_1$
\begin{equation}
|F_{t\to h^{-1}\tau} {\bar\chi}_T(t)
(1-\psi_{2\,x})u\psi_y |\le Ch^s\qquad\forall \tau\in (-\epsilon,\epsilon);
\label{2-24}
\end{equation}

\smallskip
\noindent
(ii) In frames of proposition \ref{prop-2-6}(ii) let $\psi =\psi(x_2)$ be supported in $\ell$-vicinity of ${\bar y}_2$, $\varphi=\varphi (\xi_2)$ in $\epsilon\rho$-vicinity of ${\bar \xi}_2$ and $\psi_2$ satisfy $(\ref{2-22})$. Let 
\begin{equation}
\ell \rho \ge Ch|\log h|.
\label{2-25}
\end{equation}
Then for  ${\bar T}_0 \le T\le T_1$
\begin{equation}
|F_{t\to h^{-1}\tau} {\bar\chi}_T(t)
(1-\psi_{2\,x})u\psi_y \varphi ^t(hD_{y_2})|\le Ch^s\qquad\forall \tau\in (-\epsilon,\epsilon).
\label{2-26}
\end{equation}
\end{proposition}

\begin{proof} Both statements are equivalent under condition (\ref{2-6}). 

\smallskip
\noindent
(i) Assume first that (\ref{2-6}) holds. Then  symbol 
\begin{equation}
Y=x_2-\phi^{-1}\mu ^{-1}x_1^{1-\nu}\xi_1
\label{2-27}
\end{equation}
satisfies 
\begin{equation}
|\bigl\{ a, Y\bigr\}|\le C_0\rho^{-1}
\label{2-28}
\end{equation}
with $\rho=\mu \gamma^\nu$ and symbols
\begin{equation}
\varpi \bigl(C_0\varsigma {\frac t T} \pm {\frac Y {\rho^{-1} T}}\bigr)
\label{2-29}
\end{equation}
($\varsigma = \pm 1$) are quantizable as long as 
\begin{equation}
T\ge T_0 = C\bigl(\gamma h |\log h|\bigr)^{\frac 1 2}
\label{2-30}
\end{equation}
and then standard energy estimates method (see theorem 3.1 \cite{IRO1}) implies (\ref{2-25}) with $\psi$, and $\psi_2$ replaced by $h$-pdo supported and equal 1 respectively in domains
\begin{align*}
&\bigl\{ (y, \eta): |y_1-{\bar y}_1|\le \epsilon \gamma,
|Y(y,\eta_1)-Y(\bar y,\bar\eta_1)|\le \ell\bigr\},\\
&\bigl\{ (x, \xi): |y_1-{\bar y}_1|\le \epsilon \gamma,
|Y(x,\xi_1)-Y(\bar y,\bxi_1)|\le \ell+C\mu^{-1}\gamma^{-\nu} T\bigr\},
\end{align*}
as $T_0\le T \le T_1 = C\mu \gamma^\nu $ and
\begin{equation}
\ell_0=C\bigl(\gamma h |\log h|\bigr)^{\frac 1 2} \rho^{-1} \le \ell \le \epsilon.
\label{2-31}
\end{equation}
Original estimate (\ref{2-24}) follows from this because $|Y(x,\xi_1)-x_2|\le C_0\ell$.

As (\ref{2-24}) is proven for $T=T_0$ it is valid also for 
$T\in [{\bar T}_0,T_0]$.

\smallskip
\noindent
(ii) Assume now that (\ref{2-6}) is violated. Then $x_1$ is localized only in $Ch|\log h|$-vicinity of $0$ rather than in $\gamma$-vicinity of it and $x_1$ is not separated from $0$. However a naive idea that $\mu \phi x_1^\nu/\nu$ should be replaced by $\xi_2$ and thus we should define 
\begin{equation}
Y=x_2-\nu^{-1}x_1\xi_1\xi_2^{-1}
\label{2-32}
\end{equation}
instead of (\ref{2-27}) works; (\ref{2-28}) and (\ref{2-29}) are preserved, also  (\ref{2-30}) becomes
$T_0=Ch|\log h|$ and (\ref{2-31})  becomes (\ref{2-25}).
\end{proof}

The following  proposition estimates by $\epsilon_0\mu^{-1}\gamma^{-\nu}$ (and almost equivalently) by
$\epsilon_0\rho^{-1}$ $x_2$-speed of the propagation from below; this result would fail in the inner zone unless far from periodic trajectories.

\begin{proposition}\label{prop-2-8} (i) In frames of proposition \ref{prop-2-6}(i) let $\psi $ be supported in $\gamma$-vicinity of ${\bar y}$ and $\psi_2=\psi_2(x_2)$ satisfy
\begin{equation}
\psi_2=1\qquad {\text in}\quad
\{ | y_2-{\bar y}_2|\ge
\epsilon_0 (\mu \gamma^\nu)^{-1}T-C_0\ell\}
\label{2-33}
\end{equation}
with $\ell=\mu^{-1}\gamma ^{1-\nu}$.

Then  for $T'_0= C_0 \gamma \le T\le T_1$ inequality $(\ref{2-25})$ holds.

\smallskip
\noindent
(ii) In frames of proposition \ref{prop-2-6}(ii) let $\psi =\psi(x_2)$ be supported in $\ell$-vicinity of ${\bar y}_2$, $\varphi=\varphi (\xi_2)$ in $\epsilon\rho$-vicinity of ${\bar \xi}_2$ and $\psi_2$ satisfy $(\ref{2-32})$. Let $\ell$ satisfy $(\ref{2-26})$.

Then for  ${\bar T}_0 =C_0h|\log h| \le T\le T_1$ inequality $(\ref{2-27})$ holds.
\end{proposition}

\begin{proof} One needs to apply the standard arguments with
\begin{equation}
\varpi \bigl(C_0\varsigma {\frac Y {\rho^{-1} T}}\pm {\frac t T} \bigr)
\label{2-34}
\end{equation}
assuming that inequality (\ref{2-29})  is reversed to
\begin{equation}
|\bigl\{ a, Y\bigr\}|\ge \epsilon_0\rho ^{-1}.
\label{2-35}
\end{equation}
In frames of (\ref{2-6}) one needs to modify $Y$ in the spirit of (\ref{1-28}) which would be equivalent to the following modification in frames of (ii). Namely, note that as $\gamma \ll 1$, $\rho \gg 1$
\begin{equation*}
\{a,Y\}= \xi_2^{-1}\bigl(-\nu^{-1}p_1^2 + p_2^2\bigr)+o(\rho^{-1})=
{\frac 1 2}(1+\nu^{-1})\xi_2^{-1}(p_2^2-p_1^2)+
{\frac 1 2}(1-\nu^{-1})\xi_2^{-1}(p_2^2+p_1^2)
\end{equation*}
and if we redefine $Y$ as
\begin{equation}
Y=x_2-\nu^{-1}x_1\xi_1\xi_2^{-1}+ {\frac 1 2}\nu^{-1}(1+\nu^{-1})x_1\xi_1p_2 \xi_2^{-2}
\label{2-36}
\end{equation}
we arrive to
\begin{equation}
\{a,Y\}=
{\frac 1 2}(1-\nu^{-1})\xi_2^{-1}(p_2^2+p_1^2) +o(\rho^{-1})=
{\frac 1 2}(1-\nu^{-1})\xi_2^{-1} (a+V) +o(\rho^{-1})
\label{2-37}
\end{equation}
which implies (\ref{2-36}) and also symbols (\ref{2-35}) are quantizable.
\end{proof}

From the proof of proposition \ref{prop-2-8} immediately follows that (\ref{2-26}) holds with $\psi(x_2)$, $\psi_2(x_2)$ replaced by $\psi(Y(x,\xi))^\w$, $\psi_2(Y(x,\xi))^\w$ with assumption $\ell=C_0\mu^{-1}\gamma^{1-\nu}$ replaced by (\ref{2-25}) i.e. 
\begin{equation}
\ell = C_0\mu^{-1}\gamma^{-\nu}h|\log h|.
\tag*{$(\ref{2-25})'$}\label{2-25-'}
\end{equation}
Then we immediately get

\begin{corollary}\label{cor-2-9}
(i) In frames of proposition \ref{prop-2-6}(i)
\begin{equation}
|F_{t\to h^{-1}\tau}\chi_T(t)  \Gamma \bigl(u\psi \bigr)|\le Ch^s\qquad \forall \tau\in (-\epsilon,\epsilon);
\label{2-38}
\end{equation}
as ${\bar T}_0=Ch|\log h| \le T\le T_1$ and thus
\begin{equation}
|F_{t\to h^{-1}\tau} \bigl({\bar\chi}_{T_1}(t)-{\bar\chi}_{{\bar T}_0}(t)\bigr)
\Gamma \bigl(u\psi\bigr) |\le Ch^s 
\label{2-39}
\end{equation}
where I remind that ${\bar\chi}$ is supported in $[-1,1]$ and equal $1$ in $[-{\frac 1 2},{\frac 1 2}]$.

\smallskip
\noindent
(ii) In frames of proposition \ref{prop-2-6}(ii) $(\ref{2-38}),(\ref{2-39})$ hold with $\psi (x_1,x_2)$ replaced by $\psi (x_2)\varphi (hD_2)$ with $\varphi$ supported in $\{|\xi_2|\ge C\}$.
\end{corollary}

Assume now that (\ref{2-6}) is fulfilled.  Then  we arrive to
\begin{equation}
|F_{t\to h^{-1}\tau} {\bar\chi}_T(t)\Gamma (u\psi ) |\le 
C\gamma h^{-1} \qquad \forall T\in [T_0,T_1]
\label{2-40}
\end{equation}
as $\psi=\psi'(x_1)\psi'(x_2)$ is an element of $\gamma$-admissible partition with respect to $x_1$; really  this estimate holds with $T=T_0$ due to the non-degenerating  results rescaled (\cite{IRO3}): $x\mapsto x/\gamma$, $h\mapsto h'=h/\gamma$, 
$t\mapsto t/\gamma$ and $\mu \mapsto \mu'=\mu \gamma^\nu$ with $\mu' h'\le 1$.
Then contribution of each element of $\gamma$-admissible (sub)partition with respect to $x$ contributes $C(h/\gamma)^{-1}\times \gamma $ with extra-factor $\gamma$ coming from $dt$ in Fourier transform $F_{t\to h^{-1}\tau}$. After summation with respect to all $\asymp \gamma^{-1}$ such subelements we get (\ref{2-40}),

After (\ref{2-40}) is established we can apply  Tauberian arguments and conclude immediately that the contribution of $\psi$ to the remainder estimate does not exceed
\begin{equation}
CT_1^{-1}\times h^{-1} \gamma= 
C(\mu \gamma^\nu)^{-1}h^{-1}\gamma =C\mu^{-1}\gamma^{1-\nu}
\label{2-41}
\end{equation}
while the principal part is given by formula
\begin{equation}
h^{-1}\int _{-\infty}^0 \Bigl(F_{t\to h^{-1}\tau} {\bar\chi}_T(t)\Gamma \bigl(u\psi\bigr)\Bigr)\,d\tau
\label{2-42}
\end{equation}
with arbitrary $T\in \in [T_0,T_1]$. 

Then  summation of (\ref{2-41}) over  zone $|x_1|\ge \gamma$ results in $C\mu^{-1}h^{-1}\gamma^{1-\nu}$ coinciding with $C{\bar\gamma} h^{-1}$ as $\gamma={\bar\gamma}$ and we arrive to

\begin{proposition}\label{prop-2-10} Under condition $(\ref{2-8})$ contribution of the whole outer zone $\cZ_\out=\{x:|x_1|\ge C_0{\bar\gamma}\}$ to the remainder estimate
\begin{equation}
\cR = |\Gamma (e\psi) - h^{-1}\int _{-\infty}^0 \Bigl( F_{t\to h^{-1}\tau} {\bar\chi}_T(t)\Gamma (u\psi)\Bigr)\,d\tau |
\label{2-43}
\end{equation}
does not exceed $C{\bar \gamma}h^{-1}$ where 
$\psi=1-\phi \bigl(x_1/(C_0{\bar\gamma})\bigr)\psi'(x_1)\psi''(x_2)$ with all functions regular, $\omega$ supported $[-1,1]$ and equal $1$ on $[-{\frac 1 2},{\frac 1 2}]$. 
\end{proposition}

The similar results as condition (\ref{2-8}) violated will be derived later in sections 4, 5. I will replace (\ref{2-42}) by more explicit expression in section 3. Basically it will be the answer prescribed by the non-vanishing magnetic field theory but for $\mu$ large enough some correction will be needed.

\subsection{Inner zone. I}

Inner zone is $\cZ_\inn=\{|\xi_2|\le C_0\}$ or under condition (\ref{2-8}) equivalently (around energy level 0) $\cZ_\inn =\{|x_1|\le C_0\gamma\}$; under condition (\ref{2-18}) $|x_1|\le C_0h|\log h|$ in the microlocal sense in this zone. For a sake of simplicity of notations we use $\gamma$ instead of ${\bar\gamma}$ until the end of the section.

Let us start from the confinement of propagation in the inner zone and a finite speed of propagation; proofs repeating those of the corresponding propositions of the previous subsection are left to the reader. Just to remark that in this zone we study $\xi_2$ and $x_2$ directly without going to $X$ and $Y$.

\begin{proposition}\label{prop-2-11}\footnote{\label{foot-5}Cf proposition \ref{prop-2-6}.} (i)
Let condition  $(\ref{2-8})$ be fulfilled. Let 
${\bar x}\in B(0,{\frac 1 2})\cap  \{x, |x_1|\le  \gamma\}$ and 
$\psi\in C_0^\infty (B({\bar x},\gamma))$ be a rescaling of the standard function.  Let $\psi_1$ be $\gamma$-admissible, supported in 
$\{|x_1|\le 3C_0\gamma\}$ and equal $1$ in 
$\{|x_1|\le 2C_0\gamma\}$. Let $T_1=\epsilon $. Then estimate $(\ref{2-19})$ holds with $T=T_1$.

\smallskip
\noindent
(ii) Without condition $(\ref{2-8})$ let 
$\psi\in C_0^\infty \bigl(({\frac 1 2},{\frac 1 2})\bigr)$, 
$\varphi \in C_0^\infty \bigl((-C_0,C_0)\bigr)$, 
$\varphi_1 \in C_0^\infty \bigl((-3C_0,3C_0)\bigr)$
and equal $1$ in $(-2C_0,2C_0)$. Then estimate $(\ref{2-20})$ holds with ${\bar\xi_2}=1$ and $T=T_1$.
\end{proposition}

\begin{proposition}\label{prop-2-12}\footnote{\label{foot-6}Cf proposition \ref{prop-2-7}.}
(i) In frames of proposition \ref{prop-2-11}(i)  let $\psi_1=\psi(x_2)$,
\begin{equation*}
\psi_1=1\qquad {\text in}\quad
\{  | y_2-{\bar y}_2|\le C_0 T\}.
\end{equation*}
Then for $T_0= C_0\gamma  \le T\le T_1$ inequality 
\begin{equation}
|F_{t\to h^{-1}\tau} {\bar\chi}_T (t)
(1-\psi_{1\,x})uQ_y |\le Ch^s\qquad\forall \tau\in (-\epsilon,\epsilon);
\label{2-44}
\end{equation}
holds with $Q_y=\psi_y$.

\smallskip
\noindent
(ii) Without condition $(\ref{2-8})$ this estimate holds with 
$Q_y=\psi(y_2)\varphi (-hD_{y_2})$ with 
$\varphi \in C_0^\infty \bigl((-C_0,C_0)\bigr)$.
\end{proposition}

Now we are studying the finite speed with respect to $\xi_2$ and the proof repeats those of proposition \ref{prop-2-7}:

\begin{proposition}\label{prop-2-13}
Let $\psi=\psi(x_2)\in C_0^\infty \bigl((-{\frac 1 2},{\frac 1 2})\bigr)$,
$\varphi \in C_0^\infty \bigl(B({\bar \eta}_2,\rho )\bigr)$ and
$\varphi _1\in C_0^\infty \bigl(B({\bar \eta}_2,3\rho )\bigr)$ is equal to $1$ in $B({\bar \eta}_2,3\rho )$ be standard functions rescaled. Let 
\begin{equation}
\rho \ge C_1h|\log h|+C_1\gamma
\label{2-45}
\end{equation}

Then for $T_0= C_2\gamma  \le T\le T_1(\rho)=\epsilon \rho$ inequality 
\begin{equation}
|F_{t\to h^{-1}\tau} {\bar\chi}_T (t)
\bigl(1-\varphi_1(hD_{x_2})\bigr) u\psi_y \varphi(-hD_{y_2})|\le Ch^s\qquad\forall \tau\in (-\epsilon,\epsilon);
\label{2-46}
\end{equation}
holds.
\end{proposition}

We will make more specific statements later. Now we want to prove that as $\xi_2$ is disjoint from $k^*V^{\frac 1 2}$ there is a drift with the velocity $\asymp (\xi_2-k^*V^{\frac 1 2})$ (with the correct sign). Let $\varphi=\varphi(\xi_2)\in C_0^\infty \bigl(B({\bar\eta}_2,\rho)\bigr)$ and
$\psi =\psi (x_2)\in C_0^\infty (B({\bar y}_2,\ell)$, $\ell =\epsilon_1 \rho$
be standard functions rescaled. Assuming that
\begin{equation}
\rho \ge C(h|\log h|)^{\frac 1 2}+ C_0\gamma
\label{2-47}
\end{equation}
we conclude from propositions \ref{prop-2-12}, \ref{prop-2-13} that
\begin{equation*}
|F_{t\to h^{-1}\tau } {\bar\chi}_T(t)\bigl(1- \varphi_1(hD_{x_2})\psi_2(x_2)\bigr)u\psi(y_2) \varphi(-hD_{y_2})|
\le Ch^s\qquad \forall\tau\in (-\epsilon,\epsilon) 
\end{equation*}
as long as $T\le \epsilon \rho$.

\begin{proposition}\label{prop-2-14}  
Let  conditions $(\ref{2-8})$ and $(\ref{1-21})$ be fulfilled. Let
$\psi=\psi(x_2)\in C_0^\infty\bigl(B({\bar y}_2,\ell)\bigr)$ and
$\varphi \in C_0^\infty \bigl(B( {\bar\eta}_2,\rho)\bigr)$ be also standard function rescaled, with $\rho\le C_0$, 
\begin{equation}
{\bar \eta }_2- k^*V(0,{\bar y}_2)^{\frac 1 2}=\pm \rho .
\label{2-48}
\end{equation}
Let $\psi_1= \psi_1(x_2)$ be also a standard function rescaled and equal $1$ in 
$\{|x_2-{\bar y}_2|\ge \epsilon_0\rho T-C_0\ell\}$.

Then estimate $(\ref{2-44})$ holds with $Q_y=\psi (y_2)\varphi (-hD_{y_2})$ as $|\tau |\le \epsilon\rho$,  $\ell\ge\gamma$, 
\begin{align}
&\rho\ell \ge Ch|\log h|,\label{2-49}\\
&T\ge C_0{\frac \ell \rho}\label{2-50}.
\end{align}
\end{proposition}

\begin{proof} With no loss of the generality we can assume that 
$V(0,{\bar y}_2)=1$. Let us consider 
\begin{equation}
Y= x_2- \gamma({\frac {x_1} \gamma}, \xi_1)
\label{2-51}
\end{equation}
where $Z$ is a function introduced in proposition \ref{prop-1-6} and pick up $\gamma=\mu^{-{1/ \nu}}$ exactly. Then
$|\{a, Y\} -\beta (\xi_2-k^*) - \alpha a  |\le \epsilon' \rho$ with a constant $\epsilon'$ which one can make arbitrarily small.

Then we can apply standard arguments with
\begin{equation}
\varpi \Bigl( 
C_0\varsigma {\frac 1 \ell}\bigl(Y(x,\xi)- Y(y,\eta) )\bigr)\pm {\frac t T}\Bigr)
\label{2-52}
\end{equation}
as long as $T\rho \ge C\ell $; symbol (\ref{2-52}) is obviously quantizable.
\end{proof}

It follows from the proof that estimate (\ref{2-44}) holds with 
$Q_y=\psi (Y)^\w \phi(-hD_2)$ and $\psi_1$ replaced by $\psi_1(Y)^\w$ without assumption $\ell \ge \gamma$; this estimate immediately implies

\begin{corollary}\label{cor-2-15}  In frames of proposition \ref{prop-2-14}
estimate $(\ref{2-56})$ (see below) holds with  the standard function $\psi$ and
$|\tau| \le \epsilon \rho$ as long as
\begin{align}
&\rho \ge C_0 (h|\log h|)^{\frac 1 2},\label{2-53}\\
&T_0 \Def {\frac {C_0h|\log h|}{\rho^2}}\le T\le T_1=\epsilon \rho.\label{2-54}
\end{align}
\end{corollary}

Note that $T_0\le \epsilon \gamma$ as long as 
\begin{equation}
\rho \ge {\bar\rho}_1 \Def C \bigl({\frac h \gamma}\bigr)^{\frac 1 2}
\label{2-55}
\end{equation}
and we arrive to

\begin{proposition}\label{prop-2-16} Let conditions $(\ref{2-8}),(\ref{1-21})$ and $(\ref{2-55})$ be fulfilled. 

Then 

\smallskip
\noindent
(i) Estimates 
\begin{equation}
|F_{t\to h^{-1}\tau}\chi_T(t)  \Gamma \bigl(\varphi(hD_2)\psi  u\bigr)|\le Ch^s, 
\label{2-56}
\end{equation}
and
\begin{equation}
|F_{t\to h^{-1}\tau} \bigl({\bar\chi}_{T_1}(t)-{\bar\chi}_{{\bar T}_0}(t)\bigr)
\Gamma \bigl(\varphi(hD_2)\psi  u\bigr) |\le Ch^s 
\label{2-57}
\end{equation}
hold with $|\tau|\le \epsilon\rho$ and ${\bar T}_0\le T\le T_1=\epsilon \rho$.

\smallskip
\noindent
(ii) Furthermore
\begin{equation}
|F_{t\to h^{-1}\tau} {\bar\chi}_{T_1}(t)\Gamma (\varphi(hD_2)\psi  u)|\le Ch^{-1}\gamma \rho.
\label{2-58}
\end{equation}
\end{proposition}

\begin{proof}  We need to cover $T\le \epsilon \gamma$ only.
After rescaling $x \mapsto x/\gamma$, $t\mapsto t/\gamma$, 
$h\mapsto \hbar= h/\gamma$ we find ourselves in frames of the standard propagation and therefore estimate (\ref{2-38}) with $h$ replaced by $\hbar$ would hold as long as $C_0 \hbar |\log \hbar | \le T/\gamma \le \epsilon$.
Then if $\hbar \le h^\delta$ or equivalently $\mu \le h^{\delta -\nu}$ with arbitrarily small exponent $\delta>0$ we arrive to original (\ref{2-38}) with indicated $T$, $\tau$. 

Further, these arguments imply (\ref{2-58}) with no such restriction.

For $h^{\delta-\nu}\le \mu \le \epsilon (h|\log h|)^{-\nu}$ we instead can make a fine $(\epsilon'\gamma, \epsilon')$-subpartition with respect to $(x_1,\xi_1)$ and then depending on subelement apply standard arguments with one of 
\begin{equation}
\varpi \Bigl(C_0\varsigma {\frac 1 T}(x_1-y_1) \pm {\frac t T}\Bigr),\qquad
\varpi \Bigl(C_0\varsigma {\frac 1 T}(\xi_1-\eta_1) \pm {\frac t T}\Bigr)
\label{2-59}
\end{equation}
which are quantizable. Easy details are left to the reader.
\end{proof}

Then applying Tauberian arguments we conclude that contribution of this element to the remainder estimate does not exceed 
$Ch^{-1}\gamma \rho /T_1\asymp Ch^{-1}\gamma$. After summation over partition in $\rho \ge \epsilon_1$ we get remainder estimate $Ch^{-1}\gamma$; we will cover smaller values of $\rho$ in the next subsection. Thus we arrive to

\begin{proposition}\label{prop-2-17} Under conditions $(\ref{2-8}),(\ref{1-21})$ contribution of zone 
\begin{equation*}
\bigl\{|x_1|\le \gamma, 
|\xi_2- k^*V(x_2)^{\frac 1 2}|\ge \epsilon_0\bigr\}
\end{equation*}
to the remainder estimate does not exceed 
$C_\epsilon\mu^{-{ 1 /\nu}}h^{-1}$. 
\end{proposition}

\subsection{Inner zone. II}

Now we need to consider ``near-periodic'' zone
\begin{equation*}
\bigl\{|x_1|\le \gamma, 
|\xi_2- k^*V(x_2)^{\frac 1 2}|\le \epsilon_0\bigr\}.
\end{equation*}

Let us investigate this case under condition (\ref{1-21}), assuming first that 
$\partial_{x_2} V$ is disjoint from 0. We consider a bit more general case: namely $\partial_{x_2} V \asymp \zeta$ with large enough parameter $\zeta$. We consider first propagation assuming that that $|\partial_{x_2} V |\le  \zeta$.

\begin{proposition}\label{prop-2-18} Let conditions $(\ref{2-8}),(\ref{1-21})$   be fulfilled and
\begin{equation}
 |\partial_{x_2} V|\le \zeta
 \qquad {\text as\ } x_1=0
\label{2-60}
\end{equation}
with 
\begin{equation}
\zeta \ge C_0\gamma.
\label{2-61}
\end{equation}
Let $\psi=\psi(x_2)\in C_0^\infty (B(0,{\frac 1 2}))$  be a standard function  and $\varphi  =  \varphi (\xi_2)$ be a standard function rescaled supported in $\{|\xi_2 -k^*V^{\frac 1 2}(0)|\le \rho\}$\,\footnote{\label{foot-7}Then 
under condition (\ref{2-62})\; $|\xi_2 -k^*V^{\frac 1 2}(0,x_2)|\le 2\rho$ for all $x_2$.} with
\begin{equation}
\rho \ge C_0\zeta.
\label{2-62}
\end{equation}
Let $\varphi_1=\varphi_1(\xi_2)$ be also a standard function rescaled, equal $1$ as $|\xi_2 -k^*V^{\frac 1 2}(0,x_2)|\le 3\rho$. Then 

\smallskip
\noindent
(i) Estimate $(\ref{2-46})$ holds as $|\tau |\le \epsilon \zeta$ and 
$Ch\zeta^{-1}|\log h| \le T \le \epsilon$;

\smallskip
\noindent
(ii) In addition, if
\begin{equation}
\rho\ge{\bar\rho}_1\Def (C\gamma^{-1}h|\log h|)^{\frac 1 2}+C\gamma.
\label{2-63}
\end{equation}
estimate $(\ref{2-46})$ holds as $|\tau |\le \epsilon \zeta$ and
\begin{equation}
Ch\zeta^{-1}|\log h| \le T \le 
\epsilon\min\bigl( {\frac 1 \rho},{\frac \rho \zeta}\bigr).
\label{2-64}
\end{equation}
\end{proposition}

\begin{proof} I leave to the reader a standard proof based on the axillary symbol
\begin{equation*}
\varpi\Bigl( C_0\varsigma {\frac t T} \pm {\frac 1 \rho}(\xi_2-\eta_2)\Bigr).
\end{equation*}
(which under our assumptions is quantizable) that the speed of propagation with respect to $\xi_2$ does not exceed $C\zeta$. This implies (i).

To prove (ii) one needs to show in addition that the ``averaged'' speed of propagation with respect to $x_2$ does not exceed $\rho$; under condition (\ref{2-64}) one can easily do it for $t=T^*V(0)^{1/2}$ by means of the standard microlocal analysis. I leave details to the reader. More delicate analysis of subsections 2.6--2.9 will also imply results of this subsection.
\end{proof}

It immediately implies

\begin{corollary}\label{cor-2-19} Under conditions $(\ref{2-18}),(\ref{1-21}),(\ref{2-60})-(\ref{2-63})$
propagation remains confined to
\begin{equation*}
\{|x_2|\le {\frac 3 4}, |\xi_2-k^*V(0)^{\frac 1 2}|\le \rho+\zeta +\zeta T\}
\end{equation*}
as $T\le \epsilon \min\bigl(\rho^{-1}, \zeta^{-{\frac 1 2}}\bigr)$.
\end{corollary}

Now let us estimate from below the propagation speed with respect to $\xi_2$.

\begin{proposition}\label{prop-2-20} Let conditions $(\ref{2-8}),(\ref{1-21}),(\ref{2-63})$   be fulfilled and
\begin{equation}
\epsilon_0\zeta\le  |\partial_{x_2} V|\le \zeta
 \qquad {\text as\ } x_1=0
\label{2-65}
\end{equation}
with $\zeta \ge C_0\gamma$. Let $\varphi =\varphi (\xi_2)$ be a standard function rescaled supported in  $\rho $-vicinity of 
${\bar \eta}_2=k^*V(0)^{\frac 1 2}$.
Then estimate $(\ref{2-46})$ holds as $|\tau |\le \epsilon \zeta$, 
\begin{equation}
T_0\Def C\rho \zeta^{-1} \le T \le T_1\Def 
\epsilon \min\bigl(\rho^{-1}, \zeta^{-{\frac 1 2}}\bigr),
\label{2-66}
\end{equation}
$\varphi_1=1-\varphi$.
\end{proposition}

\begin{proof}  I leave to the reader a standard proof based on axillary symbol
\begin{equation*}
\varpi\Bigl( C_0\varsigma  {\frac 1 \rho}(\xi_2-\eta_2)\pm {\frac t T} \Bigr)
\end{equation*}
which under our assumptions is quantizable.
\end{proof}

Then in frames of this proposition (\ref{2-56}) holds. 
Note that $T_0\le \epsilon \gamma$ as 
\begin{equation}
\zeta \ge C\gamma^{-1}h|\log h|
\label{2-67}
\end{equation}
and due to the same arguments as in proposition \ref{prop-2-16} we arrive to 

\begin{proposition}\label{prop-2-21}
Let conditions of proposition \ref{prop-2-20} and $(\ref{2-67})$ be fulfilled. Then estimates $(\ref{2-56}),(\ref{2-57}),(\ref{2-58})$  hold as 
${\bar T}_0=Ch|\log h|\le T\le T_1$
with $T_1$ defined by $(\ref{2-66})$. 
\end{proposition}

Then due to Tauberian arguments we arrive immediately to

\begin{proposition}\label{prop-2-22} Under conditions $(\ref{2-8}),(\ref{1-21})$ and $(\ref{2-67})$ the contribution of the inner zone to the remainder estimate $\cR$ does not exceed $C\mu^{-{1/ \nu}}h^{-1}$. 
\end{proposition}

In particular $\zeta \asymp 1$ satisfies (\ref{2-67}) under condition (\ref{2-8}).
Case of smaller $\zeta\le Ch\gamma^{-1}|\log h|$ will be considered in subsection 2.7.

\subsection{Inner zone. III}

Let us finish analysis in the ``near-periodic'' but not ``periodic'' zone. To do this I am going to analyze elements of $(x_2,\xi_2)$ partitions on which either $|\xi_2- V^{1/2}k^*|+|\partial _{x_2}V|\ge C{\bar\rho}_1$ branding the rest as ``periodic zone''.

The following results would be the results of the previous subsection rescaled but some conditions are more relaxed because after rescaling 
$x \mapsto  x/\ell$, $\gamma\mapsto\gamma/\ell$. $\zeta\mapsto \zeta\ell$ our new symbols have smaller derivatives (with a factor $\ell$) with respect to $x_1$ than it was assumed there.

\begin{proposition}\label{prop-2-23}\footnote{\label{foot-8}Cf. proposition\ref{prop-2-18} and corollary \ref{cor-2-19}.}
Let conditions $(\ref{2-8}),(\ref{1-21}),(\ref{2-63})$ be fulfilled. 
Consider point  $(\bar y_2,\bar\eta_2)$ and assume that 
\begin{align}
&|{\bar\eta}_2-k^* V^{\frac 1 2}(0,\bar y_2)|\le \rho,\label{2-68}\\
&|\partial_{x_2}^\alpha V|\le C_\alpha\zeta\ell^{1-|\alpha|}
\qquad \text{as\ }x_1=0,\ |x_2-{\bar y}_2|\le \ell\qquad \forall\alpha:|\alpha|\le K, \label{2-69}\\
&\ell \rho \ge Ch|\log h|,\label{2-70}\\
&\ell \zeta \le {\frac 1 2}\rho\label{2-71},\\
&\ell \ge C_0\gamma,\quad \rho \ge C_0\gamma, \quad \zeta\ge C_0\gamma+C_0h^{1-\delta}.  \label{2-72}
\end{align}
Then for $|\tau|\le \rho$, 
\begin{equation}
|t|\le T_1\Def \epsilon \min\bigl({\frac \ell \rho},{\frac \rho \zeta}\bigr)
\label{2-73}
\end{equation}
quantum evolution starting from 
$\{|y_2-{\bar y}_2|\le {\frac 1 3}\ell, 
|\eta_2-{\bar\eta}_2|\le {\frac 1 3}\rho\}$ is confined to
$\{|y_2-{\bar y}_2|\le\ell,|\eta_2-{\bar\eta}_2|\le {\frac 1 2}\rho\}$.
\end{proposition}

\begin{proposition}\label{prop-2-24}\footnote{\label{foot-9}Cf. proposition \ref{prop-2-20}.} Let conditions of proposition \ref{prop-2-23} be fulfilled.

\smallskip
\noindent
(i) If 
\begin{equation}
|{\bar\eta}_2-k^* V^{\frac 1 2}(0,\bar y_2)|\le \rho,\label{2-74}
\end{equation}
then estimate $(\ref{2-56})$ holds for $|\tau|\le \rho$, 
\begin{equation}
T_0\Def C{\frac {h|\log h|}{\rho^2}}\le T\le T_1;
\label{2-75}
\end{equation}

\smallskip
\noindent
(ii) If $(\ref{2-65})$ is fulfilled as $x_1=0$, $|x_2-{\bar y}_2|\le \ell$ then
estimate $(\ref{2-56})$  holds for $|\tau|\le \rho$, 
\begin{equation}
T_0\Def C{\frac {h|\log h|}{\zeta \ell}}\le T\le T_1.
\label{2-76}
\end{equation}
\end{proposition}

Now, given point $(\bar y_2,\bar\eta_2)$ let us define
\begin{equation}
\ell=\rho=\zeta=\varrho(\bar y_2,\bar\eta_2),\qquad \varrho(x_2,\xi_2)\Def \epsilon 
\Bigl(|\xi_2-k^*V(0,x_2)^{\frac 1 2}|+|\partial _{x_2}V(0,x_2)|\Bigr).
\label{2-77}
\end{equation}
Note that (\ref{2-70}), (\ref{2-72}) become equivalent to (\ref{2-63}) and thus   propositions \ref{prop-2-23}, \ref{prop-2-24} hold with $T_1=\epsilon$ under condition (\ref{2-63}). Also note that  $T_0\le \epsilon \gamma$ is equivalent to condition (\ref{2-63}).

Then using the same arguments as above we get  (\ref{2-56}) with 
${\bar T}_0=Ch|\log h|\le T\le T_1$ and also estimate
\begin{equation}
|F_{t\to h^{-1}\tau} {\bar \chi}_T(t) \Gamma \bigl(\varphi(hD_2)\psi u\bigr)|\le Ch^{-1}\rho \gamma \ell\qquad\forall \tau:|\tau|\le \epsilon \rho 
\label{2-78}
\end{equation}
and ${\bar T}_0\le T\le T_1$.

Therefore due to Tauberian arguments contribution of
$(\ell,\rho)$-vicinity of $(\bar y_2,\bar\eta_2)$ to the remainder estimate does not exceed the right-hand expression of (\ref{2-78}) multiplied by $CT_1^{-1}$, i.e. $Ch^{-1}\gamma \rho \ell $; I remind that $\ell=\rho$ here. Then for given $\rho$ the contribution of zone 
$\bigl\{ (x_2,\xi_2), {\frac 1 2}\rho \le \varrho(x_2,\xi_2)\le \rho\bigr\}$ to the remainder estimate does not exceed $Ch^{-1}\rho\gamma$.
Finally, summation over $\varrho$ results in 
$Ch^{-1}\gamma$. This yields statement (i) of the following proposition:

\begin{proposition}\label{prop-2-25} Under conditions $(\ref{2-8}),(\ref{1-21})$ 

\smallskip
\noindent
(i) Contribution of zone $\{(x_2,\xi_2):\varrho (x_2,\xi_2) \ge {\bar\varrho}_1\}$ to the remainder estimate $(\ref{2-43})$ with any $T\in [{\bar T}_0,\epsilon]$ does not exceed
$C\mu^{-{1 /\nu}}h^{-1}$;

\smallskip
\noindent
(ii) The total remainder estimate $(\ref{2-43})$ does not exceed 
\begin{equation}
C\mu^{-{\frac 1 \nu}}h^{-1} + C {\bar\rho}_1 h^{-1};
\label{2-79}
\end{equation}

\smallskip
\noindent
(iii) In particular, as ${\bar\rho}_1 \le {\bar\gamma}$ which is equivalent to 
\begin{equation}
\mu \le C(h|\log h|)^{-{\frac \nu 3}}
\label{2-80}
\end{equation}
the total remainder estimate $(\ref{2-43})$ does not exceed 
$C\mu^{-{1/ \nu}}h^{-1}$.
\end{proposition}

\begin{proof} To prove (ii) one needs to estimate the contribution of \emph{periodic zone\/} $\{\varrho \le C{\bar\rho}_1\}$. One can take 
$T_1= \epsilon \gamma$ there and thus its contribution does not exceed $Ch^{-1}{\bar\rho}_1\gamma \times \gamma^{-1}$.
\end{proof}

\subsection{Periodic orbits. I. Pilot-model} 

Therefore one needs to consider contribution of the {\sl periodic zone\/}
\begin{equation}
\cZ_\per =\bigl\{ (x_2,\xi_2), \varrho(x_2,\xi_2)\le {\bar\varrho}_1= 
C(\gamma^{-1}h|\log h|)^{\frac 1 2}\bigr\}
\label{2-81}
\end{equation}
as (\ref{2-80}) is violated i.e. as
\begin{equation}
C(h|\log h|)^{-{\frac \nu 3}}\le \mu \le \epsilon (h|\log h|)^{-\nu}.
\label{2-82}
\end{equation}

Assume that proposition \ref{prop-2-24} remains valid\footnote{\label{foot-10} Actually it does but with a twist.} as
\begin{equation}
{\bar\rho}_0 =C(h|\log h|)^{\frac 1 2}+C\gamma \le \varrho \le {\bar\rho}_1.
\label{2-83}
\end{equation}
Then for particular partition element inequality
\begin{equation}
|F_{t\to h^{-1}\tau }{\bar\chi}_T(t)\Gamma \bigl(\varphi(hD_2)\psi(x_2)u\bigr)| \le Ch^{-2}\varrho^2\gamma \times  {\frac {h|\log h|}{\varrho^2}}= 
Ch^{-1}\gamma|\log h|
\label{2-84}
\end{equation}
holds as $|\tau|\le \epsilon \rho$, $Ch\rho^{-1}\log h|\le T\le T_1=\epsilon$ where the second factor in the middle expression of (\ref{2-84}) is just $T_0$. Inequality (\ref{2-84}) leads to the estimate of the contribution of the given element to the remainder,  equal to the right-hand expression of (\ref{2-84}) (since $T_1\asymp 1$). Then the total contribution of 
zone $\{(x_2,\xi_2):{\frac 1 2}\rho \le \varrho (x_2,\xi_2)\le \rho\}$ is equal to $Ch^{-1}\gamma\rho^{-1}|\log h|$. This estimate is not only is not good enough for the sharp remainder estimate but is much worse than estimates we had before\footnote {\label{foot-11} One can increase $T_1$ and decrease ${\bar\rho}_0$ in the pilot model cases improving this estimate but (\ref{2-84}) should be improved as well.}. 

To improve the estimates above let us start from the pilot-model operator. However  let us consider first 1-dimensional operator
\begin{equation}
L(k)= {\frac 1 2}\Bigl(\hbar^2D_1^2 + \bigl(k- x_1^\nu/\nu\bigr)^2-1\Bigr)
\label{2-85}
\end{equation}
with $\hbar=h/\gamma$.

\begin{proposition}\label{prop-2-26} As $\hbar \ll 1$ 
\begin{equation}
\Spec \bigl(L(k)\bigr) \cap [-\epsilon,\epsilon]= 
{\check g}\Bigl({\frac {2\pi} {T(k)}}(n+{\frac 1 2})\hbar;k,\hbar\Bigr)+ O(\hbar^2), \qquad n\in \bZ
\label{2-86}
\end{equation}
and consists of simple eigenvalues where ${\check g}(\tau;k,\hbar)$ is an analytic function of all its arguments as 
$|\tau|< \epsilon, \hbar< \epsilon, |k-k^*|<\epsilon$, inverse (with respect to $\tau$) to
\begin{equation}
g(\tau;k,\hbar)= 
T(k)^{-1} \int T\bigl( k(1+2\tau)^{-{\frac 1 2}}\bigr) (1+2\tau)^{\frac {1-\nu}{2\nu}} \,d\tau 
 - {\frac1 2}\hbar + O(\hbar^2)
\label{2-87}
\end{equation}
where $T(k)$ is an elementary period  of section 1 calculated as $\mu=1$, $\tau=0$.

In particular
\begin{equation}
e^{i\hbar^{-1}T(k) g(L;k,\hbar)}{\bar \chi}_\epsilon (L)={\bar \chi}_\epsilon (L) I.
\label{2-88}
\end{equation}
\end{proposition}

\begin{proof} Note that for $L(k)$ Hamiltonian flow on the energy level $\tau$ is periodic with period
$T(\tau,k)=
T\bigl( k(1+2\tau)^{-{\frac 1 2}}\bigr) (1+2\tau)^{\frac {1-\nu}{2\nu}}$.

Proof now follows from Bohr-Sommerfeld formula. To derive asymptotics with $O(\hbar^2)$ error one should note that subprincipal symbol of $L(k)$ is 0 and Maslov' index $\iota_M$ of the trajectory in $(x_1,\xi_1)$ space is $2$.

\begin{figure}[h]
\centering
\subfigure[First pilot-model, $\nu=2$, $k=0.65$]{
\label{fig-6a}
\includegraphics[width=.4\textwidth]{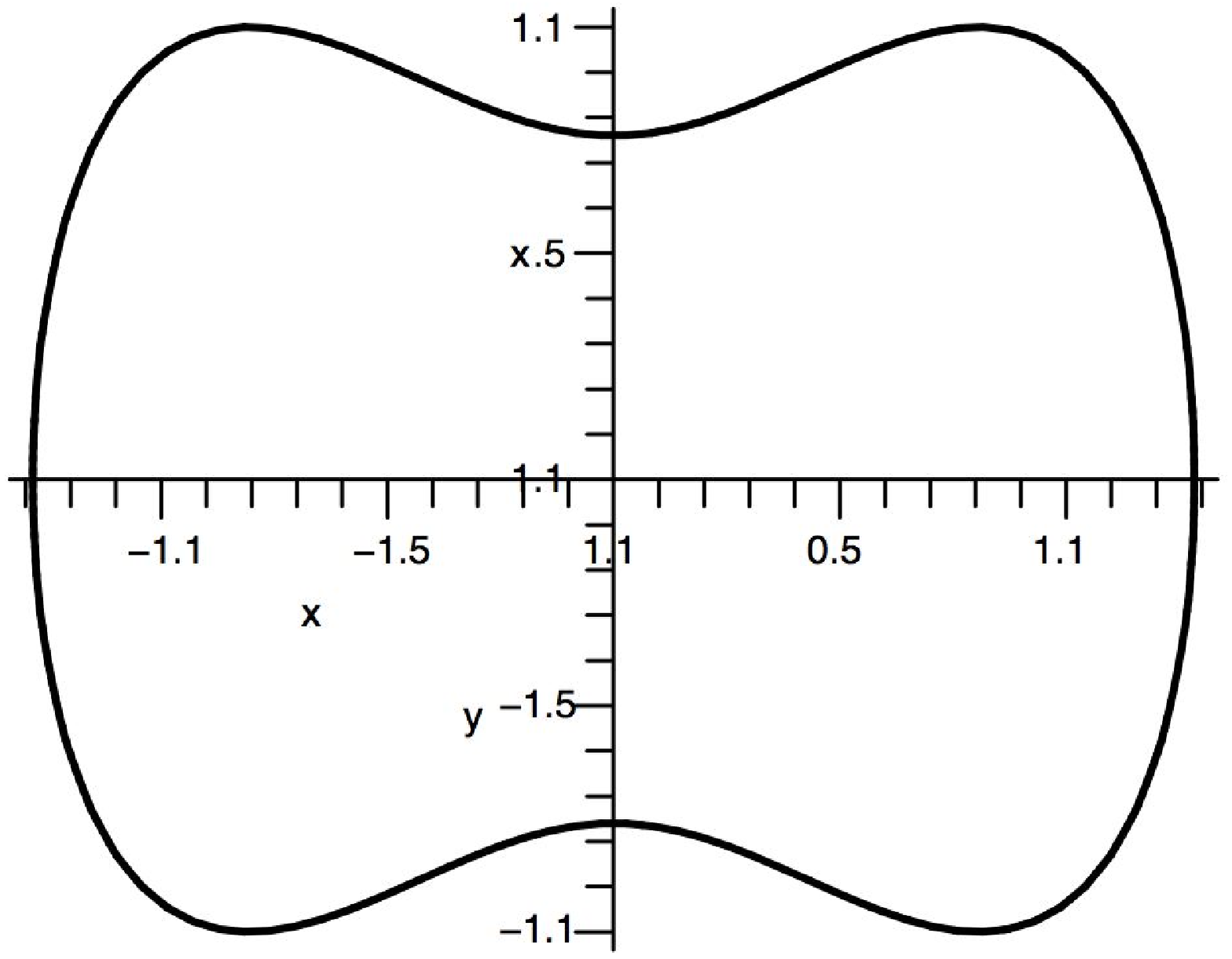}}
\subfigure[Second pilot-model, $\nu=2$, $k=0$]{
\label{fig-6b}
\includegraphics[width=.34\textwidth]{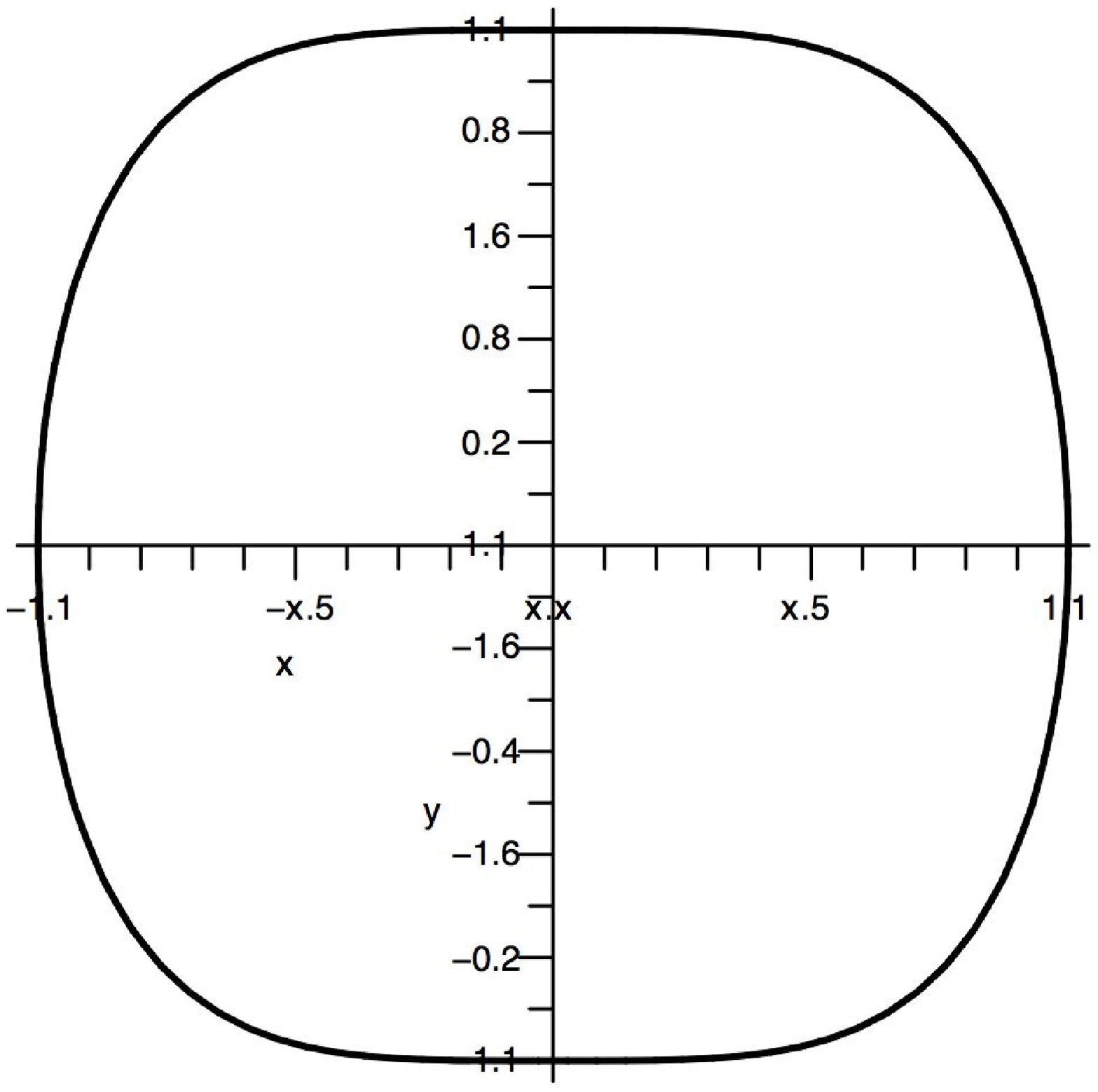}} \\
\label{fig-6}
\end{figure}

\end{proof}

\begin{proposition}\label{prop-2-27} Let condition $(\ref{2-8})$ be fulfilled and let $Q=\varphi(hD_2)$ be a partition element supported in 
$B(k^*,\rho)$. Let us  assume for simplicity that $V(0)=1$. Let $A$ be the pilot-model operator:
\begin{equation}
A= {\bar A} \Def 
{\frac 1 2}\Bigl(\hbar^2D_1^2 + (hD_2- \mu x_1^\nu/\nu)^2-1\Bigr)
\label{2-89}
\end{equation}
Assume that $\gamma \ge Ch|\log h|$. Then modulo negligible operator
\begin{align}
&e^{ih^{-1}\gamma t^* \cA}Q{\bar\chi}_\epsilon (\cA)
\equiv e^{i h^{-1}\gamma B}
Q{\bar\chi}_\epsilon (\cA),\label{2-90}\\
&\cA= g(A,k^*,\hbar), \quad T^*=T(k^*),\label{2-91}
\end{align}
$B=\beta( hD_2,A,\hbar)$ such that 
\begin{equation}
\beta ( \xi_2,\tau,\hbar)= 
\kappa ( \xi_2,\hbar)
\bigl(\xi_2-k_\hbar ^* \bigr)^2 + \kappa_1 (\hbar)+O(|\tau|)
\label{2-92}
\end{equation}
and function $\beta (\xi_2,\tau,0)$ coincides with $b$ defined by $(\ref{1-48})$.
Further, $k_\hbar^*=k^*+O(\hbar)$, $\kappa_1=O(\hbar)$ are analytic functions as $\hbar\ll 1$.
\end{proposition}

\begin{proof} Proof immediately follows from propositions \ref{prop-1-15}, \ref{2-26}. 
\end{proof}

\begin{remark}\label{rem-2-28}
From now in the analysis of periodic zone $k^*=k^*_\hbar$ and respectively $\rho=|\xi_2-k^*_\hbar|$ or $\rho=|\xi_2-k^*_\hbar V^{1/2}|$
\end{remark}

Now (\ref{2-90}) implies that
\begin{multline}
e^{ih^{-1}t\cA}Q{\bar \chi}_{\epsilon\rho} (A) \equiv 
e^{ih^{-1}t''\cA} e^{ih^{-1}t'B}Q{\bar \chi}_{\epsilon\rho} (A),\qquad
t'=\gamma \lfloor {\frac t \gamma }\rfloor, \ t''=t - t'\\
\text{as\ }|t|\le T_1 = \epsilon \rho^{-1}.
\label{2-93}
\end{multline}

Here important thing is that on the partition element in question operator $e^{ih^{-1}t'B}$ is a proper FIO as long as 
\begin{equation}
|t'|\le T'_1\Def {\frac \gamma {\rho^2}}
\label{2-94}
\end{equation}
which is greater than $T'_0=C\rho^{-2}h|\log h|$. Then repeating corresponding analysis of subsections 2.4, 2.5 one can prove easily

\begin{proposition}\label{prop-2-29} Let $\rho$ be defined as in remark \ref{rem-2-28}. Then
proposition \ref{prop-2-18}(ii), corollary \ref{cor-2-19} and proposition \ref{prop-2-24} remains valid with $(\ref{2-63})$ replaced by $(\ref{2-83})$.
\end{proposition}

These statements extended allow to prove

\begin{proposition}\label{prop-2-30} In frames of proposition \ref{prop-2-27}
assume that $\supp (\varphi)\cap B(k^*_\hbar,\rho/2)=\emptyset$ and that $\psi$ is a fixed admissible function. Then
\begin{equation}
|F_{t\to h^{-1}\tau }{\bar\chi}_T(t)\Gamma \bigl(\varphi(hD_2)\psi(x_2)u\bigr)|
\le C\rho^{-1}|\log h|^2.
\label{2-95}
\end{equation}
\end{proposition}

\begin{proof} Due to proposition \ref{prop-2-24} extended by proposition \ref{prop-2-29} it is sufficient to consider $T=T_0$ only. Consider left-hand expression of (\ref{2-95}) with ${\bar\chi}_T(t)$ replaced by $\chi_n(t)={\bar\chi}_\gamma (t-n T^*\gamma)$, 
$|n|\le T_0/T^*\gamma$. Note that then $e^{ih^{-1}\gamma n B}$ is $h$-pdo with symbol which will be regular after rescaling $x_1\mapsto x_1/(\gamma |\log h|)$,
$\xi_1\mapsto \xi_1/|\log h|$, $\xi_2\mapsto \xi_2/|\log h|$. Therefore
\begin{equation}
|F_{t\to h^{-1}\tau }\chi_n(t)\Gamma \bigl(\varphi (hD_2)\psi (x_2)u\bigr)|=
|F_{t\to h^{-1}\tau} \chi_0(t) \Gamma \bigl(Q_nu\bigr)|
\label{2-96}
\end{equation}
with $Q_n = e^{ih^{-1}\gamma n B}\varphi (hD_2)\psi (x_2)$ and therefore expression (\ref{2-96}) does not exceed 
$C\rho \gamma h^{-2}\times h|\log h|= C\rho\gamma h^{-1}|\log h|$ due to the standard theory. 

After summation by $n:|n|\le CT_0/\gamma$ we get 
$C\rho \gamma h^{-1}|\log h| \times h\rho^{-2}\gamma^{-1}|\log h|$ which is exactly the right-hand expression of (\ref{2-95}). 
\end{proof}

I believe that factor $|\log h|^2$ is superficial and one can get rid of it: one  factor $|\log h|$ is due to an inaccuracy in the estimate of (\ref{2-95}) modified and another due to the lost factor 
$\bigl(1+n\gamma/{\bar T}'_0\bigr)^{-s}$ in the same estimate (with 
${\bar T}'_0=C\rho^{-2}h$). However even such weakened estimate is sufficient for our needs unless $\mu$ is really close to $h^{-\nu}$ and the latter special case will be considered separately.

Really, after (\ref{2-95}) is proven,  Tauberian theorem implies that the contribution of this partition element to the remainder estimate does not exceed 
$C\rho^{-1}|\log h|^2 \times T_1^{-1}= C|\log h|^2$ and thus the total contribution of zone $\{|\xi_2-k^*-\hbar|\ge C\gamma\}$ does not exceed this expression integrated over $\rho^{-1}\,d\rho$ resulting in $C|\log h|^3$ while contribution of zone 
$\{|\xi_2-k^*-\hbar|\le C\gamma\}$ does not exceed $C\gamma h^{-1}$ as we already know. Thus we arrive to 

\begin{proposition}\label{prop-2-31} For operator coinciding with the pilot-model $(\ref{2-89})$ in $B(0,1)$ remainder estimate\footnote{\label{foot-12}Which still means $\cR$  given by $(\ref{2-43})$.} is $C\gamma h^{-1}$ as 
$\mu \le Ch^{-\nu}|\log h|^{-2\nu}$.
\end{proposition}

\subsection{Periodic orbits. II. Another example}

Consider now the pilot-model (\ref{2-89}) perturbed by potential $V(x_2)$ satisfying
\begin{equation}
\epsilon_0\zeta \le |\partial_{x_2}V|\le \zeta \qquad\forall x.
\label{2-97}
\end{equation}
Then we can apply the same arguments as above for partition elements with $\rho^2 \ge \zeta$ thus covering completely case $\zeta\le C\gamma^2$. On the other hand, as $\zeta\ge \gamma^2$ we replace $\rho$ by
$\varrho = \max (\zeta^{\frac 1 2}, \rho)$ and on partition elements with $\varrho= \zeta^{\frac 1 2}$ we again can apply the same arguments with $T_1=\epsilon \zeta^{-{\frac 1 2}}$, $T'_0=\epsilon h\zeta^{-1} |\log h|$.

So, we arrive to

\begin{proposition}\label{prop-2-32} For operator coinciding in $B(0,1)$ with pilot-model perturbed by potential $V(x_2)$ satisfying $(\ref{2-90})$ remainder estimate does not exceed $O(\gamma h^{-1})$ as long as 
$\zeta^{\frac 1 2}+\gamma\ge Ch|\log h|^3$ i.e. either 
$\mu \le Ch^{-\nu}|\log h|^{-3\nu}$ or 
$\zeta \ge Ch^2|\log h|^6$.
\end{proposition}

\subsection{Periodic orbits. III. General case}

Consider now general case. Let $\varphi (\xi_2)$ and $\psi(x_2)$ are two $\rho$- and $\ell$-admissible functions supported in $\rho/3$- and $\ell/3$-vicinities of ${\bar \xi}_2$ and ${\bar x}_2$ respectively. Let us assume that conditions (\ref{2-68})--(\ref{2-72}) are fulfilled.

Then due to proposition \ref{prop-2-24} extended by proposition \ref{prop-2-29} as 
\begin{equation}
|t|\le T_1 \Def \epsilon \min \bigl({\frac \ell \rho}, {\frac \rho \zeta}\bigr)\asymp \epsilon \rho \ell \bigl(\rho^2 +\ell\zeta\bigr)^{-1}
\label{2-98}
\end{equation}
propagation started from $\supp \varphi\times \supp \psi$ is confined to $(\rho/3, \ell/3)$-vicinity of it. 

Now we can easily generalize proposition \ref{prop-2-27}:

\begin{proposition}\label{prop-2-33} Let conditions $(\ref{2-8}),(\ref{1-21}),(\ref{2-68})-(\ref{2-72})$ be fulfilled. Assume for simplicity that $V(0)=1$

\smallskip
\noindent
(i) Identity $(\ref{2-90})$ holds as $\cA$ is defined by $(\ref{2-91})$;

\smallskip
\noindent
(ii) Therefore $(\ref{2-93})$ holds.
\end{proposition}

\begin{proof} Proof  repeats one of proposition \ref{prop-2-27}. I leave details to the reader.\end{proof}

\begin{proposition}\label{prop-2-34} Let conditions $(\ref{2-8}),(\ref{1-21}),(\ref{2-68})-(\ref{2-72})$ be fulfilled. Assume for simplicity that $V(0)=1$. Then estimate
\begin{equation}
|F_{t\to h^{-1}\tau} {\bar\chi}_T(t)\Gamma \bigl(\varphi (hD_2)\psi (x_2)u\bigr)|\le CT\rho  \ell h^{-1} |\log h|
\label{2-99}
\end{equation}
holds as 
\begin{equation}
\gamma \le T\le T'_0 \Def C_1 h|\log h| (\rho^2+\zeta \ell)^{-1}.
\label{2-100}
\end{equation}
\end{proposition}

\begin{proof} Let us  note that for $|t'|\le T'_0$ operator $e^{ih^{-1}t'B}$ remains legitimate $h$-PDO and therefore contribution of the time interval 
$[t_1,t_2]$ with $|t_j|\le T'_0$, $|t_2-t_1|=\epsilon\gamma$ does not exceed $C\rho\ell \gamma h^{-1}|\log h|$ and therefore the contribution of 
$\asymp T/\gamma$ of such intervals does not exceed the right-hand expression of (\ref{2-99}).
\end{proof}

It immediately implies

\begin{corollary}\label{cor-2-35} Let in frames of the previous propositions 
either 
\begin{align}
&|\xi_2-k^*_\hbar V(0,{\bar x}_2)|\ge \rho,\label{2-101}\\
\intertext{or}
&|\partial_{x_2}V|\ge \epsilon_0\zeta\label{2-102}
\end{align}
in $\ell$-vicinity of ${\bar x}_2$. Then 

\smallskip
\noindent
(i) Left-hand expression $(\ref{2-99})$ with $T\le T_1$ does not exceed $CT'_0\rho \ell h^{-1}|\log h|^2$;

\smallskip
\noindent
(ii) Contribution of the partition element $\varphi(hD_2)\psi(x_2)$ to the remainder estimate does not exceed $C|\log h|^2$.
\end{corollary}

\begin{proof}
Statement (i) immediately follows from propositions \ref{prop-2-24} extended and \ref{prop-2-34}; statement (ii) follows from the Tauberian arguments: contribution of this element to the remainder estimate does not exceed 
$Ch^{-1}\rho\ell |\log h|\times T'_0/T_1= C|\log h|^2$
since $T'_0/T_1= Ch|\log h|/(\rho\ell)$.
\end{proof}

Assume now that 
\begin{equation}
|\partial _{x_2}V|+|\partial^2_{x_2}V|\ge \epsilon_0.
\label{2-103}
\end{equation}
Let us introduce functions $\varrho$ by (\ref{2-76}) and let us consider $\varrho$-admissible partition in $(x_2,\xi_2)$ and apply the same arguments with $\rho=\zeta=\ell=\varrho$; then on each partition element either
$|\xi_2M- k^*_\hbar V(0,x_2)|\ge \epsilon_0 \varrho$ or 
$|\partial_{x_2}V(0,x_2)|\ge \epsilon \varrho$ and everything works as long as 
\begin{equation}
\varrho \ge C\gamma, \qquad \varrho \ge C(h|\log h|)^{\frac 1 2}
\label{2-104}
\end{equation}
where the second inequality is equivalent to $C\ell\rho\ge Ch|\log h|$.

We know that the contribution of each such  element to the remainder estimate does not exceed $C|\log h|^2$ and therefore the total contribution of all such elements does not exceed $C\log h|^2 I$ with 
$I=\int \varrho^{-2}\,dx_2 d\xi_2\asymp C|\log h|$ due to condition (\ref{2-103}).

On the other hand, as $\varrho \le {\bar\varrho}=
C\gamma + C(h|\log h|)^{\frac 1 2}$ let us redefine $\varrho$ as ${\bar\varrho}$. Then contribution of each such partition element to the remainder estimate does not exceed 
$C{\bar\varrho}^2h^{-1}= C\gamma^2h^{-1}+C|\log h|^2$
and  again due to condition (\ref{2-103}) there is no more than $C$ of such elements. Thus we arrive to

\begin{proposition}\label{prop-2-36} For operator satisfying in $B(0,1)$ conditions  $(\ref{1-21})$ and $(\ref{2-103})$ the remainder estimate is $O(\gamma h^{-1})$ as long as  $\mu \le Ch^{-\nu}|\log h|^{-3\nu}$.
\end{proposition}

On the other hand, exactly the same approach but without condition $(\ref{2-103})$ results in the remainder estimate
\begin{equation*}
C{\bar\varrho}^{-1}|\log h|^2 + C{\bar\varrho}h^{-1}
\end{equation*}
where the first term is an upper estimate  of $C|\log h|^2 I$ with integral over $\{\varrho \ge {\bar\varrho}\}$ and the second term is the contribution of elements with $\varrho\asymp {\bar\varrho}$. Picking up 
${\bar\varrho}=C\max(\gamma, h^{\frac 1 2}|\log h|)$ we arrive to
estimate 
\begin{equation*}
 C\gamma h^{-1}+ Ch^{-{\frac 1 2}}|\log h|
\end{equation*}
which is our target $C\gamma h^{-1}$ as long as 
$\gamma \ge h^{\frac 1 2}|\log h|$ i.e.
$\mu \le Ch^{-{\frac \nu 2}}|\log h|^{-\nu}$. Thus we arrive to

\begin{proposition}\label{prop-2-37} For operator satisfying in $B(0,1)$ condition  $(\ref{1-21})$  

\smallskip
\noindent
(i) The remainder estimate is $O(\gamma h^{-1})$ as long as  
$\mu \le h^{-{\frac \nu 2}}|\log h|^{-\nu}$; 

\smallskip
\noindent
(ii) The remainder estimate is $O(h^{-{\frac 1 2}}|\log h|)$ as long as $h^{-{\frac \nu 2}}|\log h|^{-\nu}\mu \le Ch^{-\nu}|\log h|^{-\nu}$.

\subsection{Periodic orbits. IV. General case (continuation)}

Finally  I am going to prove the most general

\begin{proposition}\label{prop-2-38} Under conditions $(\ref{2-8}),(\ref{1-21})$ and
\addtocounter{equation}{1}
\begin{equation}
\sum_{1\le k\le m}|\partial_{x_2}^kW|\ge \epsilon_0.
\label{2-105}
\tag*{$(\theequation)_m$}
\end{equation}
with $W(x_2)=V(0,x_2)$ the remainder does not exceed 
$C\bigl(|\log h|^{m+1}+\gamma h^{-1}\bigr)$ and thus the remainder estimate $O(\bigl(\gamma h^{-1}\bigr)$ as long as $\mu \le Ch^{-\nu}|\log h|^{-(m+1)\nu}$;

\medskip
\noindent 
(ii) Without condition \ref{2-105} the remainder  does not exceed 
$C\bigl(\gamma ^{-\delta_1}+\gamma h^{-1}\bigr)$ and thus the remainder estimate is $O(\bigl(\gamma h^{-1}\bigr)$ as long as $\mu \le Ch^{-\nu +\delta_1}$ with arbitrarily small $\delta_1>0$.
\end{proposition}

\begin{proof} Let us consider the scaling function
\addtocounter{equation}{1}
\begin{equation}
\ell (x_2) = \epsilon_1 \Bigl(\sum_{1\le k\le m-1}
|\partial_{x_2}^k W(x_2)|^{\frac m {m-k}}\Bigr)^{\frac 1 m}+{\bar\ell}
\tag*{$(\theequation)_m$}
\label{2-106}
\end{equation}
with small enough constant $\epsilon_1>0$  such that 
$|\partial_{x_2} \ell|\le {\frac 1 2}$. Let us consider $\ell$-admissible covering of $(-1,1)$.

Then 
\begin{equation}
|\partial_{x_2}^kW|\le C_0\ell(x_2) ^{m-k}\qquad \forall k:k=1,\dots,m-1;
\label{2-107}
\end{equation} 
further, on each element with $\ell \ge 2{\bar\ell}$  one of these inequalities could be reversed with $C_0$ replaced by $1$:
\begin{equation}
\addtocounter{equation}{1}
|\partial_{x_2}^kW|\le C_0\ell(x_2) ^{m-k}\qquad k=1,\dots,m-1.
\tag*{$(\theequation)_k$}
\label{2-108}
\end{equation}

Furthermore, under condition \ref{2-105} there is no more than $C_0$ of $\ell$-elements with $\ell \in (l,2l)$ for each $l$ and thus no more than 
$C|\log {\bar\ell}|$ elements in total while without   condition \ref{2-105} the corresponding estimates are $C l^{-1}$ and $C{\bar\ell}^{-1}$ respectively.

On each $\ell$-element oscillation of $W(x_2)$ does not exceed $C\ell^m$ and thus we can ponder condition 
\begin{equation*}
|\xi_2-k^*_\hbar W({\bar y})^{\frac 1 2}|\asymp \rho
\end{equation*}
with  $\rho \ge C\ell^m$ and ${\bar y}$ the center of this element. As this condition is fulfilled and also $\rho\ell\ge Ch|\log h|$, $\rho \ge h^{1-\delta}$ we can apply corollary \ref{cor-2-35}(ii) and conclude that the contribution of this element to the remainder estimate does not exceed $C|\log h|^2$. 

Now we need to sum contributions of all such elements. To do it in the most efficient way we redefine 
\addtocounter{equation}{1}
\begin{equation}
 \ell(x_2,\xi_2)= \epsilon_1 \Bigl(|\xi_2-k^*_\hbar W^{\frac 1 2}(x_2)|+
\sum_{1\le k\le m-1}
|\partial_{x_2}^k W(x_2)|^{\frac m {m-k}}\Bigr)^{\frac 1 m}+{\bar\ell}.
\tag*{$(\theequation)_m$}
\label{2-109}
\end{equation}
Then for each given $\rho$ the total contribution of all elements of this type would not exceed $C|\log h|^2$ and $C\rho^{-{\frac 1 m}}|\log h|^2$ under condition  \ref{2-105} and without it respectively; after integration with respect to $\rho^{-1}d\rho$ we get $C|\log h|^3$ and 
$C{\bar\ell}^{-1}|\log h|^2$ respectively and both these expressions do not exceed stated in (i),(ii) remainder estimates\footnote{\label{foot-13}Provided $m$ is large enough in (ii) which we assume.} as we pick up
\begin{equation}
{\bar \ell}= C_1\gamma^{1/(m-1)} +C_1(h|\log h|)^{1/(m+1)}.
\label{2-110}
\end{equation}

\smallskip
On the other hand, the contribution of each  element with 
$\ell\asymp {\bar\ell}$ does not exceed $C{\bar\ell}^{m+1}h^{-1}|\log h|^2$ and the total contribution  of such elements does not exceed 
$C{\bar\ell}^{m+1}|\log h|^2h^{-1}$ and 
$C{\bar\ell}^m|\log h|^2h^{-1}$ under condition \ref{2-105} and without it respectively and these expressions do not exceed remainder estimates announced  in (i),(ii) respectively.

\smallskip
Thus we need to consider only elements with $\rho \le C\ell^m$ and 
$\ell\ge 2{\bar\ell}$ with $\ell$ given by \ref{2-106} rather than \ref{2-109}. Let us consider one such element 
$({\bar y}-{\frac 1 2}\ell,{\bar y}+{\frac 1 2}\ell)$; after rescaling 
\begin{equation*}
x\mapsto (x-{\bar y})/\ell,\; 
W(x_2)-W({\bar y})\mapsto \bigl(W(x_2)-W({\bar y})\bigr)/\ell^m,\;
\rho \mapsto \rho \ell^{-m}
\end{equation*}
we find ourselves in frames of $(2.105)_{m-1}$.

So in the rescaled coordinates let us define scaling function $\ell'_1$ by $(2.109)_{m-1}$ and redefine $W$ as $\ell^{-m}W$. However, as we return to the original scale this new rescaling function $\ell_1 =\ell'_1\ell$ will be defined by $(2.109)_{m-1}$ with original $W$. Now we can repeat the same arguments as before with ${\bar\ell}_1$ defined by $(\ref{2-110})$ with $m$ replaced by $(m-1)$. However now in the total partition there is no more than $C_1$ of $\ell_1$ elements with $\ell_1\in (l,2l)$ if original condition \ref{2-105} was fulfilled and no more than $C_1l^{-1}$ of them otherwise because as I mentioned after rescaling of the elements in question $(2.105)_{m-1}$ always is fulfilled. Then we eliminate elements of $\ell_1$ partition with $\ell_1\le 2{\bar\ell}_1$ and with $\rho \ge C_1\ell_1^{m-1}$ estimating their total contribution by 
$C|\log h|^{m+1}+C\gamma h^{-1}$ and 
$C|\log h|^{m+1} + C {\bar\ell}_1^{m+1}{\bar\ell}^{-1}h^{-1}$ in frames of (i), (ii) respectively which do not exceed  remainder estimates announced there.

We will continue this process as long as our scaling function $\ell_k$ is defined just by $(2.109)_{m-k}$ i.e. until $k=m-2$. Then all elements remaining are eliminated by corollary \ref{cor-2-35}(ii) under condition 
$|\partial_{x_2} W|\asymp \zeta$.
\end{proof}

Thus I have covered the case 
\begin{equation}
\mu \le h^{-\nu}|\log h|^{-\nu N}
\label{2-111}
\end{equation}
leaving the case of superstrong magnetic field
$h^{-\nu}|\log h|^{-\nu N}\le \mu \le Ch^{-\nu}$ for section 4.
\end{proposition}

\sect{Calculations}

So far we recovered good remainder estimates but the asymptotical expression (\ref{2-42}) I derived was implicit. In this section I am going to rewrite it in more explicit form as either 
\begin{equation}
\int \cE^\MW (x,0)\psi (x)\,dx   
\label{3-1}
\end{equation}
which is of magnitude $h^{-2}\bigl(1+\mu h\bigr)^{-1/(\nu -1)}$ or the same expression but with a smaller correction term. Again I consider only the case (\ref{2-111}) leaving the case of superstrong magnetic field for section 4.

\begin{remark}\label{rem-3-1} With no loss of the generality one can assume that $\psi=\psi_2(x_2)\psi_1(x_1)$ with fixed admissible functions $\psi_j$; moreover, one can take $\psi_=1$ as $\mu\ge Ch^{-1}$. Really, it follows from the analysis of subsection 0.1 that for fixed $\psi $ vanishing as $x_1=0$ the remainder estimate is $O({\bar\gamma}h^{-1})$ while the principal part is given by (\ref{3-1}).
\end{remark}

\subsection{Outer Zone}

I remind that the outer zone 
$\cZ_\out = \bigl\{C{\bar\gamma}\le |x_1|\le {\bar\gamma}_1\bigr\}$ and in this subsection $\gamma$ ranges from ${\bar\gamma}$ to ${\bar\gamma}_1$ as well. I also I remind that according to subsection 2.2 we can take here $T=Ch|\log h|$ in (\ref{3-1}).

\begin{proposition}\label{prop-3-2} Let  condition $(\ref{2-8})$ be fulfilled. Then
\begin{multline}
|h^{-1}\int_{-\infty}^0 
\Bigl(F_{t\to h^{-1}\tau}{\bar\chi}_T(t)\Gamma 
\bigl(\psi (1-\psi_1)u\bigr)\Bigr)\,d\tau -\\
\int  \cE^\MW (x,0)\bigl(\psi (1-\psi_1)\bigr)\,dx|\le
C(\mu h+1)^{1/(\nu-1)}
\label{3-2}
\end{multline}
for $T=Ch|\log h|$, fixed function $\psi$ and $\psi_1(x_1)={\bar\psi}_1(x_1/{\bar\gamma})$ with a fixed admissible function ${\bar\psi}$, supported in $(-2,2)$ and equal $1$ in $(-1,1)$.

In particular we get $C{\bar\gamma}h^{-1}$-estimate as $\mu \le Ch^{-\nu^2/(2\nu-1)}$.
\end{proposition}

\begin{proof} It follows from analysis of \cite{IRO3} that for 2-dimensional Schr\"odinger operator with non-degenerate magnetic field, binding parameter $\mu$ and semiclassical parameter $h$ under non-degeneracy condition $|\nabla V/F|\ge \epsilon_0$ the following asymptotic expansion holds
\begin{equation}
h^{-1}\int_{-\infty}^0 
\Bigl(F_{t\to h^{-1}\tau}{\bar\chi}_T(t)\Gamma 
\bigl(\psi u\bigr)\Bigr)\,d\tau \sim \sum_{m,n\ge 0} 
\kappa_{mn} h^{-2+2n+2m } \mu^{2m}
\label{3-3}
\end{equation}
as $\mu \le h^{\delta-1}$ where 
\begin{equation}
\sum_{m\ge 0} \kappa_{m0} h^{-2}(\mu h)^{2m}\sim 
\int  \cE^\MW (x,0)\psi (x)\,dx
\label{3-4}
\end{equation}
and therefore one can replace (\ref{3-3}) by (\ref{3-4}) with an error $\kappa_{02}  +O((\mu h)^2)$.

It also follows that in the same settings but with $h^{\delta-1}\le \mu\le Ch^{-1}$ another asymptotic expansion holds  
\begin{equation*}
h^{-1}\int_{-\infty}^0 
\Bigl(F_{t\to h^{-1}\tau}{\bar\chi}_T(t)\Gamma 
\bigl(\psi u\bigr)\Bigr)\,d\tau \sim
\sum_{m,n\ge 0} \kappa'_{m,n} h^{-2+2n+2m } +
\sum_{m,n\ge 0} \kappa''_{m,n} \mu^{-2-2m} h^{-2+2n}
\end{equation*}
and comparison of this two expression with (\ref{3-3}) where they overlap implies that 
$\kappa''_{**}=0$ and an error in question is given again by $\kappa_{02}  +O((\mu h)^2)$.

Therefore in our setting  the contribution (to this error) of any element of $\gamma$-partition in the outer zone will be of magnitude 1. To calculate the total contribution of the outer zone one must integrate the latter expression  over $\gamma^{-2}d\gamma$  resulting in $O({\bar\gamma}^{-1})$ which is not as good as 
$O({\bar\gamma}_1^{-1})$ with ${\bar\gamma}_1=\min (1, (\mu h)^{-1/(\nu-1)}$ 
which I claimed. However plugging $\mu\mapsto \mu\gamma^\nu$, 
$h\mapsto h/\gamma$ into $O((\mu h)^2)$ and integrating produces $O({\bar\gamma}_1^{-1})$ which is exactly what I claimed and one needs just to calculate explicitly the correction appearing from $\kappa_{02}$-term.

In the non-degenerate settings this latter correction does not depend on $\mu$   and therefore one can calculate it as $\mu=0$. But  this is non-magnetic case and then due to the standard theory and condition $|V|\ge \epsilon_0$
\begin{equation*}
h^{-1}\int_{-\infty}^0 \Bigl(F_{t\to h^{-1}\tau} \bigl({\bar\chi}_T(t)u\bigr)\Bigr)\,d\tau \sim
\sum_{n\ge 0} \varkappa _n(x) h^{-2+2n}
\end{equation*}
(without integration with respect to $dx$) with 
$\varkappa_0 (x) h^{-2}= \cE^\W (x,0)$ and therefore the total contribution of this correction term is $O(1)$.
\end{proof}

Now I am going to derive a bit sharper remainder estimate; to do this I must appeal to the spectral projector of the corresponding one-dimensional Schr\"odinger operator:

\begin{proposition}\label{prop-3-3} Let condition $(\ref{1-21})$ be fulfilled and let $\psi_1$ be an even function, $\psi=\psi(x_2)$. Then as 
$h^{-\nu^2/(2\nu-1)}\le \mu \le \epsilon (h|\log h|)^{-\nu}$, $T=\epsilon {\bar\gamma}$  the following estimate holds:
\begin{multline}
|h^{-1}\int_{-\infty}^0 
\Bigl(F_{t\to h^{-1}\tau}{\bar\chi}_T(t)\Gamma 
\bigl(\psi (1-\psi_1)u\bigr)\Bigr)\,d\tau -\\
\int \Bigl( \cE^\MW(x,0)-\cE_0^\MW(x,0)\Bigr) \psi (x_2)(1-\psi_1(x_1))\,dx -\\
(2\pi h)^{-1}\int {\rm e}_0(x_1, x_1;x_2, \xi_2, 0, \hbar) \psi (x_2) \bigl(1-{\bar\psi}_1(x_1)\bigr)\,dx d\xi_2| \le \\
Ch^{-1}{\bar\gamma}+C|\log h|^K
\label{3-5}
\end{multline}
where ${\rm e}_0(x_1,y_1;x_2,\xi_2,\tau,\hbar)$ is the Schwartz kernel of the spectral projector ${\bf e}_0(x_2,\xi_2,\tau,\hbar)$ for 1-dimensional Schr\"odinger operator
\begin{equation}
{\bf a}_0(x_2,\xi_2,\hbar)= {\frac 1 2}\Bigl(\hbar^2 D_1^2 + (\xi_2- x_1^\nu/\nu )^2 -W(x_2)\Bigr)
\label{3-6}
\end{equation}
with $\hbar=h/{\bar\gamma}$, while 
\begin{equation}
\cE_0^\MW (x,0)= {\frac 1 {2\pi}} \sum_{n\ge 0}
\theta \Bigl(2\tau +V(x_2)-(2n+1)\mu h x_1^{\nu-1}\Bigr) \mu h ^{-1}x_1^{\nu-1}
\label{3-7}
\end{equation}
is   Magnetic Weyl approximation for corresponding pilot-model operator.
\end{proposition}

\begin{proof} (i) To calculate $u(x,y,t)$ let us apply first 
the method of  successive approximations with unperturbed operator ${\bar A}$ obtained from $A$ by freezing $x_2=y_2$ in $V$, $\sigma$, $\phi$:
\begin{equation}
{\bar A}=
{\frac 1 2}\Bigl(h^2D_1^2 + {\bar\sigma}(x_1)^2 \bigl(\xi_2-\mu {\bar\phi }(x_1) x_1^\nu )/\nu \bigr)^2 -{\bar V}(x_1)\Bigr),
\label{3-8}
\end{equation}
${\bar\sigma}=\sigma(x_1,y_2)$, ${\bar\phi}=\phi(x_1,y_2)$, 
${\bar V}=V(x_1,y_2)$.

Then the first term of  
\begin{equation}
h^{-1}\int_{-\infty}^0 
\Bigl(F_{t\to h^{-1}\tau}{\bar\chi}_T(t)\Gamma 
\bigl(\psi (1-\psi_1)u\bigr)\Bigr)\,d\tau
\label{3-9}
\end{equation}
is also given by formula (\ref{3-9}) but for operator ${\bar A}$ instead of $A$ and it is $O\bigl(h^{-3}{\bar\gamma}_1T_0\bigr)=O\bigl(h^{-2}{\bar\gamma}_1|\log h|\bigr)$ \footnote{\label{foot-14}Actually factor $|\log h|$ is superficial and the first term is of magnitude $h^{-2}{\bar\gamma}_1$; further, all logarithmic factors below are superficial}.

The estimate for the second term gains the factor $CT_0^2/h\asymp Ch|\log h|^2$ 
provided $\phi =1$ identically or factor 
\begin{equation*}
CT_0^2 \bigl(1+ \mu {\bar\gamma}_1^{\nu+1}h\bigr)/h\asymp C\bigl(h+{\bar\gamma}_1^2\bigr)|\log h|^2
\end{equation*}
in the general case (I remind that then $\phi=1$ as $x_1=0$); therefore the second term does not exceed  
$Ch^{-2}{\bar\gamma}_1\bigl(h+{\bar\gamma}_1^2\bigr)|\log h|^3$.

Finally, in two-term approximation the remainder gains one more such factor and thus does not exceed 
$Ch^{-2}{\bar\gamma}_1(h^2+{\bar\gamma}_1^4)|\log h|^5$.  One can check easily that this is less than $Ch^{-1}{\bar\gamma}$ unless $\nu = 2$ and 
$\mu \le h^{-\nu^2/(2\nu -1)}|\log h|^K$ in which \emph{special case\/}  there will be an extra 
$|\log h|^K $ factor (rather superficial); then we consider three-term approximation and one can see easily that the remainder estimate will be $Ch^{-1}{\bar\gamma}$.

So we are left only with expression (\ref{3-9})  for
\begin{align}
&u\mapsto {\bar u}=-ih\sum_{\varsigma=\pm }
{\bar G}^{\varsigma}\delta(t)\delta(x_2-y_2)\delta (x_1-y_1)\label{3-10}\\
\intertext{where here and below ``$\mapsto$'' means ``replaced by'' and two extra terms given by (\ref{3-9}) again for} 
&u\mapsto {\bar u}_1=-ih\sum_{\varsigma=\pm } 
{\bar G}^{\varsigma}A_1{\bar G}^{\varsigma}
\delta(t)\delta(x_2-y_2)\delta (x_1-y_1),\label{3-11}\\
&u\mapsto {\bar u}_2=-ih\sum_{\varsigma=\pm } 
{\bar G}^{\varsigma}A_1{\bar G}^{\varsigma}
A_1{\bar G}^{\varsigma}\delta(t)\delta(x_2-y_2)\delta (x_1-y_1)\label{3-12}
\end{align}
where ${\bar G}^\pm$ are corresponding parametrices and $A_1=A-{\bar A}$. 

Also, note that replacing in (\ref{3-11}) $A_1$ by 
${\bar A}_1=(x_2-y_2)\cdot B_1$, $B_1=(\partial_{x_2} A)|_{x_2=y_2}$ one gets a proper remainder estimate with the exception of the special case in which in addition to terms due to (\ref{3-10})-(\ref{3-12}) one must consider also an extra term given by (\ref{3-9}) for 
\begin{equation}
u\mapsto {\bar u}'_2=-ih\sum_{\varsigma=\pm } 
{\bar G}^{\varsigma}A_2{\bar G}^{\varsigma}
\delta(t)\delta(x_2-y_2)\delta (x_1-y_1)
\label{3-13}
\end{equation}
with ${\bar A}_2= {\frac 1 2}(x_2-y_2)^2\cdot (\partial^2_{x_2} A)|_{x_2=y_2}$.

\smallskip
\noindent
(ii) I will tackle  terms due to (\ref{3-11})-(\ref{3-13}) later; now let us concentrate on the main term due to (\ref{3-10}). After this term is derived with $T=Ch|\log h|$ one can increase $T$ to $T_1$ where 
$T_1=\epsilon \mu \gamma^\nu $; however arguments of subsection 2.2 applied to ${\bar A}$  (whose coefficients does not depend on $x_2$) imply that now one can take even $T=+\infty$. 

With $T=+\infty$  one can rewrite this term (before integration by $\psi$ and integration over $x_2$) as
\begin{equation}
(2\pi h)^{-1} \int {\rm e} (x_1, x_1;x_2, \xi_2, 0, \hbar) 
\bigl(1-{\bar\psi}_1(x_1)\bigr)\,dx_1 d\xi_2
\label{3-14}
\end{equation}
where $ {\rm e} (x_1, y_1;x_2, \xi_2, 0, \hbar)$ is the Schwartz kernel of 1-dimensional Schr\"odinger operator ${\bf a}(x_2,\xi_2,\hbar)$
obtained from ${\bar A}$ by change of variables $x_1\mapsto {\bar\gamma}x_1$ (and thus $\mu\mapsto 1$) and $hD_2\mapsto \xi_2$ while $e^\W (x_1;x_2,\xi_2,\hbar)$ would be   Weyl approximation for its restriction to diagonal.

However I would like to refer to a canonical operator ${\bf a}_0$ rather than to more general operator ${\bf a}$ and thus few extra steps are needed. I will delay them slightly.

\smallskip
\noindent
(ii) Consider now the second term corresponding to (\ref{3-11}); I remind it does not exceed $Ch\gamma |\log h|^2$ with superficial logarithmiic factor factor. I claim that

\begin{claim}\label{3-15}
The total contribution of the second term is 0 modulo 
$C\mu^{-1/\nu}h^{-1}+|\log h|^K$.
\end{claim}

Really, note first that one can assume that $\mu \le h^{-\nu}|\log h|^{-K}$, otherwise the total contribution $Ch^{-1}{\bar\gamma}_1|\log h|^2$ is less than announced. Let us apply rewrite  the term in question as
\begin{align}
u \mapsto {\bar u}_1=&-ih\sum_{\varsigma=\pm }{\bar G}^\varsigma  B_1
{\bar G}^\varsigma [{\bar A}, x_2-y_2] {\bar G}^\varsigma \delta(t)\delta (x_2-y_2)\delta(x_1-y_1)\label{3-16}\\
&-ih\sum_{\varsigma=\pm }{\bar G}^\varsigma [x_2-y_2, B_1] {\bar G}^\varsigma \delta(t)\delta (x_2-y_2)\delta(x_1-y_1).\notag
\end{align}
Let us apply to (\ref{3-16}) the method of successive approximations with unperturbed operator frozen as $x_1=y_1$; then each next term will acquire factor $C \mu \gamma^{\nu-1}h |\log h|\asymp 
(\gamma^/{\bar \gamma}_1)^{\nu-1}|\log h|$ and therefore the second term does not exceed $Ch^{-1}{\bar\gamma}|\log h|^{-K_1}$ as long as  
$\gamma \le {\bar\gamma}|\log h|^{K_2}$. On the other hand the leading term of is just a Weyl expression without any auxillary operators; this expression is odd with respect to $(\xi_2-\mu x_1^\nu/\nu)$ produces 0 after integration with respect to $\xi_2$. 

On the other hand, as $\gamma \ge C{\bar\gamma}|\log h|^{K_2}$ one can apply successive approximation method with unperturbed operator ${\hat A}$ is obtained from ${\bar A}$ by $\mu x_1\nu/\nu  \mapsto \xi_2 + \mu z_1^{\nu-1}(x_1-z_1)$,
$z_1= (\mu^{-1}\nu \xi_2)^{1/\nu}$ and the same is true for commutator 
$[{\bar A}, x_2-y_2]$ and $[B_1, x_2-y_2]$. Then since both $x_1$ and $y_1$ are confined to interval of the length $\beta = (\mu \gamma^{\nu-1})^{-1})$, each next term acquires factor 
$\beta \gamma^{-1}|\log h|\asymp \mu^{-1} h|\log h|\gamma^{-\nu}$ and therefore the second term does not exceed $C\mu^{-1}h^{-1}|\log h|^2 \gamma^{1-\nu}$ and 
its total contribution does not exceed what is announced as 
$\gamma \ge {\bar\gamma}|\log h|^{K_2}$. 

On the other hand the leading term corresponds to harmonic oscillator and obviously sums to $0$ as $x_1=y_1$ runs through $\cR$.

So I will keep $T=Ch|\log h|$ for a while. In the  analysis below I am going to consider (\ref{3-9}) with $\Gamma$ replaced by $\Gamma'$ which contains integration with respect to $x_1$ but not $x_2$; then one can assume with no loss of the generality that 

\begin{claim}\label{3-17}
$V$, $\sigma$, $\phi$ do not depend on $x_2$. 
\end{claim}

It is convenient to have $\phi=1$; one can always reach it now changing $x_1\mapsto x_1+O(x_1^2)$ and paying the price: now
\begin{equation}
{\bar A}= {\frac 1 2}\Bigl(\sigma_1(x_1)^2 (hD_1)^2 + 
\sigma (x_1)^2\bigl(hD_2 -\mu x_1^\nu/\nu \bigr)^2-V(x_1)\Bigr)
\label{3-18}
\end{equation}
with $\sigma=\sigma_1=1$ as $x_1=0$.

Let us apply the method of successive approximations again this time using as unperturbed operator ${\bar A}_0$ obtained from ${\bar A}$ by freezing $x_1=0$ in $V$, $\sigma$ and $\sigma_1$:
\begin{equation}
{\bar A}_0= {\frac 1 2}\Bigl((hD_1)^2 + \bigl(hD_2 -\mu x_1^\nu/\nu \bigr)^2-V(0)\Bigr).
\label{3-19}
\end{equation}
In these successive approximations the first term is still $O\bigl(h^{-2}{\bar\gamma}_1|\log h|\bigr)$  and an estimate for each subsequent term is gaining an extra factor $CT_0{\bar\gamma}_1 /h\asymp {\bar\gamma}_1 |\log h|$;  therefore $(n+1)$-th term does not exceed
$Ch^{-2}{\bar\gamma}_1^n |\log h|^n$ which is small for large $n$  since ${\bar\gamma}_1\le h^{(\nu-1)/(2\nu-1)}$ due to our assumption $\mu \ge h^{-\nu^2/(2\nu-1)}$. Note that we can write
\begin{equation*}
A \sim \sum_{k=0}^K x_1 ^k {\bar A}_k \qquad \text{with\ }  {\bar A}_k = {\bar A}_k(x_1, x_2,hD_1,hD_2;h,\mu)\quad \text{as\ }\ k\ge 1,
\end{equation*}
${\bar A}_0$ is already defined by (\ref{3-19}). Plugging it into successive approximation formula we get many terms but the first term  results in
\begin{equation}
(2\pi h)^{-1} \int {\rm e}_0 (x_1, x_1;x_2, \xi_2, 0, \hbar) \bigl(1-{\bar\psi}_1(x_1)\bigr)\,dx_1 d\xi_2
\label{3-20}
\end{equation}
(after finally I take $T=+\infty$ instead of $T=Ch|\log h|$) and the second term will be one for the approximation 
$A={\bar A}_0+x_1{\bar A}_1$. 

But then the contribution to  (\ref{3-14}) of this refined second term would be 0 because it is equal to (before integration over $\psi \,dx_2$)
\begin{equation}
\int \Tr_1 \bigl(b(\xi_2){\bar\psi}_1\bigr)\,d\xi_2
\label{3-21}
\end{equation}
where $b(\xi_2)$ is an operator in $\bH=L^2(\bR)$, $\Tr_1$ is the trace in this auxiliary space and $\cT b(\xi_2)\cT = -b (\xi'_2)$ with 
$\xi'_2=(-1)^\nu \xi_2$, $(\cT v)(x_1)=v(-x_1)$.

\smallskip
\noindent
(iii) Now let us consider terms in this successive approximation with $n\ge 3$; they will be still too large just to be skipped.  Let us apply to them the  successive approximations method  taking as unperturbed operator 
${\hat A}_0$ with $x_1^\nu/\nu$ replaced by $y_1^\nu/\nu + y_1^{\nu-1}(x-x_1)$:
\begin{equation}
{\hat A}_0= {\frac 1 2}\Bigl((hD_1)^2 + \bigl(hD_2 -
\mu \bigl(y_1^\nu/\nu + y_1^{\nu-1}(x-x_1)\bigr) \bigr)^2-V(0)\Bigr).
\label{3-22}
\end{equation}
One can see easily that each successive term of this approximation is gaining factor 
$\mu \gamma^{\nu-2} T^3/h \asymp \mu \gamma^{\nu-2}h^2 |\log h|^3$. Further, contributions of the second terms to (\ref{3-14}) will be equal 0 again after easy calculations. Finally the remainder will be then less than
\begin{equation*}
Ch^{-2}{\bar\gamma}_1^3 \times (\mu h^2 {\bar \gamma}_1^{\nu -2})^2|\log h|^K
\asymp C{\bar\gamma}_1 |\log h|^K\le 1
\end{equation*}
where the first factor is the estimate of the 3-rd term in the previous approximations and the second factor is what we gained now. So we can leave only the first terms in these last approximations.

Now when we sum all these terms with $n\ge 3$, we get
$\bigl(\cE^\MW -\cE_0^\MW\bigr)$ (multiplied by $\psi (1-\psi_1)$ and integrated in the end of the day) where
$\cE_0^\MW$ is what we would get for the first term if we ran last round of successive approximations for it; and this is exactly
$\cE^\MW_0(x_1,x_2,\tau,\mu,h)$.

However since we did not do it we have the same expression but with a correction term
\begin{equation}
(2\pi h)^{-1}\int {\rm e}_0 (x_1/{\bar\gamma}; \xi_2, \tau,\hbar)\, d\xi_2 - 
\cE_0^\MW (x_1,x_2,\tau,\mu,h).
\label{3-23}
\end{equation}

\smallskip
\noindent
(iv)  Consider now the second term given by (\ref{3-9}) with $u$ replaced by ${\bar u}_1$ defined by (\ref{3-11}). This term does not exceed 
$Ch^{-1}\gamma |\log h|^2$ (with a superficial logarithmic factor).
Again applying the successive approximations of (ii) one can see that replacing ${\bar G}^\varsigma$ by ${\bar G}_0^\varsigma$ brings the  extra factor 
$\gamma |\log h|^K$ and therefore all the terms but the first one do not exceed
$Ch^{-1}{\bar\gamma}_1^2 |\log h|^K$ which is less than the right-hand expression of of (\ref{3-5}) (unless in frames of the special case). Similarly one can replace $A_1$ by
\begin{equation*}
{\bar A}_1=-\mu (x_2-y_2)\beta(y_2)x_1^{\nu+1} \bigl( hD_2 - \mu x_1^\nu/\nu\bigr) -
{\frac 1 2}(x_2-y_2) w(y_2).
\end{equation*}

Or even let us consider the  original second term; I remind it does not exceed
$Ch^{-1}\gamma |\log h|$. Let

\smallskip
\noindent
(v) Now there is a special case when one needs to consider extra terms due to (\ref{3-12}), (\ref{3-13}) and the second term in approximation of (iv); all these terms do not exceed $Ch^{-1}{\bar\gamma}_1^2|\log h|^K$ and one needs only to ret rid off this logarithmic factor. 

Replacing ${\bar G}^\zeta$, 
${\bar G}_0^\zeta$ by ${\hat G}_0^\zeta$ linked to ${\hat A}_0$ defined by (\ref{3-22}) brings errors less than the right-hand expression of (\ref{3-5}). On the other hand calculations after this substitution show that these terms do not exceed $Ch^{-1}{\bar\gamma}_1^2$ (so there is no logarithmic factor) which even in the special case does not exceed $Ch^{-1}{\bar\gamma}$. Easy but tedious calculations I leave to the reader.
\end{proof}

\begin{proposition}\label{prop-3-4} (i) Let condition $(\ref{1-21})$ be fulfilled and let $\psi_1$ be an even function, $\psi=\psi(x_2)$. Then as 
$Ch^{-1}\le \mu \le h^{-\nu^2/(2\nu-1)}$
\begin{multline}
|(2\pi h)^{-1} \int {\rm e}_0(x_1, x_1;x_2, \xi_2, 0, \hbar) \bigl(1-{\bar\psi}_1(x_1)\bigr)\,dx_1 d\xi_2\\
-\int \cE_0^\MW(x_1,x_2, \xi_2, 0, \hbar)  \bigl(1-\psi_1(x_1)\bigr)\,dx_1| \le 
C{\bar\gamma}h^{-1}.
\label{3-24}
\end{multline}

\noindent 
(ii) Therefore estimate $(\ref{3-5})$ holds for $Ch^{-1}\le \mu \le h^{-\nu^2/(2\nu-1)}$ as well.
\end{proposition}

\begin{proof} (i) instantly follows from proposition \ref{prop-3-2} for 
$A={\bar A}_0$ and (ii) follows from proposition \ref{prop-3-2} and statement (i).\end{proof}

\subsection{Inner zone. I}

Now I need to consider an inner zone $\cZ_\inn=\{|x_1|\le {\bar\gamma}\}$. In this subsection I consider 
\begin{equation}
h^{-1}\int_{-\infty}^0 \Bigl( F_{t\to h^{-1}\tau}{\bar\chi}_T\Gamma\bigl(\psi \psi_1 u\bigr)\Bigr)\,d\tau
\label{3-25}
\end{equation}
with $T=Ch|\log h|$ even if this formula does not represent
$\Gamma (\psi\psi_1e)|_{\tau=0}$ with a desired precision, leaving the necessary correction for the rest of the section. In the end of this subsection 
I will list cases (due to section 2) when  provides the desired precision.

Running an absolutely standard method of successive approximations (i.e. taking as unperturbed operator with $x=y$) one can recover Weyl asymptotics with $C{\bar\gamma}h^{-2}\times\hbar ^2=C{\bar\gamma}^{-1}$ error with 
$\hbar= h/{\bar\gamma}$ while the contribution of the whole inner zone to the asymptotics is of magnitude $h^{-2}{\bar\gamma}$\,\footnote{\label{foot-15}It is easy to get rid off superficial logarithmic factors in this case.}. This $C{\bar\gamma}^{-1}$ error does not exceed $Ch^{-1}{\bar\gamma}$ as long as
${\bar\gamma}\ge \epsilon h^{\frac 1 2}$ which means exactly that
$\mu \le Ch^{-{\frac \nu 2}}$.

This is not extremely strong condition, much weaker than condition 
${\bar\gamma}\ge (h|\log h|)^{\frac 1 3}$ arising as 
$h^{-1}(\hbar|\log h|)^{\frac 1 2}\le h^{-1}{\bar\gamma}$ to provide remainder estimate $Ch^{-1}{\bar\gamma}$  for the spectral asymptotics in the general case. Moreover, it is even weaker a bit than condition 
${\bar\gamma}\ge (h|\log h|)^{\frac 1 2}$ arising as 
$h^{-1}(\hbar|\log h|)\le h^{-1}{\bar\gamma}$ to provide remainder estimate $Ch^{-1}{\bar\gamma}$ for the spectral asymptotics under condition $(\ref{2-103})$. 

And in the inner zone I can write a standard Weyl approximation or magnetic Weyl approximation to my sole discretion\footnote{\label{foot-16}Which would not be the case in the outer zone where we must use magnetic Weyl at least as $\mu \ge h^{-1}$.}. 

Note that both Weyl expressions are actually more precise in the inner zone than I have written before, giving an error not exceeding $Ch^{-2}{\bar\gamma}\times \hbar^4$. Really, when we write the complete decomposition in powers of $\mu $ and $h$ fixing ${\bar\gamma}\asymp {\bar\mu}^{-1/\nu}$ but considering $\mu\le {\bar\mu}$ as variable, we realize that terms with $\mu^{2k}h^{-2+2k}$ all are in magnetic Weyl approximation, terms without $\mu$ are produced by  Weyl expression for non-singular operator and thus are estimated by $Ch^{-2+2k}{\bar\gamma}$ rather than by $Ch^{-2}\hbar^{2k}{\bar\gamma}$ and unaccounted terms contain factor $\mu^{2k}h^{2l}$ with $l\ge k\ge 1$ and thus an error term does not exceed 
$Ch^{-2}\hbar ^4{\bar\gamma}=Ch^2{\bar\gamma}^{-3}$. 

Note that $Ch^2{\bar\gamma}^{-3}\le Ch^{-1}{\bar\gamma}$ as ${\bar\gamma}\ge h^{3/4}$ or equivalently $\mu \le Ch^{-3\nu/4}$. Further, we get 
\begin{equation*}
\int {\frac T h}{\hat{\bar\chi}}\bigl(-{\frac T h}\tau\bigr)\int \cE^\MW (x,0)\psi_1(x)\psi(x)\,dx\,d\tau
\end{equation*}
which can be rewritten as (\ref{3-26}) with $O(1)$ error.
Therefore we arrive to

\begin{proposition}\label{prop-3-5}\footnote{\label{foot-17}Cf. proposition \ref{prop-3-2}.}  With an error $Ch^2{\bar\gamma}^{-3}$ one can rewrite $(\ref{3-25})$ with $T=Ch|\log h|$ as
\begin{equation}
\int \cE^\MW (x,0)\psi_1(x)\psi (x)\,dx.
\label{3-26}
\end{equation}
\end{proposition}

Therefore due to subsection 2.3-2.5 we arrive to 

\begin{corollary}\label{cor-3-6} (i) In the general case the contribution of the inner zone to the spectral remainder estimate with the principal part $(\ref{3-26})$
does not exceed 
$Ch^{-1}{\bar\gamma}+Ch^{-1}(h|\log h|/{\bar\gamma})^{\frac 1 2}$;
in particular as $\mu\le (h|\log h|)^{-{\frac \nu 3}}$ the remainder estimate is $Ch^{-1}{\bar\gamma}$;

\smallskip
\noindent
(ii) Under conditions $(\ref{1-21})$, $(\ref{2-103})$  the contribution of the inner zone to the remainder estimate with the principal part $(\ref{3-26})$ does not exceed  
$Ch^{-1}{\bar\gamma}+C |\log h|{\bar\gamma}^{-1}$; in particular as $\mu\le (h|\log h|)^{-{\frac \nu 2}}$ this remainder estimate is $Ch^{-1}{\bar\gamma}$;

\smallskip
\noindent
(iii) Under conditions $(\ref{1-21})$, $(\ref{2-97})$ with 
$\zeta \ge Ch |\log h|/{\bar\gamma}+C{\bar\gamma}$ the contribution of the inner zone to the remainder estimate with the principal part $(\ref{3-26})$ does not exceed  
$Ch^{-1}{\bar\gamma}+ Ch^2{\bar\gamma}^{-3}$; in particular as $\mu\le h^{-{\frac 3 4}\nu}$ this remainder estimate is $Ch^{-1}{\bar\gamma}$.
\end{corollary}

Then combining with the result of the outer zone analysis we arrive to

\begin{theorem}\label{thm-3-7} Let condition $(\ref{2-8})$  be fulfilled. Then

\smallskip
\noindent
(i) In the general case
$\cR$ defined by $(\ref{0-12})$ does not exceed 
$Ch^{-1}{\bar\gamma}+Ch^{-1}(h|\log h|/{\bar\gamma})^{\frac 1 2}$; in particular $\cR$ does not exceed $Ch^{-1}{\bar\gamma}$ as $\mu \le C(h|\log h|)^{-\nu/3}$;

\smallskip
\noindent
(ii) Under conditions  $(\ref{1-21})$ and  $(\ref{2-103})$  $\cR$ does not exceed 
$Ch^{-1}{\bar\gamma}+C{\bar\gamma}|\log h|$; in particular $\cR$ does not exceed $Ch^{-1}{\bar\gamma}$ as $\mu \le C(h|\log h|)^{-\nu/2}$;

\smallskip
\noindent
(iii) Under conditions $(\ref{1-21})$ and $(\ref{2-97})$ with 
$\zeta \ge Ch |\log h|/{\bar\gamma}+C{\bar\gamma}$ $\cR$ does not exceed 
$Ch^{-1}{\bar\gamma}+C(\mu h+1)^{1/(\nu -1)}$; in particular $\cR$ does not exceed $Ch^{-1}{\bar\gamma}$ as $\mu \le Ch^{-\nu^2/(2\nu-1)}$.
\end{theorem}

On the other hand repeating (almost all) arguments of the proof of proposition 
\ref{prop-3-3} one can prove easily

\begin{proposition} \label{prop-3-8} In frames of proposition \ref{prop-3-3}
\begin{multline}
|h^{-1}\int_{-\infty}^0 
\Bigl(F_{t\to h^{-1}\tau}{\bar\chi}_T(t)\Gamma 
\bigl(\psi \psi_1u\bigr)\Bigr)\,d\tau -\\
\int \Bigl( \cE^\MW(x,0)-\cE_0^\MW(x,0)\Bigr) \psi (x_2)\psi_1(x_1)\,dx -\\
(2\pi h)^{-1}\int {\frac T h} {\hat{\bar\chi}}\bigl(-{\frac T h}\tau)\bigr)
\int {\rm e}_0(x_1, x_1;x_2, \xi_2, \tau, \hbar) \psi (x_2) {\bar\psi}(x_1)\,dx d\xi_2\,d\tau| \le \\
Ch^{-1}{\bar\gamma}+C|\log h|^K
\label{3-27}
\end{multline}
\end{proposition}

Combining with the results of the previous subsection we arrive to

\begin{theorem}\label{thm-3-9} Let conditions $(\ref{1-21})$  and $(\ref{2-97})$ with $\zeta \ge Ch |\log h|/{\bar\gamma}+C{\bar\gamma}$ be fulfilled, $\psi=\psi(x_2)$. Then for $h^{-\nu^2/(2\nu-1)}\le \mu \le h^{-\nu}|\log h|^{-K}$
\begin{equation}
|\int \Bigl( e(x,x,0)-
 \cE^\MW(x,0)\Bigr) \psi (x )\,dx - 
 \int  \cE_\corr^\MW(x_2,0)\psi (x_2)\,dx_2|\le 
Ch^{-1}{\bar\gamma}+C|\log h|^K
\label{3-28}
\end{equation}
with
\begin{equation}
\cE^\MW_\corr =
 (2\pi h)^{-1}
\int {\rm e}_0(x_1, x_1;x_2, \xi_2, \tau, \hbar) \,dx_1 d\xi_2-
\int \cE^\MW_0 (x,\tau)\,dx_1.
\label{3-29}
\end{equation}
\end{theorem}

\subsection{Inner zone. II}

In section 2  we derived asymptotics with the principal part which is the sum of 
\begin{equation}
\sum_m h^{-1}\int_{-\infty}^0 \Bigl(F_{t\to h^{-1}\tau} 
{\bar\chi}_{T_m}(t)
\Gamma (\psi_1 \psi Q_m u)\Bigr) \,d\tau
\label{3-30}
\end{equation}
where $Q_m$ are elements of partition (in $hD_2$ and may be $x_2$) and 
$T_m = C h |\log h|  /\varrho_m^2$. We can rewrite (\ref{3-30}) as the sum of
(\ref{3-25}) with $T={\bar T}=\epsilon{\bar\gamma}$ and 
\begin{equation}
\sum_m h^{-1}\int_{-\infty}^0 \Bigl(F_{t\to h^{-1}\tau} 
\bigl({\bar\chi}_{T_m}(t)- {\bar\chi}_{\bar T}(t)\bigr)
\Gamma (\psi_1 \psi Q_m u)\Bigr) \,d\tau.
\label{3-31}
\end{equation}

Consider now (\ref{3-25}) with $T=\epsilon{\bar\gamma}$. Also due to propagation results of section 2 one can replace there $T=\epsilon{\bar\gamma}$  by $T=Ch|\log h|$; then all the results of the
previous subsection are applicable. It includes  (\ref{3-27}) as well where one can replace back $T=Ch|\log h|$ by $T=\epsilon {\bar\gamma}$ in  the first and/or the last term of the left-hand expression.

Moreover again due to propagation results of section 2, all terms in (\ref{3-31}) with $\varrho_m \ge C{\bar\rho}_1$ are negligible and therefore one can rewrite (modulo negligible) expression (\ref{3-30}) as the sum of (\ref{3-25}) with $T=Ch|\log h|$ and (\ref{3-31}) with ${\bar T}=\epsilon{\bar\gamma}$ with summation with respect to $m$ such that $\varrho_m \le C{\bar\rho}_1$.

\begin{proposition}\label{prop-3-10} (i) With an error not exceeding $Ch^{-1}{\bar\gamma}+C{\bar\gamma}_1^{-1}$ one can rewrite $(\ref{3-30})$ as
the sum of $(\ref{3-26})$ and $(\ref{3-31})$ with ${\bar T}=\epsilon {\bar\gamma}$ and $Q_m$ partition of unity in the periodic zone only.

\smallskip
\noindent
(ii) In particular, as $\mu \le Ch^{-\nu^2/(2\nu-1)}$ this error  does not exceed $Ch^{-1}{\bar\gamma}$.
\end{proposition}

As $\mu \ge Ch^{-\nu^2/(2\nu-1)}$ I need a bit more complicated analysis:

\begin{proposition}\label{prop-3-11} (i) With an error not exceeding $Ch^{-1}{\bar\gamma}+C|\log h|^K$ one can rewrite $(\ref{3-30})$ as
the sum of the following three expressions:
\begin{multline}
\int \Bigl( \cE^\MW(x,0)-\cE_0^\MW(x,0)\Bigr) \psi (x_2)\psi_1(x_1)\,dx +\\
\int {\rm e}_0(x_1, x_1;x_2, \xi_2, 0, \hbar) \psi (x_2) {\bar\psi}(x_1)\,dx d\xi_2
\label{3-32}
\end{multline}
where the first term is actually less than the remainder estimate,
\begin{equation}
(2\pi h)^{-1}\int \Bigl({\frac T h} {\hat{\bar\chi}}\bigl(-{\frac T h}\tau)\bigr)-\delta (\tau)\Bigr)
\int {\rm e}_0(x_1, x_1;x_2, \xi_2, \tau, \hbar) \psi (x_2) {\bar\psi}(x_1)\,dx d\xi_2\,d\tau 
\label{3-33}
\end{equation}
with $T={\bar T}=\epsilon {\bar\gamma}$, Dirac $\delta$-function \underbar{and}  $(\ref{3-31})$ where now  $\{Q_m\}$ is a partition of unity in the periodic zone only;

\smallskip
\noindent
(ii) Further, one can rewrite $(\ref{3-33})$ as
\begin{multline}
\sum_m (2\pi h)^{-1}\int \Bigl({\frac T h} {\hat{\bar\chi}}\bigl(-{\frac T h}\tau)\bigr)-\delta (\tau)\Bigr)\times \\
\int {\rm e}_0(x_1, x_1;x_2, \xi_2, \tau, \hbar) \psi (x_2) {\bar\psi}(x_1)Q_m(x_2,\xi_2)\,dx d\xi_2\,d\tau 
\label{3-34}
\end{multline}
where now  $\{Q_m\}$ is a partition of unity in the periodic zone only.
\end{proposition} 

\begin{proof} Proposition follows from proposition \ref{prop-3-8} and propagation results of section 2 with the only exception that in (\ref{3-34}) summation is taken over all $m$, and one needs to prove that the corresponding terms are negligible in non-periodic zone where condition
\begin{equation}
\bigl(|\xi_2-V^{\frac 1 2}k^*_\hbar |+|\partial_{x_2}V|\bigr)\asymp \rho \le C{\bar\rho}_1
\label{3-35}
\end{equation}
is violated.

Note that the expression in question is the sum of
\begin{multline}
(2\pi T)^{-1}{\hat\chi} \bigl(-{\frac h T}\tau\bigr)\int {\rm e}_0(x_1, x_1;x_2, \xi_2, \tau, \hbar) \psi (x_2) {\bar\psi}(x_1)Q_m(x_2,\xi_2)\,dx d\xi_2\,d\tau =\\
\int h^{-1}\int_{-\infty}^0 \Bigr(F_{t\to h^{-1}\tau} \chi_T(t) \bigl(\Gamma'\psi_1(x_1) Q_m (y_2,hD_2) {\bar u}_0\bigr)\Bigr) \psi(y_2)\,dy_2
\label{3-36}
\end{multline}
with $T$ running from $\epsilon{\bar\gamma}$ to $+\infty$ and ${\bar u}_0$
defined for operator ${\bar A}_0$.

If $|\xi_2-V^{\frac 1 2}k^*_\hbar |\ge \rho$ on partition element $Q_m$ then (\ref{3-36}) does not exceed $Ch^s(T+1)^{-s}$  as $T\ge C_0h|\log h|/\rho^2$ just due to propagation results of section 2 applied to ${\bar u}_0$ and ${\bar A}_0$ and since for $C_0 h|\log h|/\rho^2 \le \epsilon {\bar\gamma}$ as 
$\rho \ge {\bar\rho}_1$, these elements are covered.

On the other hand if $|\partial_{x_2}V|\asymp \rho$ on the partition element in question, one can notice that ${\bar A}_0$ depends on $y_2$ via 
$-{\frac 1 2}W(y_2)$ only and introducing new variable $y'_2=-{\frac 1 2}W(y_2)$
we get mollification with respect to the spectral parameter and thus (\ref{3-36}) does not exceed $Ch^s(T+1)^{-s}$ as well. I leave details to the reader.
\end{proof}

\subsection{Periodic and near-periodic orbits. I. Pilot-model}

Thus I need to analyze periodic zone $\cZ_\per$ more accurately defining it by (\ref{3-35}) rather than by 
\begin{equation}
|\xi_2-V^{\frac 1 2}k^*_\hbar|\asymp \rho \le C{\bar\rho}_1
\label{3-37}
\end{equation}
in the case when magnetic field is strong enough to prevent corollary \ref{cor-3-6} from presenting a sharp remainder estimate.

Namely I need to consider  term
\begin{equation}
h^{-1}\int_{-\infty}^0 \Bigl(F_{t\to h^{-1}\tau} 
\bigl({\bar\chi}_T(t)- {\bar\chi}_{\bar T}(t)\bigr)
\Gamma (\psi_1 \psi Q u)\Bigr) \,d\tau
\label{3-38}
\end{equation}
with partition element $Q=Q(x_2,hD_2)$, $T=T_1$ where $T_0$ and $T_1$ are defined in subsections 2.8-2.9 and  ${\bar T}=C_0h|\log h|$. Here due to propagation results of these subsections I can take any $T\in [T_0,T_1]$ to my discretion, in particular $T=T_0$.

Let consider a pilot-model first; more precisely let us assume that operator under consideration coincides with a pilot-model in $B(0,1)$. Then 
$T_0= C\rho^{-2}h|\log h|$ and $T_1=\epsilon \rho$ with 
$|\xi_2-k^*_\hbar|\asymp \rho$
on $\supp Q$.

Then (\ref{3-38}) will be the same (modulo negligible) if one replaces $u$ by ${\bar u}$ constructed for operator coinciding with pilot-model in $\bR^2$. But then one can replace $T=\epsilon \rho^{-1}$ by $T=\infty$. Also one can remove restriction from below on $\rho: \rho \ge {\bar\gamma}$. 

Also in this case one can replace ${\bar\psi}_1$ by $1$ since contribution to (\ref{3-38}) of the outer  zone is negligible due to propagation results there and contribution of the forbidden zone is negligible as well. Further, one does not need to integrate over $\psi(x_2)dx_2$ in (\ref{3-38}) anymore.

After all these modifications there is no need to consider a partition and then correction (\ref{3-38}) is transformed to  
\begin{equation}
 h^{-1}\int_{-\infty}^0 \Bigl(F_{t\to h^{-1}\tau} 
\bigl(1- {\bar\chi}_{\bar T}(t)\bigr)\Gamma'  {\bar u}\Bigr)\, d\tau
\label{3-39}
\end{equation}
where I remind $\Gamma' v= \int v(x_1,x_2;x_1,x_2)\,dx_1$; in our case $\cE^\MW_\corr$ does not depend on $x_2$. 

Obviously
\begin{equation}
h^{-1}\int_{-\infty}^0 \Bigl(F_{t\to h^{-1}\tau} 
\Gamma'  {\bar u}\Bigr)\, d\tau=(2\pi h)^{-1}\int  {\bf n}_0(\xi_2,\hbar)\,d\xi_2
\label{3-40}
\end{equation}
with ${\bf n}_0(\xi_2,\hbar)$ the number of negative eigenvalues of operator 
${\bf a}={\bf a}_0(\xi_2,\hbar)$  with $W=1$. On the other hand, 
\begin{equation}
h^{-1}\int_{-\infty}^0 \Bigl(F_{t\to h^{-1}\tau} {\bar\chi}_{\bar T}(t)
 \Gamma'  {\bar u}\Bigr)\, d\tau \equiv\int \cE_0^\MW (x,0)\, dx_1 \qquad
\mod O\bigl({\bar\gamma}_1^{-1}\bigr)
\label{3-41}
\end{equation}
and therefore
\begin{equation}
\cE^\MW_\corr\equiv (2\pi h)^{-1}\int  {\bf n}_0(\xi_2,\hbar)\,d\xi_2 - 
\int \cE_0^\MW (x,0)\, dx_1.
\label{3-42}
\end{equation}
The error in (\ref{3-42}) does not exceed $Ch^{-1}{\bar\gamma}$ even if 
${\bar\gamma}_1^{-1}\ge h^{-1}{\bar\gamma}$ because in this case I already made
an error replacing left-hand expression of (\ref{3-41}) by its right-hand-expression in the analysis of subsection 3.1 and now I just compensated  it by adding the skipped term. One can consider (\ref{3-42}) as a definition of the correction term for the model operator. Thus we arrive to

\begin{proposition}\label{prop-3-12} For operator coinciding with pilot model in 
$(-1,1)\times \bR$
\begin{equation}
\int e(x,x,0)\,dx_1 \equiv \int \cE^\MW_0\,dx_1 +\cE^\MW_\corr
\qquad \mod O\bigl(h^{-1}{\bar\gamma}\bigr)\label{3-43}
\end{equation}
in $(-{\frac 1 2},{\frac 1 2})$ where $\cE^\MW_\corr$ is defined by $(\ref{3-42})$.
\end{proposition}

Expression (\ref{3-42}), multiplied by $\psi(x_2)dx_2$ and integrated will be used in two next subsections in the general case as well. However to get more explicit even if less precise expression one can replace  (\ref{3-39})  by 
\begin{equation}
 h^{-1}\int_{-\infty}^0 \Bigl(F_{t\to h^{-1}\tau} 
\bigl(1- {\bar\chi}_{\bar T}(t)\bigr)\Gamma'  \phi (hD_2){\bar u}\Bigr)\, d\tau
\label{3-44}
\end{equation}
with $\phi $ supported in $[-2C_0,2C_0]$ and equal 1 in $[-C_0,C-0]$, making 
$O\bigl({\bar\gamma}_1^{-1}\bigr)$ error.

Using
\begin{align}
&h^{-1}\int_{-\infty}^0 \Bigl(F_{t\to h^{-1}\tau} 
\Gamma'  \phi(hD_2){\bar u}\Bigr)\, d\tau=
(2\pi h)^{-1}\int  {\bf n}_0(\xi_2,\hbar)\phi(\xi_2)\,d\xi_2,\label{3-45}\\
&h^{-1}\int_{-\infty}^0 \Bigl(F_{t\to h^{-1}\tau} {\bar\chi}_{\bar T}
\Gamma'  \phi(hD_2){\bar u}\Bigr)\, d\tau\equiv
(2\pi h)^{-1}\int  n^\W_0(\xi_2,\hbar)\phi(\xi_2)\,d\xi_2
\label{3-46}
\end{align}
with the second equality modulo $O\bigl({\bar\gamma}^{-3}h^2\bigr)$ with 
\begin{equation}
n_0^\W (\xi_2,\hbar)=(\pi \hbar)^{-1}\int \Bigl(1-\bigl(\xi_2-x_1^\nu/\nu\bigr)^2\Bigr)_+^{\frac 1 2}\,dx_1
\label{3-47}
\end{equation}
the Weyl approximation of ${\bf n}_0(\xi_2,\hbar)$, we arrive to

\begin{proposition}\label{prop-3-13} For operator coinciding with pilot model in 
$(-1,1)\times \bR$
\begin{equation}
\cE^\MW_\corr \equiv (2\pi h)^{-1} \int \Bigl( {\bf n}_0 (\xi_2,\hbar)- 
n^\W_0 (\xi_2,\hbar)\Bigr)\phi (\xi_2)\,d\xi_2
\qquad \mod O\bigl({\bar\gamma}_1^{-1}\bigr)
\label{3-48}
\end{equation}
in $(-{\frac 1 2},{\frac 1 2})$.
\end{proposition}

To calculate (\ref{3-48}) more explicitly I need to calculate eigenvalues $\lambda_n(\xi_2,\hbar)$ of operator ${\bf a}_0(\xi_2,\hbar)$ and I am interested in those eigenvalues which are close to 0 as $\xi_2$ is close to $k^*_\hbar$.

These eigenvalues are defined modulo $O(\hbar^2)$ from Bohr-Sommerfeld condition
\begin{align}
&{\frac 1 {2\pi \hbar} }S(\xi_2,\tau) + {\frac 1 4}\iota_ M \bigl(\cL(\xi_2,\tau)\bigr)\in \bZ\qquad \text{as\ } \tau=\lambda_n\label{3-49}\\
\intertext{where}
&S(\xi_2,\lambda)= \oint_{\cL(\xi_2,\lambda)} \xi_1\,dx_1
\label{3-50}
\end{align}
and $\iota_M \bigl(\cL(\xi_2,\tau)\bigr)$ is Maslov' index of the closed  trajectory $\cL(\xi_2,\tau)$ (on 2-dimensional phase plane) on the energy level $\tau$; $\iota_M \bigl(\cL(\xi_2,\tau)\bigr)=2$ as $\tau\approx 0$. Since
$\partial_\tau S(\xi_2, \tau) |=T (\xi_2,\tau)$ due to Hamiltonian mechanics, 
the spacing between two consecutive eigenvalues is $2\pi \hbar/T(\xi_2,\tau) +O(\hbar^2)$. Further, all these eigenvalues $\lambda_n(\xi_2,\hbar)$ are uniformly analytic functions of 
$\xi_2, |\xi_2|\ll 1$, $\hbar\ll 1$ as $|\lambda_n |\ll 1$ and also
$\partial_{\xi_2}^2 \lambda_n\ge \epsilon_0$. This yields immediately

\begin{proposition}\label{prop-3-14} As $\mu \le \epsilon (h|\log h|)^{-\nu}$

\smallskip
\noindent 
(i) $|\cE^\MW_\corr |\le Ch^{-1}\hbar^{\frac 1 2}$; in particular, 
$|\cE^\MW_\corr |\le Ch^{-1}{\bar\gamma}$ as $\mu \le Ch^{-\nu/3}$;

\smallskip
\noindent 
(ii) One can calculate $\cE^\MW_\corr$ with an error not exceeding $Ch^{-1}\hbar$ (which does not exceed $Ch^{-1}{\bar\gamma}$ as $\mu \le Ch^{-\nu/2}$) replacing $\lambda_n(\xi_2,\hbar)$ by its approximation  modulo 
$O\bigl( \eta^4+\hbar^2\bigr)$  where here and below $\eta = \xi_2-k^*_\hbar$;

\smallskip
\noindent 
(iii) One can calculate $\cE^\MW_\corr$ with an error not exceeding $Ch^{-1}\hbar^{3/2}$ (which does not exceed 
$C\bigr(h^{-1}{\bar\gamma}+ {\bar\gamma}_1^{-1}\bigl)$ as $\nu \ge 3$) replacing $\lambda_n(\xi_2,\hbar)$ by its approximation modulo  $O\bigl( \eta^6+\hbar^3\bigr)$;

\smallskip
\noindent 
(iv) As $\nu=2$ one can calculate $\cE^\MW_\corr$ with an error not exceeding $Ch^{-1}\hbar^2$ (which does not exceed $C{\bar\gamma}_1^{-1})$) replacing $\lambda_n(\xi_2,\hbar)$ by its approximation  modulo  $O\bigl( \eta^8+\hbar^8\bigr)$.
\end{proposition}

To exploit (ii) one can replace $\lambda_n(\xi_2,\hbar)$ by its approximate value 
\begin{equation}
\lambda_n(\xi_2,\hbar) \equiv {\bar\lambda}(\xi_2) + {\frac {(2 n+1)\pi\hbar} {T (\xi_2)}}.
\label{3-51}
\end{equation}
  Here ${\bar\lambda}(\xi_2)$ is a  ``classical eigenvalue'' corresponding to 
$n=0$ and Maslov index 0 (instead of 2) and one can find it  from (\ref{3-49}) as $\tau=0$;  then 
${\bar\lambda}(\xi_2)= -S(\xi_2, \tau)/\partial_\tau S\bigr|_{\tau=0}$.

I remind that $(\partial_\tau S)\bigr|_{\tau=0}=T(\xi_2)$,  
$\partial_{\xi_2} S =I(\xi_2)$ with $I(\xi_2)$ defined by (\ref{1-4}) for even $\nu$ and similar formula for odd $\nu$; further, 
$I(k^*)=0$, $\partial_{\xi_2} I(k^*)=\kappa >0$. Then
\begin{align}
\cE^\MW_\corr\equiv &(2\pi h)^{-1}\int \Bigl(
{\frac 1 {2\pi\hbar}}  S(\xi_2) - \Bigl\lfloor {\frac 1 {2\pi\hbar}}  S(\xi_2)+{\frac 1 2} \Bigr\rfloor
\Bigr)\,d\xi_2
\equiv 
\label{3-52}\\
&(2\pi h)^{-1}\int \Bigl(
{\frac 1 {2\pi\hbar}}  \bigl(S_0 + {\frac 1 2}\kappa \eta^2\bigr)
- \Bigl\lfloor {\frac 1 {2\pi\hbar}}  \bigl(S_0 + {\frac 1 2}\kappa \eta^2\bigr) +{\frac 1 2} \Bigr\rfloor\Bigr)\,d\eta =
\notag\\
&(2\pi)^{-{\frac 3 2}} h^{-1}\hbar^{\frac 1 2}\kappa^{-{\frac 1 2}} 
G\bigl({\frac {S_0}{2\pi\hbar}}\bigr)\qquad \mod O(h^{-1}\hbar)\notag
\end{align}
with function $G$ defined by
\begin{equation}
G(t) = \int_\bR \Bigl(t+{\frac 1 2}\eta^2 - \bigl\lfloor t+{\frac 1 2}\eta^2 +{\frac 1 2}\bigr\rfloor \Bigr)\,d\eta
\label{3-53}
\end{equation}
with the converging integral in the right-hand expression.

One can prove easily that
\begin{equation}
G\not\equiv 0,\quad G(t+1)=G(t),\qquad\int_0^1G(t)\,dt=0,\qquad 
G \in C^{\frac 1 2}.
\label{3-54}
\end{equation}

To exploit proposition  \ref{prop-3-14}(iii),(iv) one can use more precise version of (\ref{3-51})
\begin{equation}
\lambda_n(\xi_2,\hbar) \equiv {\bar\lambda}(\xi_2) + {\frac {(2 n+1)\pi\hbar} {T (\xi_2)}} + \varkappa_1 (\xi_2 )\hbar^2 + \varkappa_2 (\xi_2 )\hbar^3
\label{3-55}
\end{equation}
and after obvious calculations (\ref{3-52}) is adjusted to
\begin{align}
&\cE^\MW_\corr \equiv 
(2\pi)^{-{\frac 3 2}} h^{-1}\hbar^{\frac 1 2}\kappa^{-{\frac 1 2}} 
G\bigl({\frac {S_0}{2\pi\hbar}}+\kappa _1\hbar\bigr)
 \qquad\qquad\mod O\bigl(h^{-1}\hbar^{\frac 3 2}\bigr), \tag*{$(3.52)^*$}\label{3-52-*}\\
&\cE^\MW_\corr \equiv (2\pi)^{-{\frac 3 2}} h^{-1}\hbar^{\frac 1 2}\kappa^{-{\frac 1 2}} 
\Bigl( G(t)+ \kappa_3 \hbar G(t) + \kappa_4 \hbar G_1(t)\Bigr)
\Bigr|_{\displaystyle t={\frac {S_0}{2\pi\hbar}}+\kappa _1\hbar+\kappa_2\hbar^2}\tag*{$(3.52)^{**}$}\label{3-52-**}\\
&\hskip300pt\mod 
O\bigl(h^{-1}\hbar^2\bigr)
\notag
\end{align}
with some constants $\kappa_1,\dots,\kappa_4$ and 
\begin{equation}
G_1(t) =\int_0^t G(t')\,dt' - \int_0^1 (1-t')G(t')\,dt',
\label{3-56}
\end{equation}
$G_1\in C^{3/2}$ satisfying (\ref{3-54}). Thus we arrive to

\begin{proposition}\label{prop-3-15}
For operator coinciding in $B(0,1)$ with pilot-model with potential $V=1$ and fixed $\psi$, $\psi_1$ ($\psi_1=1$ as $|x_1|\le \epsilon$) in asymptotics 
$(\ref{3-47})$ 

\smallskip
\noindent
(i) Modulo $O\bigl(h^{-1}\hbar^{1/2}\bigr)$, $O\bigl(h^{-1}\hbar\bigr)$, $O\bigl(h^{-1}\hbar^{3/2}\bigr)$ and $O\bigl(h^{-1}\hbar^2\bigr)$ one can respectively skip correction term $\cE^\MW_\corr$ or to define it by $(\ref{3-52})$,  \ref{3-52-*}, \ref{3-52-**};

\smallskip
\noindent
(ii) In particular, uncorrected asymptotics has sharp remainder estimate if and only if $\hbar^{\frac 1 2}\le C{\bar\gamma}$ i.e. $\mu \le Ch^{-\nu/3}$. Further, one can skip define $\cE^\MW_\corr$ by $(\ref{3-52})$  without spoiling remainder estimate as $Ch^{-\nu/3}\le \mu \le Ch^{-\nu/2}$.
\end{proposition}

\begin{remark}\label{rem-3-16}
(i) One can then generalize proposition \ref{prop-3-15} immediately to a  pilot-model with potential $V=W=\const$ instead of $V=1$ by replacing $h\mapsto h W^{-1/2}$, 
$\mu \mapsto \mu W^{-1/2}$ and therefore 
$\hbar \mapsto \hbar W^{-(\nu+1)/2\nu}$ just modifying (\ref{3-52}) to
\begin{equation}
\cE^\MW_\corr \equiv 
(2\pi)^{-{\frac 3 2}}h^{-1} \hbar^{\frac 1 2}\kappa^{-{\frac 1 2}}
W^{\frac {\nu - 1}{4\nu}}
G\Bigl({\frac {S_0W^{\frac {\nu +1}{2\nu}}}{2\pi\hbar} }\Bigr)
\label{3-57}
\end{equation}
and similarly modifying \ref{3-52-*}, \ref{3-52-**}.

\smallskip
\noindent
(ii) One should expect the similar correction for the general operators as well but one can see easily that under condition $|\partial_{x_2}W|\ge \epsilon$ integrated correction term 
\begin{equation}
\int \cE^\MW_\corr (x)\psi \,dx_2
\label{3-58}
\end{equation}
is negligible, and under condition (\ref{2-103}) it should not exceed $Ch^{-1}\hbar$ (and be of this magnitude if there are critical points of $W$) with integrated versions of (\ref{3-52}), \ref{3-52-*} valid modulo $O(h^{-1}\hbar^{\frac 3 2})$, $O(h^{-1}\hbar^2)$ respectively.
\end{remark}

\subsection{Periodic orbits. II. General settings}

Now let us consider inner zone in frames of subsections 2.8--2.9 when we managed to prove sharp remainder estimate $C{\bar\gamma}h^{-1}$ as 
$\mu \le h^{\delta -\nu}$ producing the final result as the sum of expressions (\ref{3-30})
taken over partitions $\psi=\psi_m(x_2)$ and $Q=Q_m(hD_2)$ of unity; I remind that $\psi_m$ is $\ell$-admissible and $Q_m$ is $\rho$-admissible
with $\ell $, $\rho$ satisfying some conditions the most important of which are (\ref{2-70})-(\ref{2-72}); here $\ell$ denotes the function associated with the finest of subpartitions and we assume that on this element 
\begin{equation}
|\nabla V|\le C_0\zeta, \qquad |\xi_2- k^*_\hbar V(0,x_2)|\le C_0\rho
\label{3-59}
\end{equation}
with one of these inequalities being reversible (with $\epsilon_0$ instead of $C_0$);  then one can take in (\ref{3-30}) any $T$ ranging from 
\begin{equation}
T_0= Ch|\log h| \min \bigl({\frac 1 {\rho^2}}, {\frac 1{\ell\zeta}}\bigr) =
C\epsilon^{-1}h|\log h| \times T_1
\label{3-60}
\end{equation}
to $T_1$.

There are also exceptional elements with $T={\bar T}=Ch|\log h|$ but they will be treated easily in the same manner. 

Therefore I need to consider expression $(\ref{3-30})$ with 
$T={\bar T}=Ch|\log h|$ and also correction terms (\ref{3-38}) with arbitrarily $T\in [T_0,T_1]$ depending on the partition element; index $m$ indicating partition element I am skipping. 

One can replace then 
$\bigl({\bar\chi}_T(t)- {\bar\chi}_{\bar T}(t)\bigr)$ by the sum of $\chi_T(t)$
with $T$ running from ${\bar T}$ to $T_0$ and with $\chi$ supported in $(-1,-{\frac 1 2})\cup ({\frac 1 2},1)$ thus replacing (\ref{3-38}) with the sum of terms
\begin{equation}
h^{-1}\int_{-\infty}^0 \Bigl(F_{t\to h^{-1}\tau} 
\chi_T(t)
\Gamma (\psi_1 \psi Q u)\Bigr) \,d\tau = iT^{-1}
F_{t\to h^{-1}\tau} 
{\check\chi}_T(t) \Gamma (\psi_1 \psi Q u)
\label{3-61}
\end{equation}
with ${\check\chi}(t)=t^{-1}\chi(t)$ and $T\in [{\bar T}, T_0]$.

Further, one needs to consider only $T$ ranging from 
${\bar T}'=\epsilon {\bar\gamma}$ to $T_0$ because terms (\ref{3-61}) with $T$ ranging from ${\bar T}$ to ${\bar T}'$ are negligible; therefore partition elements with 
$T_0\le {\bar T}'$\,\footnote{\label{foot-18} Or equivalently with 
$(\zeta\ell +\rho^2){\bar\gamma}\ge Ch|\log h|$.} will be eliminated.

\begin{remark}\label{rem-3-17}
(i) Obviously (\ref{3-61}) does not exceed $Ch^{-2}{\bar\gamma}\rho\ell$ and thus   (\ref{3-38}) does not exceed $Ch^{-2}{\bar\gamma}\rho\ell|\log h|$;

\smallskip
\noindent
(ii) Further, one can decompose (\ref{3-61}) in the same way as in the proof of proposition \ref{prop-2-30}; therefore (\ref{3-61}) does not exceed 
$Ch^{-2}{\bar\gamma}\rho\ell\times h|\log h|/{\bar\gamma}= Ch^{-1}\rho\ell\times |\log h|$ and thus (\ref{3-38}) does not exceed $Ch^{-1}\rho\ell|\log h|^2$. 

Even if this is the remainder estimate obtained by taking $T=\epsilon{\bar\gamma}$ in the Tauberian arguments (with a superficial logarithmic factor), it provides a solid ground for the further estimates.
\end{remark}

Let us consider regular elements and apply successive approximation method freezing in unperturbed operator $x_2=y_2$ as usual. However to estimate approximation terms I will not commute $(x_2-y_2)$ with $G^\pm, {\bar G}^\pm$ but simply remember that $|x_2-y_2|\le C\rho T_0+C{\bar \gamma}$ on the partition element in question in propagation as $T\le T_0$; this is due to proposition \ref{2-24} extended according to proposition \ref{prop-2-33}. 

\emph{In this subsection I consider  the case when condition  $(\ref{2-103})$ is fulfilled leaving the most general case for the next one\/}. 

Then $\rho=\zeta=\ell$ and $T_0=C\rho^{-2}h|\log h|$ and hence
$|x_2-y_2|\le C\rho T_0+{\bar\gamma}$. Therefore perturbation $R_1=(x_2-y_2)\partial_{x_2}A+O\bigl(|x_2-y_2|^2\bigr)$ does not exceed 
\begin{equation*}
C\zeta (\rho T_0 +{\bar\gamma})+  C(\rho T_0 +{\bar\gamma})^2\asymp r \Def \bigl(h|\log h|+\rho {\bar\gamma})+ \rho^{-2}h^2|\log h|^2+{\bar\gamma}^2\bigr).
\end{equation*}
Therefore while the first term of what is  obtained when one plugs approximations into (\ref{3-61}) does not exceed 
$Ch^{-1}\rho\ell|\log h| =Ch^{-1}\rho^2|\log h|$, each next term acquires a factor \begin{equation}
r\times T/h \asymp  \bigl(h|\log h|+\rho {\bar\gamma}+ \rho^{-2}h^2|\log h|^2+{\bar\gamma}^2\bigr)\times \rho^{-2}|\log h|;
\label{3-62}
\end{equation}
to keep this factor less than 1 one needs to consider
\begin{equation}
\rho \ge \rho^*_0\Def Ch^{\frac 1 2}|\log h|
\label{3-63}
\end{equation}
which is greater than $C{\bar\gamma}$ because under condition (\ref{2-103}) and  ${\bar\gamma}\ge \rho^*_0$ remainder estimate $O\bigl(h^{-1}{\bar\gamma}\bigr)$
is already proven without any correction terms. Under condition (\ref{3-63}) expression (\ref{3-62}) is of magnitude of $\rho^{-2}h|\log h|^2$.

Now the second term of approximations does not exceed 
\begin{equation*}
Ch^{-1}\rho^2|\log h|\times \rho^{-2}h|\log h|^2 \asymp |\log h|^3.
\end{equation*}
Finally, summation with respect to all partition results in $C|\log h|^K$ which does not exceed $Ch^{-1}{\bar\gamma}$ under condition (\ref{2-111}).

On the other hand, contribution of \emph{core} $\{\rho \le \rho^*_0\}$ to the remainder estimate due to Tauberian arguments with $T=\epsilon{\bar\gamma}$ does not exceed $Ch^{-1}\rho^{*\,2}_0\le C|\log h|^K$.

\smallskip

So, one needs to consider only the first term of approximation.
Therefore, under conditions $(\ref{2-103})$ and $\mu \le h^{-\nu}|\log h|^{-K}$ remainder estimate $Ch^{-1}{\bar\gamma}$ still holds with $u$ replaced in the principal part (\ref{3-30}) by the first term of this successive approximation procedure. 

Thus we can take $T_0=\infty$ thus arriving to 
\begin{multline}
\cE^\MW_{\corr,Q} = 
\int \Bigl( (2\pi h)^{-1} {\bf e}(x_1,x_1,0;x_2,\xi_2,\hbar,\mu) -\\ 
h^{-1}\int_{-\infty}^0 F_{t\to h^{-1}\tau}{\bar\chi}_{\bar T}(t)
\Gamma _x\bigl ({\bar U}\bigr)\,d\tau\Bigr)\psi(x_2){\bar\psi}_1(x) \varphi(\xi_2)\,dx_1dx_2d\xi_2
\label{3-64}
\end{multline}
where ${\bf e}={\bf e}(x_1,y_1,\tau ; x_2,\xi_2,\hbar,\mu)$ is a Schwartz kernel of the spectral projector, associated with one-dimensional Schr\"odinger operator 
${\bf a}={\bf a}(x_2,\xi_2,\hbar,\mu)$ which is obtained from $A$ by replacement $hD_2\mapsto \xi_2$ and subsequent change of variables $x_1\mapsto x_1/{\bar\gamma}$ but not setting $\sigma=\phi=1$, $V=W$ and $U=U(x_1,y_1,t;x_2,\xi_2,\hbar,\mu)$ is a corresponding propagator, i.e. Schwartz kernel of $e^{i\hbar^{-1}t{\bf a}}$.
Here and below $Q=\psi(x_2)\varphi(\xi_2)$ denotes partition element.

However after summation over $x_1$-partition remembering that as $|x_1|\ge C$
on $\supp {\bar\psi}_1$ we get negligible terms, we can replace ${\bar\psi}_1$ by $1$ thus arriving to
\begin{multline}
\cE^\MW_{\corr,Q} = 
\int \Bigl( (2\pi h)^{-1}{\bf n}(0;x_2,\xi_2,\hbar,\mu) -\\ 
h^{-1}\int_{-\infty}^0 F_{t\to h^{-1}\tau}{\bar\chi}_{\bar T}(t)
 \Gamma '\bigl ({\bar U}\bigr)\,d\tau \Bigr)\psi(x_2) 
\varphi(\xi_2)\, dx_2d\xi_2
\label{3-65}
\end{multline}
with ${\bf n}$   eigenvalue counting function associated with the same operator ${\bf a}$; I remind that $\Gamma'$ includes integration with respect to $x_1$.

Without the last term formula (\ref{3-64}) delivers the correct principal part of asymptotics which after summation over all partitions (including non-periodic and outer zones) is
\begin{equation}
(2\pi h)^{-1}\int    {\bf e}(x_1,x_1,0;x_2,\xi_2,\hbar,\mu) \psi(x_2) \psi_1(x_1)\,dx_1dx_2d\xi_2
\label{3-66}
\end{equation}
which as we show in the next section is a correct answer without condition (\ref{2-103}) as well.

\smallskip
However, I would like to replace the reference to operator ${\bf a}$ by the reference to operator ${\bf a}_0$. To do this I need the following

\begin{proposition}\label{prop-3-18}  Let condition $(\ref{2-103})$ be fulfilled. Let $\rho \ge C\rho^*_0$. Then

\smallskip
\noindent
(i) Eigenvalues $\lambda_n$ of ${\bf a}$  satisfy inequality
\begin{equation}
|\nabla _{x_2,\xi_2}\lambda_n|\asymp \rho
\label{3-67}
\end{equation}
as  $|\lambda_n|\le \epsilon \rho^2$;

\smallskip
\noindent
(ii) As $\rho \ge C{\bar\gamma}$ 
\begin{align}
&|\lambda_n(x_2,\xi_2) -\lambda^0_n(x_2,\xi_2)|\le C {\bar\gamma}^2\qquad 
&&\nu \text{ is even},\label{3-68}\\
&|\lambda_n(x_2,\xi_2) -\lambda^0_n(x_2,\xi_2)-\beta (x_2,\xi_2)\xi_2 |\le C {\bar\gamma}^2, \quad |\beta|\le C{\bar\gamma}\qquad 
&&\nu \text{ is odd},\label{3-69}
\end{align}
where $\lambda^0_n$ are corresponding eigenvalues of ${\bf a}_0$.
\end{proposition}

\begin{proof} Note that ${\bf a}={\bf a}_0 + t{\bf a}'$ with ${\bf a}'$ bounded 
by $C{\bar\gamma}\bigl({\bf a}+C_0\bigr)$ and $t=1$ therefore due to Rellich decomposition one needs to prove that 
$|\partial_t \lambda^t_n|\le C{\bar\gamma}^2$ as $t=0$. As $\nu$ is even eigenfunctions $\Upsilon_n$ of ${\bf a}_0$ are either even or odd and therefore $\langle {\bf a}'_{{\rm odd}} \Upsilon _n,\Upsilon_n\rangle =0$ while
$|\langle {\bf a}'_{{\rm even}} \Upsilon _n,\Upsilon_n\rangle |\le C{\bar\gamma}^2$.

Consider odd $\nu$. Then $k^*_\hbar=0$ just due to symmetry and the above arguments are applicable for $k=0$. On the other hand 
$\lambda_n \equiv \lambda_{n,\text{sc}}$, 
$\lambda^0_n \equiv \lambda^0_{n,\text{sc}}$ $\mod O\bigl(h^2\bigr)$ with semiclassical expressions $ \lambda_{n,\text{sc}}$, $\lambda^0_{n,\text{sc}}$ and we need to prove (\ref{3-69}) for these expressions; due to analyticity etc these equalities are true as ${\bar\gamma}\ge Ch|\log h|$ i.e. in our assumptions. Since $\lambda^t_{n,\text{sc}}$ is analytic with respect to $\hbar, t,\xi_2$, coincides with $\lambda^0_{n,\text{sc}}$ modulo $O({\bar\gamma})$ 
and also coincides with $\lambda^0_{n,\text{sc}}$ modulo $O({\bar\gamma}^2)$ as $\xi_2=0$, we get (\ref{3-69}).
\end{proof}

I want to remind that also 
\begin{equation}
\cE^\MW_{\corr, Q}=
h^{-1}\int \Bigl(\int_{-\infty}^0 F_{t\to h^{-1}\tau}
\bigl({\bar\chi}_T (t)-{\bar\chi}_{\bar T}(t)\bigr) \Gamma ' ({\bar U})\,d\tau\Bigr)
\varphi(\xi_2)\psi(x_2)\,dx_2d\xi_2
\label{3-70}
\end{equation}
and one can rewrite it as the sum of
\begin{equation}
T^{-1}\int \Bigl(\int_{-\infty}^0 F_{t\to h^{-1}\tau}
{\check\chi}(t) \Gamma ' ({\bar U})\,d\tau\Bigr)
\varphi(\xi_2)\psi(x_2)\,dx_2d\xi_2.
\label{3-71}
\end{equation}

\begin{proposition}\label{prop-3-19} In frames of proposition \ref{prop-3-18}
as $h^{-\nu/2}\le \mu \le h^{-\nu}|\log h|^{-K}$

\smallskip
\noindent
(i) With an error not exceeding $C\rho {\bar\gamma} h^{-1}$ one can replace in the right-hand expressions of $(\ref{3-70})$ ${\bar U}$ by ${\bar U}_0$ associated with operator ${\bf a}_0(x_2,\xi_2,\hbar)$; therefore in  $(\ref{3-65})$ one can replace ${\bf n}$  by ${\bf n}_0$ and simultaneously ${\bar U}$ by ${\bar U}_0$;

\smallskip
\noindent
(ii) $|\cE^\MW_{\corr,I} |\le C{\bar\gamma}^{-1}$ and this estimate cannot be improved unless $V(0,x_2)$ has no critical points. In particular, as $\mu\le C h^{-\nu/2}$ and only then one can skip $\cE^\MW_\corr$ without deteriorating remainder estimate $C{\bar\gamma}h^{-1}$.
\end{proposition}

 \begin{proof}  I remind that under condition (\ref{2-103})  the contribution of the partition element to the correction term does not exceed $Ch^{-1}\rho^2|\log h|^K$. Then running successive approximation method with unperturbed operator  ${\bf a}_0$ and using proposition \ref{prop-3-18} one can see easily that each next term gains factor 
 $C\rho  {\bar \gamma}\times T_0 /h= C\rho ^{-1}{\bar \gamma}|\log h|$
 and thus the error does not exceed
\begin{equation*}
 Ch^{-1}\rho^2|\log h|^K \times \rho ^{-1}{\bar \gamma}|\log h|= Ch^{-1} \rho {\bar\gamma}|\log h|^K
\end{equation*}
and summation over partition results in 
$Ch^{-1}{\bar\gamma} {\bar\rho}_1|\log h|^{K+2}\ll h^{-1}{\bar\gamma}$. 

Then (ii) follows from (i) and the calculations similar to those of the previous subsection. I leave easy details to the reader.
\end{proof}

Then we arrive to

\begin{proposition}\label{prop-3-20} Under conditions $(\ref{2-103})$ as
$h^{-\nu/2}\le \mu \le h^{-\nu}|\log h|^{-K}$ 

\smallskip
\noindent
(i) Asymptotics 
\begin{equation}
\int e(x,x,0)\psi (x)\,dx \equiv \int \cE_\corr^\MW (x)\psi (x)\,dx +
\int \cE^\MW (x_2)\psi(x_2)\,dx_2 \quad \mod O\bigl(h^{-1}{\bar\gamma}\bigr)
\label{3-72}
\end{equation}
holds with $\cE^\MW_\corr$ defined by $(\ref{3-42})$;

\smallskip
\noindent
(ii) Moreover, with 
$O\bigl(h^{-1}\hbar^{3/2}+{\bar\gamma}_1^{-1} \bigr)$-error one can replace $\cE^\MW_\corr$ by expression $(\ref{3-57})$ and with $O\bigl({\bar\gamma}_1^{-1} \bigr)$-error one can replace $\cE^\MW_\corr$ by  expression \ref{3-52-*} modified in the same way.
\end{proposition}

\subsection{Periodic orbits. III. General settings (continuation)}

Let us consider more general settings. Then,  as I have shown in the previous subsection, for any regular partition element described in subsection 2.9 its contribution to the correction term does not exceed $Ch^{-1}\ell\rho |\log h|^K$ (where actually $K$ should be equal to 1).

Let us consider a successive approximation method with unperturbed operator obtained by fixing $x_2=y_2$ in the coefficients. Again as before $|x_2-y_2|$
does not exceed $C(\rho T_0 + {\bar\gamma})$ and  perturbation is $(x_2-y_2)B$
with $B=(\partial_{x_2}A|_{x_2=z_2}$ calculated at some point $z_2$ of the same partition element (containing both $y_2$ and $x_2$) and thus can be estimated by $C\zeta$; thus the perturbation is estimated by 
$C\zeta (\rho T_0 + {\bar\gamma})$ and the contribution of the second term to the correction term is estimated by

\begin{multline}
Ch^{-1}\rho\ell |\log h|^K\times \zeta\bigl(\rho  T_0 +{\bar\gamma}\bigr)\times T_0/h \le \\
C\ell \zeta \Bigl(\rho^2  (\rho^2+\zeta\ell)^{-2} + h^{-1}{\bar\gamma}\rho (\rho^2+\zeta\ell)^{-1}\Bigr) |\log h|^{K+2} \le \\
C\Bigl(1+ h^{-1}{\bar\gamma}\rho\Bigr)|\log h|^{K+2}\le 
Ch^{-{\frac 1 2}}{\bar\gamma}^{\frac 1 2}|\log h|^{K+3}
\label{3-73}
\end{multline}
as $T_0=Ch(\rho^2+\zeta\ell)^{-1}|\log h|$. 

This was a contribution of the regular $(\ell,\rho)$-element. This element was actually a subelement of some more rough element and then the contribution of this larger element would be a sum over finer subpartition. So in frames of subsection 2.9 as long as we use nondegeneracy condition \ref{2-105} contribution of this larger element would be 
$Ch^{-{\frac 1 2}}{\bar\gamma}^{\frac 1 2}|\log h|^{K+4}$ as well. Continue this process  until the very top we get 
$Ch^{-{\frac 1 2}}{\bar\gamma}^{\frac 1 2}|\log h|^{2K}$ with 
$\ell_m\asymp 1$ thus resulting in $C|\log h|^K$. So, under nondegeneracy condition \ref{2-105} we arrive to a proper error estimate $O\bigl(h^{-1}{\bar\gamma}\bigr)$ as 
$\mu \le h^{-\nu}|\log h|^{-N\nu}$.

On the other hand, without nondegeneracy condition \ref{2-105} we still have this condition fulfilled until the very top when we join regular elements to $B(0,1)$; then we arrive to the error estimate $O\bigl(h^{-1}{\bar\gamma}+h^{-\delta}\bigr)$ with arbitrarily small exponent $\delta>0$.

On the other hand contribution of the irregular elements to the correction terms are less than $C|\log h|^K$, $Ch^{-\delta}$ with/without non-degeneracy condition \ref{2-105} respectively.

Thus we arrive to

\begin{proposition}\label{prop-3-21} Let $\mu \le Ch^{-\nu}|\log h|^{-N}$. Then

\smallskip
\noindent
(i) Under nondegeneracy condition \ref{2-105} asymptotics holds with the principal part defined by $(\ref{3-66})$ and the remainder estimate $O\bigl(h^{-1}{\bar\gamma}\bigr)$;

\smallskip
\noindent
(ii) Without nondegeneracy condition \ref{2-105} asymptotics holds with the principal part defined by $(\ref{3-66})$ and the remainder estimate $O\bigl(h^{-1}{\bar\gamma}+h^{-\delta}\bigr)$.
\end{proposition}

Our next step is

\begin{proposition}\label{prop-3-22}\footnote{\label{foot-19} Cf proposition \ref{prop-3-18}}  Let $\rho \ge C{\bar\gamma}$. Then
eigenvalues $\lambda_n$ of ${\bf a}$  satisfy inequalities
\begin{equation}
\rho^{-1}|\nabla _{\xi_2}\lambda_n|+\zeta^{-1}|\nabla _{x_2}\lambda_n| \le C
\label{3-74}
\end{equation}
and $(\ref{3-68}$ or $(\ref{3-69})$ as $|\lambda_n|\le \epsilon \rho^2$.
\end{proposition}

\begin{proof}
Proof of (\ref{3-74}) is obvious and (\ref{3-68}),(\ref{3-69}) are already proven.
\end{proof}

\begin{proposition}\label{prop-3-23} Let $\mu \le Ch^{-\nu}|\log h|^{-N}$. Then

\smallskip
\noindent
(i) Under nondegeneracy condition \ref{2-105} asymptotics $(\ref{3-72}$ holds with the correction term defined by $(\ref{3-42})$ and the remainder estimate $O\bigl(h^{-1}{\bar\gamma}\bigr)$;

\smallskip
\noindent
(ii) Without  nondegeneracy condition \ref{2-105} asymptotics $(\ref{3-72}$ holds with the correction term defined by $(\ref{3-42})$ and the remainder estimate $O\bigl(h^{-1}{\bar\gamma}+h^{-\delta}\bigr)$.

\smallskip
\noindent
(iii)  Moreover, with 
$O\bigl(h^{-1}\hbar +{\bar\gamma}_1^{-1} \bigr)$-error one can replace $\cE^\MW_\corr$ by expression $(\ref{3-52})$ and with $O\bigl(h^{-1}\hbar^{3/2}+{\bar\gamma}_1^{-1} \bigr)$-, $O\bigl({\bar\gamma}_1^{-1} \bigr)$-error  one can replace $\cE^\MW_\corr$ by  expressions \ref{3-52-*}, \ref{3-52-**} respectively modified in the same way
\footnote{\label{foot-20} Surely under condition \ref{2-104} better estimates  (depending on $m$) holds.}.
\end{proposition}

\begin{proof} I remind that contribution of the regular $(\ell,\rho)$ partition element to the correction term is $O\bigl( h^{-1}\rho\ell\bigr)$.
For even $\nu$ application of the successive approximation method brings in every step factor 
$C{\bar\gamma}^2\times T_0/h =C{\bar\gamma}^2\rho^{-2}|\log h|$ and the contribution of the second term is 
$O\bigl(h^{-1}{\bar\gamma}^2\rho ^{-1}\ell |\log h| \bigr)$. Then summation over partition with $\rho \ge {\bar\rho}_0=C{\bar\gamma}|\log h|^N$ results in $O\bigl(h^{-1}{\bar\gamma}^2 {\bar\rho}_0^{-1}|\log h|^K\bigr)= O\bigl(h^{-1}{\bar\gamma}\bigr)$. 

For odd $\nu$ application of the successive approximation method brings in every step factor $C{\bar\gamma}\rho^{-1}|\log h|$ and all above estimates remain true for the third term in successive approximations and even for the the second term with the single exception of its part which is of magnitude 
$O\bigl(h^{-1} {\bar\gamma}\ell |\log h|\bigr)$ and which emerges from the linear with respect to $\xi_2$ term in formula (\ref{3-69}). One can see easily that this exceptional term becomes 0 after integration with respect to $\xi_2$.

Now in the both cases one needs to look at contribution of zone $\{|\xi_2-k^*_\hbar|\le {\bar\rho}_0\}$ which is $O(h^{-1}{\bar\rho}_0L)$ where $L$ is the total length of the corresponding $x_2$ intervals and thus is $O(h^{-1}{\bar\gamma})$ under condition \ref{2-105}; otherwise it has a superficial logarithmic factor we need to get rid of.

When looking at such narrow zone due to above arguments one can skip all the $O(\rho^3)$ terms in $\lambda_n-\lambda^0_n$ and thus we can assume with no loss of the generality that both $\lambda_n$ and $\lambda^0_n$ are quadratic forms with respect to $\eta=\xi_2-k^*_\hbar$ and ${\bar\gamma}$ with coinciding coefficient at $\eta^2$ and after obvious transformations we can assume\footnote{\label{foot-21}This assumption is valid only if we are interested in the difference of correction terms.} that 
$\lambda^0_n= \varkappa \eta^2 +\varkappa_0$, 
$\lambda_n= \varkappa \eta^2 +\varkappa_1$ with $\varkappa_1-\varkappa_0=O({\bar\gamma}^2)$.  Then we can take integration with respect to  $\xi_2$ all over $\bR$. But then both correction terms would be of the form $c_1 h^{-1}\hbar ^{1/2} G (c_2+ c_3 \varkappa _j/\hbar)$ with $j=1,2$ and equal coefficients $c_k$, $k=1,2,3$; since $G\in C^{1/2}$  the difference does not exceed $Ch^{-1}\hbar ^{1/2}\times ({\bar\gamma}^2/\hbar)^{1/2}= Ch^{-1}{\bar\gamma}$. 
\end{proof}

Thus all main theorems are proven as $\mu \le h^{-\nu}|\log h|^{-K}$.

\sect{Superstrong Magnetic Field}

In this section I assume  that $\mu$ is close to $h^{-\nu}$; more precisely I assume first that 
\begin{equation}
 h^{\delta -\nu}\le \mu \le \epsilon h^{-\nu}
 \label{4-1}
\end{equation}
with very small exponent $\delta>0$ constant $\epsilon$ and I leave the case $\epsilon h^{-\nu}\le \mu\le C_0h^{-\nu}$ for subsection 4.5. Under assumption (\ref{4-1})  parameter $\hbar=h/{\bar\gamma}$ is small but not very small: 
$h^\delta \le \hbar\ll 1$. Moreover, as 
$\hbar \ge \epsilon |\log h|^{-1}$ variable $x_1$ is no more microlocal.

I am going to consider operator in frames of pdo theory with operator valued symbols in the axillary space $\bH=L^2(\bR)$ with inner product $\blangle.,.\brangle$ and norm $\bv .\bv$. 

\subsection{Outer zone}

Propagation estimates in the outer zone are basically done: as 
$\mu \le \epsilon h^{-\nu}|\log h|^{-\nu}$ I  already proved that the propagation trace $F_{t\to h^{-1}\tau}\chi_T(t)\Gamma (\psi Qu)$ is negligible  as $T\in [T_0,T_1]$ with $T_0= Ch|\log h|$ and  
$T_1=\epsilon' \rho $, where $\psi=\psi(x_2)$ and $Q=Q(hD_2)$ with $\rho\le |\xi_2|\le 2\rho$ on the support of $Q$ and 
$C\le \rho \le C\bigl(\mu h^\nu\bigr)^{-1/(\nu -1)}$. Further, as 
$\rho \ge C\bigl(\mu h^\nu\bigr)^{-1/(\nu -1)}$ the standard ellipticity arguments imply that the propagation trace is negligible even as $\chi_T(t)$ is replaced by ${\bar\chi}_T(t)$.

Also I proved that as $\epsilon h^{-\nu}|\log h|^{-\nu}\le \mu \le C h^{-\nu}$
propagation trace  $F_{t\to h^{-1}\tau}\chi_T(t)\Gamma (\psi Qu)$ is negligible  as $T\in [T_0,T_1]$ with the same $T_0,T_1$ and $C\le \rho \le C|\log h|^\nu$ 
while ellipticity arguments work as $\rho \ge C|\log h|^\nu$. 

But then one can prove estimate 
\begin{equation}
|F_{t\to h^{-1}\tau}{\bar\chi}_T(t)\Gamma (\psi Qu)|\le 
Ch^{-1}{\bar\gamma }\rho^{1/\nu}|\log h|
\label{4-2}
\end{equation}
which is ``almost perfect'':  the problem here lies only with the logarithmic factor.

To improve  estimate (\ref{4-2}) let us  launch launch successive approximation method, fixing in the coefficients of unperturbed operator $x_2=y_2$ as usual and also fixing $x_1=0$ there but not in $\mu x_1^\nu/\nu$. Then one can estimate easily the final contribution of the first term  by the same
$Ch^{-1}{\bar\gamma }\rho^{1/\nu}|\log h|$ as before and also estimate by $C$ contribution of all other terms. Therefore \emph{in the estimate part  only the first term should be counted}.  But then one need to consider
\begin{equation}
h^{-1}\int \Bigl(F_{t\to h^{-1}\tau} \bigl(\Gamma' U(t;\xi_2)\bigr)\Bigr)Q(\xi_2)\,dx_2d\xi_2
\label{4-3}
\end{equation}
where 
\begin{equation}
U=U(x_1,y_1,t;x_2,\xi_2)=\sum_n e^{ih^{-1}\lambda_n(x_2,\xi_2,\hbar)t}\Upsilon_n(x_1;x_2,\xi_2,\hbar)
\Upsilon_n(y_1;x_2,\xi_2,\hbar)
\label{4-4}
\end{equation}
is a propagator for 1-dimensional operator which after change of variable $x_1\mapsto x_1/{\bar\gamma}$ is just ${\bf a}_0={\bf a}_0(x_2,\xi_2,\hbar)$. In the below arguments let us skip integration over $x_2$ and $x_2$ in notations . So, $\lambda_n(\xi_2,\hbar)$ are just eigenvalues of 
${\bf a}_0(\xi_2,\hbar)$ and $\Upsilon_n (x_1;\xi_2,\hbar)$ are orthonormalized eigenfunctions. Plugging (\ref{4-4}) into (\ref{4-3}) I arrive to 
\begin{equation}
 \sum_n h^{-1}T\int
{\hat{\bar\chi}}\Bigl({\frac T h}\bigl(\tau -\lambda_n(\xi_2,\hbar)\bigr)\Bigr)
Q(\xi_2)\,d\xi_2
\label{4-5}
\end{equation}
and in the case of even $\nu$ I will consider separately sums with respect to eigenvalues corresponding to eigenfunctions odd and even with respect to $x_1$.

\begin{proposition}\label{prop-4-1} As $\hbar\ll 1$ eigenvalues $\lambda_n(\xi_2,\hbar)$ of 
${\bf a}_0(\xi_2,\hbar)$ have the following properties as $|\xi_2|\ge C_1$:

\smallskip
\noindent
(i) $\lambda_n(\xi_2,\hbar) \ge \Delta \Def 
\epsilon_0\hbar |\xi_2|^{(\nu -1)/\nu}-C_0$;

\smallskip
\noindent
(ii) For even $\nu$ and $\xi_2\le -\epsilon$ all eigenvalues are larger than $C_0\epsilon^2$;

\smallskip
\noindent
(iii) Spacing between two consecutive eigenvalues\footnote{\label{foot-22} As $\nu$ is even I consider two series of eigenvalues separately: those with even with respect to $x_1$ eigenfunctions and those with odd eigenfunctions.} $\lambda_n(\xi_2,\hbar)$ and $\lambda_{n+1}(\xi_2,\hbar)$ belonging to interval $(-\epsilon_0,\epsilon_0)$ is not less than $\Delta$; 

\smallskip
\noindent
(iv) Also for these eigenvalues 
\begin{equation}
C_0\ge \xi_2\partial_{\xi_2}\lambda_n(\xi_2,\hbar)\ge \epsilon_0,
\label{4-6}
\end{equation}
\end{proposition}

\begin{proof} Proof of (i)-(iii)  follows from Bohr-Sommerfeld theory with the semiclassical parameter $h'=\hbar\rho^{(\nu-1)/\nu}$ when it is small enough. 
Note that zone $\{|x_1|\le \epsilon {\bar\gamma}\rho^{1/\nu}\}$ is classically forbidden with the smaller parameter $h''=\hbar/\rho$.

 The easiest way to prove (\ref{4-6}) is to note that propagation results of subsection 2.2 basically are equivalent to these inequalities as 
 $\mu \le \epsilon h^{-\nu}|\log h|^{-\nu}$; however this condition could be replaced just by $\hbar\le \epsilon'$ because $e_0(x_1,y_1,.,.)$ depends  on $\hbar$, $\xi_2$ rather than $\mu, h$ separately.
\end{proof}

Let us consider (\ref{4-5}). Since ${\hat{\bar\chi}}$ is fast-decaying at infinity, each term in (\ref{4-5}) does not exceed
\begin{equation*}
h^{-1}T\int \bigl(1+ T|\tau-\lambda_n(\xi_2)|h^{-1}\bigr)^{-s}\,d\xi_2 \asymp
\rho
\end{equation*}
due to (\ref{4-6}); on the other hand the number of terms in (\ref{4-5}) does not exceed $C\Delta^{-1}\asymp \hbar^{-1}\rho^{(1-\nu)/\nu}$. Therefore 

\begin{proposition}\label{prop-4-2}As $\hbar\ll 1$ and $\rho\ge C$ expression $(\ref{4-5})$ does not exceed $C\hbar^{-1}\rho^{1/\nu}$.
\end{proposition}

So, I got rid of the logarithmic factor in (\ref{4-2}). Then dividing by $T_1=\epsilon \rho$  I conclude that the contribution of this partition element $Q$ to the remainder estimate does not exceed $C\hbar^{-1}\rho^{(1-\nu)/\nu}$. After integration over $d\rho/\rho$ I conclude that  the contribution of the whole outer zone to the remainder estimate does not exceed $C\hbar^{-1}=Ch^{-1}{\bar\gamma}$. So I arrive to 

\begin{proposition}\label{prop-4-3} Under condition $(\ref{4-1})$ contribution of the outer zone $\cZ_\out=\{|\xi_2|\ge C\}$ to the remainder estimate does not exceed $Ch^{-1}{\bar\gamma}$.
\end{proposition}

\subsection{Inner zone. I. Estimates}

So far I have no results in this zone 
$\cZ_\inn=\{|\xi_2|\le C_0\}$ as $\mu \ge \epsilon h^{-\nu}|\log h|^{-\nu}$. Thus  I cannot even refer to operator ${\bf a}_0$ directly. However let us start from the properties of ${\bf a}_0={\bf a}_0(\xi_2,\hbar)$ with potential $V=1$. Let us first to move through $\xi_2=1$ or $\xi_2=\pm 1$ for even/odd $\nu$ respectively.

\begin{proposition}\label{prop-4-4} Let $\hbar \le \epsilon$ and 
$1-\epsilon\le |\xi_2 |\le C_0$ with small enough constant $\epsilon >0$. Then 
operator ${\bf a}_0(\xi_2,\hbar)$ is microhyperbolic with respect to $\xi_2$ on energy level $0$; more precisely
\begin{equation}
\xi_2 \blangle \partial_{\xi_2}{\bf a}_0(\xi_2,\hbar)v,v\rangle \ge 
\epsilon \bv v\bv^2 - C\hbar^{-r}\bv{\bf a}_0(\xi_2,\hbar)v\bv^2
\qquad \qquad \forall v\in D({\bf a}_0)
\label{4-7}
\end{equation}
with large enough exponent $r$.
\end{proposition}

\begin{proof} Let us consider $\xi_2\ge 1-\epsilon$ in the both cases; case 
$\xi_2\le -1+\epsilon$ and odd $\nu$ is treated due to the symmetry. It is sufficient to prove (\ref{4-7}) for $v\in 
\bigl({\bf E}(\xi_2,\hbar; \epsilon \hbar^r)-
{\bf E}(\xi_2,\hbar; -\epsilon\hbar^r)\bigr)\bH$ where ${\bf E}(\xi_2,\hbar;\lambda,\hbar)$ is the spectral projector of 
${\bf a}_0(\xi_2,\hbar)$. Then the standard WKB approach implies that for any $v$ described above with 
$\bv v\bv =1$ and $\varepsilon >0$ one can approximate $v$ by WKB solution as $x_1>\varepsilon$:
$\bv v - w\bv_{\{x_1\ge \epsilon_1\}} \le C\hbar ^{r-2}$ where $w$ is a WKB solution constructed for $\lambda=0$, $C=C(\epsilon_1)$, and 
$\bv w\bv_{\{x_1\ge \epsilon_1\}} \le 2$. Then 
\begin{equation*}
|w|^2 = \varkappa \Bigl( {\frac 1 {\sqrt{1-(\xi_2-x_1^\nu/\nu)^2}}}-o(1)\Bigr)
\end{equation*}
with $0\le \varkappa \le c$.

Then, from explicit the WKB calculations it follows that 
\begin{equation*}
\blangle \partial_{\xi_2}{\bf a}_0(\xi_2,\hbar)w,w\rangle_{\{x_1\ge \varepsilon\}} = \const \Bigl(\int_\varepsilon^{b(\xi_2)}  {\frac {(\xi_2-x_1^\nu/\nu)\, dx_1}{\sqrt{1-(\xi_2-x_1^\nu/\nu)^2}}}-o(1)\Bigr)\ge 
3\epsilon \bv w\bv_{\{x_1\ge \varepsilon\}} ^2\end{equation*}
{and therefore}
\begin{equation}
\blangle \partial_{\xi_2}{\bf a}_0(\xi_2,\hbar)v,v\rangle_{\{x_1\ge 
\varepsilon \}} \ge 
2\epsilon \bv v\bv_{\{x_1\ge \varepsilon \}} ^2 -C\hbar^{r-2}.
\label{4-8}
\end{equation}
The last inequality holds for domain $\{x_1\le -\varepsilon \}$ as well because it is either classically forbidden (for odd $\nu$) or just due to symmetry (for even $\nu$). Furthermore (\ref{4-8}) obviously holds for domain 
$\{|x_1|\le \varepsilon \}$. After summation I arrive to (\ref{4-8}) with integration over $\bR\ni x_1$.

Easy details linked to WKB calculations I leave to the reader. 
\end{proof}

\begin{proposition}\label{prop-4-5} $\hbar \le \hbar_0(\epsilon_0)$.  

\medskip
\noindent
(a) Further, let  $\nu $ be even. Then

\smallskip
\noindent
(i) As $1+\epsilon \le  \xi_2\le C_0$ spacing between eigenvalues\footnote{\label{foot-23} In this proposition I consider only eigenvalues belonging to $[-\epsilon, \epsilon]$.} $\lambda_n(\xi_2,\hbar)$ and $\lambda_{n+1}(\xi_2,\hbar)$ (since $\nu$ is even I separate two series of eigenvalues; see footnote $^{\ref{foot-22}}$) is at least $\epsilon \hbar $ where here and below $\epsilon$ is a small positive constant;

\smallskip
\noindent
(ii) As $\epsilon_0\le \xi_2 \le 1-\epsilon_0$ spacing between eigenvalues $\lambda_n(\xi_2,\hbar)$ and $\lambda_{n+1}(\xi_2,\hbar)$ (I do not separate two series of eigenvalues anymore) is at least $\epsilon \hbar$;

\smallskip
\noindent
(iii) As $| \xi_2 |\le \epsilon_0$ spacing between eigenvalues $\lambda_n(\xi_2,\hbar)$ and $\lambda_{n+1}(\xi_2,\hbar)$  is at least 
$\epsilon \hbar^{2\nu/(\nu+1)}$;

\medskip
\noindent
(b) On the other hand, let $\nu$ be odd. Then

\smallskip
\noindent
(iv) As $1+\epsilon \le  |\xi_2|\le C_0$ spacing between eigenvalues $\lambda_n(\xi_2,\hbar)$ and $\lambda_{n+1}(\xi_2,\hbar)$  is at least 
$\epsilon \hbar $;

\smallskip
\noindent
(v) As $|\xi_2 |\le 1-\epsilon_0$ spacing between eigenvalues $\lambda_n(\xi_2,\hbar)$ and $\lambda_{n+1}(\xi_2,\hbar)$  is at least 
$\epsilon \hbar$.
\end{proposition}

\begin{proof} All cases but (iii) follow easily from the standard WKB method. Moreover, (iii) follows from the standard WKB method and rescaling.
\end{proof}

So spacing between eigenvalues is rather large - larger than  $\epsilon \hbar^2$ (as $\hbar^2\ge h^{-\delta}$) and therefore their derivatives of order $\alpha$ could be estimated by $C_\alpha\hbar^{-2\alpha}$. Further, $|\partial_{\xi_2,\hbar}\lambda_n(\xi_2)|$ is bounded.

Then decomposing $u(x,y,t)$ into $\Upsilon_n(x_1/{\bar\gamma};hD_2,\hbar)$\,\footnote{\label{foot-24} Where $\Upsilon_n(x_1,\xi_2,\hbar)$ are eigenfunctions of ${\bf a}_0(\xi_2,\hbar)$.} which are functions of $x_1$ but also $h$-pdo theory with respect to $x_2$ 
\begin{equation}
u(x,y,t)\equiv \sum_{n,l}{\bar\gamma}^{-1}
\Upsilon_n\bigl({\frac {x_1}{\bar\gamma}};hD_{x_2},\hbar\bigr) u_{n l}(x_2,y_2,t)
\Upsilon_l\bigl({\frac {y_1}{\bar\gamma}};hD_{y_2},\hbar\bigr)^\dag
\label{4-9}
\end{equation}
one can rewrite the basic equation $(hD_t-A)u=0$ as a system
\begin{align}
&\Bigl(hD_t -\Lambda_n(x_2,hD_2,\hbar) \Bigr) u_{nl}\equiv\sum_k \cB_{kl}u_{kl},\label{4-10}\\
&\Lambda_n(x_2,hD_2,\hbar)=\lambda_n(hD_2,\hbar) - {\frac 1 2}W(x_2)+{\frac 1 2}
\label{4-11}
\end{align}
with operators $\cB_{kl}$ estimated by $Ch \hbar^{-K}$ where I increase $K$ if needed.

Let us  consider $\epsilon \hbar^2$-sized intervals with respect to $x_2$ and $\tau$.  On each such interval  at most one of operators $hD_t -\Lambda_n(x_2,hD_2,\hbar)$ fails to be elliptic there (with the ellipticity constant $\epsilon \hbar^2$). There is one exception: if $\nu$ is even and $\xi_2\ge 1-\epsilon$ ``odd'' and ``even'' eigenvalues\footnote{\label{foot-25} I.e. eigenvalues corresponding to eigenfunctions even or odd with respect to $x_1$} could be pretty close to one another and then there are at most two such operators with numbers $n$ and $n+1$ (I cannot separate eigenvalues anymore; $n$ corresponds to even eigenfunctions and $(n+1)$ to odd ones).

Let $n$ be such exceptional number.
Then one can rewrite system (\ref{4-10}) either as a single equation
\begin{equation}
\Bigl(hD_t -\Lambda_n(x_2,hD_2,\hbar)  - \cB'_{n}\Bigr) u_{nl}\equiv 0
\label{4-12}
\end{equation}
where $u_{kn}$ with $k\ne n$ are expressed via $u_{nl}$ and 
$\cB'_n \equiv \cB_{nn}$ modulo operators not exceeding $h^2\hbar^{-K}$, or as a $2\times2$-system with $\Lambda_n$ replaced by 
$\begin{pmatrix}\Lambda_n &0\\ 0&\Lambda_{n+1}\end{pmatrix}$ and with $2\times2$ matrix operator $\cB'_n$.

\begin{proposition}\label{prop-4-6} In frames of proposition \ref{prop-4-5} 
\;$\pm \partial_{\xi_2}\lambda_n(\xi_2,\hbar)\ge \epsilon_1$\; as long as 
$\pm (\xi_2-k^*)\ge\epsilon_0$ and 
$|\lambda_n(\xi_2)|\le \epsilon_1$.
\end{proposition}

\begin{proof} Proof immediately follows from Bohr-Sommerfeld approximation for $\lambda_n$.
\end{proof}

So far propositions \ref{prop-4-5}, \ref{prop-4-6} were proven for $W=1$ but they obviously hold for any $W\ge \epsilon $ with the critical values $\pm 1, k^*$ replaced by $\pm W^{1/2}$, $k^*W^{1/2}$ respectively. Then one can easily recover all the results similar to those of subsection 2.3 and thus one can set $T_1=\epsilon$, $T_0=C\hbar^{-r}h|\log h|$ in inner but not periodic zone 
$\cZ_\inn\setminus\cZ_\per=\{\epsilon_0\le |\xi_2-k^*W^{\frac 1 2}|\le C_0\}$. These propagation results imply immediately estimate
\begin{equation*}
|F_{t\to h^{-1}\tau}{\bar\chi}_T(t)\Gamma (Q\psi u)|\le 
C\hbar^{-r}h^{-1-\delta M} {\bar\gamma}|\log h|\le Ch^{-1}{\bar\gamma}.
\end{equation*}
However one can  engage the  successive approximation method as in subsection 4.1 (with $\rho=1$ now); then one needs to consider only the first term of it  and following arguments of subsection 4.1 one can see easily that  
$|F_{t\to h^{-1}\tau}{\bar\chi}_T(t)\Gamma (Q\psi u)|\le Ch^{-1}{\bar\gamma}$. Thus estimate immediately leads to 

\begin{proposition}\label{prop-4-7} For fixed $\epsilon_0$ contribution of zone $\cZ_\inn\setminus\cZ_\per=\{\epsilon_0\le |\xi_2-k^*W^{\frac 1 2}|\le C_0\}$ to the remainder estimate does not exceed $Ch^{-1}{\bar\gamma}$ as $\mu \le \epsilon h^{-\nu}$ with small enough constant $\epsilon=\epsilon (\epsilon_0)$.
\end{proposition}

\subsection{Periodic zone. I. Estimates}

Now I need to treat periodic zone 
$\cZ_\per=\{|\xi_2-k^*W^{\frac 1 2}|\le \epsilon_0\}$.
However now I have just one equation due to proposition \ref{prop-4-5} and the construction after it. 

We need

\begin{proposition}\label{prop-4-8} Let ${\bf a}_0(\xi_2,\hbar)$ be a model operator with the potential $V=1$ and $\hbar \le \epsilon=\epsilon(\epsilon_0)$. Then
\begin{equation}
 \partial^2_{\xi_2}\lambda_n(\xi_2,\hbar)\ge \epsilon_1\qquad\qquad 
 \forall \xi_2: |\xi_2-k^*|\le \epsilon_0\quad 
 \forall n: |\lambda_n(\xi_2,\hbar)|\le \epsilon.
 \label{4-13}
 \end{equation}
\end{proposition}

\begin{proof} Proof follows from Bohr-Sommerfeld approximation and propositions \ref{prop-1-4} and \ref{prop-1-8} for even and odd $\nu$ respectively.
\end{proof}

Now one can apply our standard rescaling technique as in section 3 (but in the simpler form)  with unperturbed operator $\lambda_n(hD_2,\hbar)+{\frac 1 2}-{\frac 1 2}W(x_2)$ leading to

\begin{proposition}\label{prop-4-9} Contribution of zone 
$\cZ_\per=\{ |\xi_2-k^*W^{1/ 2}|\le \epsilon_0\}$ to the remainder estimate does not exceed $Ch^{-1}{\bar\gamma}$ as non-degeneracy condition \ref{2-105} is fulfilled and $Ch^{-1+\delta}{\bar\gamma}$ otherwise with arbitrarily small exponent $\delta>0$.
\end{proposition}

\subsection{Calculations}

Now let us apply our standard successive approximation technique as in section 3  with an unperturbed operator 
${\bf a}_0(hD_2,\hbar)+{\frac 1 2}-{\frac 1 2}W(y_2)$ in the outer and inner non-periodic zones and with an unperturbed operator $\lambda_n(hD_2,\hbar)+{\frac 1 2}-{\frac 1 2}W(y_2)$ in 
the inner periodic zone leading exactly to

\begin{proposition}\label{prop-4-10} Under condition $(\ref{4-1})$ with the remainder estimate  equal to $Ch^{-1}{\bar\gamma}$ as non-degeneracy condition \ref{2-105} is fulfilled and $Ch^{-1+\delta'}{\bar\gamma}$ otherwise with arbitrarily small exponent $\delta'>0$ asymptotics 
\begin{equation}
 \int e(x,x,0)\psi(x_2)\,dx \equiv
 (2\pi h)^{-1}\int {\bf n}_0(x_2,\xi_2,0) \psi(x_2)\, dx_2d\xi_2.
\label{4-14}
\end{equation}
\end{proposition}

Since  $\cE^\MW \equiv \cE^\MW_0$ $\mod O(1)$ under condition (\ref{4-1}) this proposition finalize the proof of Theorem \ref{thm-0-1} under condition (\ref{4-1}).

\subsection{Ultra-strong magnetic field}

Now we are left with the case 
\begin{equation}
\epsilon  h^{-\nu}\le \mu \le C_0h^{-\nu}
\label{4-15}
\end{equation}
when $\hbar$ is not a small parameter anymore and when all the eigenvalues of ${\bf a}_0(\xi_2,\hbar)$ which are less than $C$ are separated by a small constant.

Then one can reduce the analysis to those of 1-dimensional scalar operators $\Lambda_n(x_2,hD_2,\hbar)$ and the problem is that we don't know how $\lambda_n(\xi_2,\hbar)$ depends on $\xi_2$.

So far I can prove  unconditionally only

\begin{proposition}\label{prop-4-11}
Under condition 
\begin{equation}
|\lambda_n(\xi_2,\hbar) +{\frac 1 2} -{\frac 1 2}W|+|\partial_{x_2}W|\ge \epsilon_1
\qquad\forall n,\xi_2,x_2
\label{4-16}
\end{equation}
asymptotics $(\ref{4-14})$ holds modulo $O(1)$.
\end{proposition}

Note that for even $\nu$ perturbations $\cB_{nn}(x_2,\xi_2)$ are $O(h^2)$ since terms of the order $h$ could appear only from multiplication by $x_1$ in the first power but due to the parity of eigenfunctions these terms are 0. For odd $\nu$ we can conclude that $\cB_{nn}(x_2,\xi_2)=O(h^2)$ as $\xi_2=0$ but I want to make conjecture

\begin{equation}
\cB_{nn}(x_2,\xi_2)=O(h^2)\qquad\text{as}\;\partial_{\xi_2}\lambda_n=\partial_{x_2}\lambda_n=0
\label{4-17}
\end{equation}
since I don't know if there are other critical points but $0$. 

Then according to  \cite{Ivr1}, section 4.4 the remainder estimate should be $O(1)$ under conditions
\begin{equation}
\sum_{1\le k\le l} |\partial_{\xi_2}^k\lambda_n| \ge \epsilon_1\qquad \forall \xi_2, n
\label{4-18}
\tag*{$(4.18)_l$}
\end{equation}
with $l=2$ and  nondegeneracy condition $(\ref{2-103})$. Further, note that condition \ref{4-18} is fulfilled with some $l$ due to the analyticity of $\lambda_n(\xi_2,\hbar)$ with respect to $\xi_2$ and the fact that $\lambda_n(\xi_2,\hbar)\to \infty$ as $|\xi_2|\to \infty$. Then according to  \cite{Ivr1}, section 4.4 remainder estimate   $O(h^{-\delta'})$ in the general case is guaranteed but to avoid contribution of junior symbols into the principal part one should make additional assumptions.

\begin{proposition}\label{prop-4-12} (i) Under conditions $(\ref{4-15}),(\ref{4-17})$, $(4.18)_2$  and $(\ref{2-103})$ asymptotics $(\ref{4-14})$ holds with the remainder estimate $O(1)$.

\smallskip
\noindent
(ii) Under conditions $(\ref{4-15}),(\ref{4-17})$ and  $(4.18)_2$  asymptotics $(\ref{4-14})$ holds with the remainder estimate $O(h^{-\delta})$.
\end{proposition}

However analysis of to  \cite{Ivr1}, section 4.4 shows that the the remainder estimate must be $O(h^{-\delta})$ under \ref{4-18} with any $l$ and non-degeneracy condition \ref{2-105} with any $m$ but there would be indefinite correction of magnitude $O(h^{2/l+2/m-1})$. Further, under condition $(2.105)_2$ remainder estimate is $O(1)$. Therefore we arrive to

\begin{proposition}\label{prop-4-13} (i)
Under conditions $(\ref{4-15})$, $(\ref{4-17})$, \ref{4-18}  with $l\ge 3$ and $(2.105)_2$  asymptotics $(\ref{4-14})$ holds modulo $O(1)$;

\smallskip
\noindent
(ii) Under conditions $(\ref{4-15})$, $(\ref{4-17})$, \ref{4-18}  with $l\ge 3$ and \ref{2-105} with $m\ge 3$  asymptotics $(\ref{4-14})$ holds modulo $O\bigl(h^{2/l+2/m-1}\bigr)$.
\end{proposition}

\sect{Generalizations}

\subsection{Vanishing $V$}

We can also generalize all our results to the case of $V$ vanishing:

\begin{theorem}\label{thm-5-1} Theorem \ref{thm-0-1} remains true if condition $(\ref{0-2})$ is replaced by
\begin{equation}
|\partial_{x_2}V|\ge \epsilon_0.
\label{5-1}
\end{equation}
\end{theorem}

\begin{proof} Consider $\ell$-admissible partition with $\ell=\epsilon |V|$.

\smallskip
\noindent
(i) First, consider elements on which $\ell \le \epsilon_1 |x_1|$. Then rescaling  $x\mapsto x/\ell$ and diving by $\ell$ we arrive to the situation of non-vanishing magnetic field $\mu_\eff=\mu \gamma_1^{\nu-1} \ell^{1/2}$ with 
$F_\eff = F\gamma^{-\nu}$ and $V_\eff = V\ell^{-1}$ with a parameter 
$\gamma \asymp |x_1|$. 

Further, after rescaling
$V_\eff/F_\eff = (V/F) \cdot (\gamma^\nu/\ell)\asymp 1$ and also
$\ell |\nabla V_\eff/F_\eff|\asymp 1$.

 To apply ``non-vanishing'' theory with $h_\eff = h\ell^{-3/2}$ and $\mu_\eff$ one must assume that $h_\eff\le 1$, $\mu_\eff\ge 1$ or equivalently
\begin{equation}
\ell \ge C_0\max \bigl(h^{\frac 2 3}, 
(\mu \gamma^{\nu-1})^{-2}\bigr)={\bar\ell}.
\label{5-2}
\end{equation}
Then the  contribution of such element to remainder estimate does not exceed 
$C\mu_\eff^{-1}h_\eff^{-1}=C\mu ^{-1}h^{-1}\ell \gamma ^{1-\nu}$.

After summation over partition with fixed  $\gamma$ one gets due to condition (\ref{5-1}) the same expression but with $\ell=\gamma$ i.e.
$C\mu ^{-1}h^{-1}  \gamma ^{2-\nu}$. One must assume that condition (\ref{5-2}) is fulfilled for this top $\ell$, which is equivalent to 
\begin{equation}
\gamma \ge  C_0\max \bigl(h^{\frac 2 3}, \mu^{-2/(2\nu-1)}\bigr)=L.
\label{5-3}
\end{equation}
Then summation over $x_1$-partition satisfying this condition results in 
$\mu^{-1}h^{-1}L^{2-\nu}$ which as one can check easily does not exceed
$C\mu^{-1/\nu}h^{-1}$.

\smallskip
\noindent
(ii) One needs also to consider elements satisfying (\ref{5-3}) but not (\ref{5-2}). Let us  unify such elements, so in fact I redefine $\ell$ setting it  $\ell =\epsilon |V|+{\bar\ell}$. Then on  elements with $\ell\asymp {\bar\ell}$ either $\mu_\eff \asymp 1$ or $h_\eff\asymp 1$. In the former case 
the condition of a potential being disjoint from 0 is not needed 
while in the latter case  element in question is  forbidden as 
$\mu_\eff \ge C_1$ and the condition of a potential being disjoint from 0 is not needed again. The contribution to the remainder estimate of of each of these elements is either $C h_\eff^{-1}$ (as $\mu_\eff \asymp 1$) or $C\mu_\eff^{-s}$ (as $h_\eff\asymp 1$)  and in any case it does not exceed 
$C\mu^{-1}_\eff h_\eff $ again.  Repeating summation procedure of (a) one can see that the total contribution of such elements does not exceed $C\mu^{-1/\nu}h^{-1}$.

So one needs to consider two other types of elements.

\smallskip
\noindent 
(iii) Consider first elements on which 
\begin{equation}
|x_1|\le C \ell, \qquad \ell=\epsilon |W|.
\label{5-4}
\end{equation}
Then rescaling $x\mapsto x/\ell$ and dividing by $\ell$ brings us to the frames of our main theory with $\mu_\eff = \mu \ell^{\nu -1/2}$, $h_\eff=h\ell^{-3/2}$ and one must assume that 
\begin{equation}
\ell \ge C_0\max\bigl( h^{\frac 2 3}, \mu ^{-2/(2\nu -1)}\bigr).
\label{5-5}
\end{equation}
Then contribution of each element to the remainder estimate does not exceed $C\mu_\eff^{-1/\nu}h_\eff^{-1}= C\mu^{-1/\nu}h^{-1}\ell^{(\nu+1)/2}$ and the sum over partition results in $C\mu^{-1/\nu}h^{-1}$. The exceptional elements will be covered in the next paragraph.

\smallskip
\noindent 
(iv) So, we are left with the analysis of the zone $\{|x_1|+|W|\le L\}$ with $L$ defined by (\ref{5-3}). Due to condition (\ref{5-1}) with no loss of the generality this zone could be replaced by 
$\{|x|\le L\}$. Rescaling $x\mapsto x/L$ and dividing by $L$ we arrive to the same situation as we started but with $h_\eff= hL^{-3/2}$, 
$\mu_\eff = \mu L^{\nu -{\frac 1 2}}$ instead of $h,\mu$ i.e. with either $h_\eff\le 1$, $\mu_\eff\asymp 1$ or with $h_\eff\asymp 1$, $\mu_\eff\ge 1$.

In the former case I just refer to the classical theory: contribution of this zone to the remainder estimate does not exceed 
$O\bigl(h_\eff^{-1}=Ch^{-1}L^{3/2}\le Ch^{-1}\mu^{-1/\nu}$. In the latter case contribution of this zone to the whole asymptotics does not exceed $C\mu_\eff^{-s}$ as this zone is ``forbidden''; one can prove it easily going to the axillary space $L^2(\bR_{x_1})$.
\end{proof}

\subsection{Remark about higher dimensions}

3-dimensional case is rather boring because the best remainder estimate if magnetic field does not degenerate is $O(h^{-2})$ and applying the standard rescaling arguments one will arrive to $O\bigl(h^{-2}\int \gamma^{-1}\,dx\bigr)$
and in the generic case $\gamma $ is the distance to the finite set, thus resulting in $O(h^{-2})$ again (as $\mu \le h^{-1})$). As $\mu \ge h^{-1}$ and $|F|\asymp \gamma$ the principal part will be of magnitude $h^{-3}{\bar\gamma}_1^3$ and the remainder estimate of magnitude $h^{-2}{\bar\gamma}_1^2$ with ${\bar\gamma}_1=(\mu h)^{-1}$ as $h^{-1}\le \mu \le Ch^{-2}$. One can expect the similar results in higher dimensions too.

Even-dimensional case (I will treat it in one of forthcoming papers) is much more interesting but still less than the special 2-dimensional case. First of all as it follows from Martinet \cite{Ma} in the generic case as $d\ge 2$ magnetic matrix $F$ never vanishes and thus $\Tr^+g^{-1}F>0$ (where $\Tr^+g^{-1}F$ is the sum of positive eigenvalues of $ig^{-1}F$ and $g$ is a metrics matrix); therefore unless we replace $V$ by $V+\mu h\Tr^+F$ the principal part will be just 0 and the remainder estimate $O(\mu^{-s})$ as $\mu\ge Ch^{-1}$ thus cutting us from the case in which one expects the periodic zone to be important.

Canonical form of $F$ in the generic case is known due to \cite{Ma} and Roussarie; first of all if $f_1$ is the smallest positive eigenvalue then $f_1\asymp x_1$ in the appropriate coordinate system. Further, in dimension 4 $F$ has no more than 2 degenerating eigenvalues while 2 others are disjoint from 0 and thus 4-dimensional case seems to be nothing but the lift-off of 2-dimensional (as $\mu\le Ch^{-1}$). It is not that simple in particular because there appears a very special submanifold of codimension 3,  but likely the above assertion is correct.

In dimension $d\ge 6$ more than 2 eigenvalues can degenerate but the degeneration happens on the manifold of codimension 6 and it seems to be too thin to have any effect. It does not mean that the proofs are going to be easy but that one should not expect any surprises.

\bibliographystyle{amsalpha}

\providecommand{\bysame}{\leavevmode\hbox to3em{\hrulefill}\thinspace}

\vglue .06truein

\hfill\hfill {\sl   \today \/}

\vglue .06truein

\begin{tabular}{rrl}
&{\hskip 220 pt} &Department of Mathematics,\cr
&&University of Toronto,\cr
&&100, St.George Str.,\cr
&&Toronto, Ontario M5S 3G3\cr
&&Canada\cr
&&ivrii@math.toronto.edu\cr
&&Fax: (416)978-4107\cr
\end{tabular}

\end{document}